\documentclass[a4paper,oneside]{article}
\usepackage{amssymb}
\usepackage{amsfonts}
\usepackage{amsmath}

\input{tcilatex}

\begin{document}

\title{The isomorphism type of the centralizer of an element in a Lie group }
\author{Haibao Duan\thanks{%
Supported by 973 Program 2011CB302400.} and Shali Liu \and Institute of
Mathematics, Chinese Academy of Sciences \and dhb@math.ac.cn}
\maketitle

\begin{abstract}
Let $G$ be an $1$--connected simple Lie group, and let $x\in G$ be a group
element. We determine the isomorphism type of the centralizer $C_{x}$ in
term of a minimal geodesic joinning the group unit $e\in G$ to $x$.

This result is applied to classify the isomorphism types of maximal
subgroups of maximal rank of $G$ \cite{[BS]}, and the isomorphism types of
parabolic subgroups of $G$.

\begin{description}
\item \textsl{2000 Mathematical Subject Classification: }22E15; 53C35

\item \textsl{Key words and phrases:} Lie groups; centralizer; homogeneous
spaces
\end{description}
\end{abstract}

\section{Introduction}

Let $G$ be a compact connected semisimple Lie group with a given element $%
x\in G$. The \textsl{centralizer} $C_{x}$ of $x$ and the \textsl{adjoint
orbit} $M_{x}$ through $x$ are the subspaces of $G$

\begin{quote}
$C_{x}=\{g\in G\mid gx=xg\}$,$\quad M_{x}=\{gxg^{-1}\in G\mid g\in G\}$,
\end{quote}

\noindent respectively. The map $G\rightarrow G$ by $g\rightarrow gxg^{-1}$
is constant along the left cosets of $C_{x}$ in $G$, and induces a
diffeomorphism from the \textsl{homogeneous space} $G/C_{x}$ onto the orbit
space $M_{x}$

\begin{quote}
$f_{x}:G/C_{x}\overset{\cong }{\rightarrow }M_{x}$ by $[g]\rightarrow
gxg^{-1}$.
\end{quote}

\noindent The manifolds $G/C_{x}$ arising in this fashion have offered many
important subjects in geometry, as shown in the next two examples.

\bigskip

\noindent \textbf{Example 1.1. }Let $\mathcal{Z}(G)$ be the center of $G$
and let $x\in G$ be an element with $x^{r}\in \mathcal{Z}(G)$ for some power 
$r\geq 2$. The homogeneous space $G/C_{x}$ possesses a canonical periodic $r$
automorphism

\begin{quote}
$\sigma _{x}:G/C_{x}\rightarrow G/C_{x}$ by $\sigma _{x}[g]=[xgx^{-1}]$.
\end{quote}

\noindent For this reason the pair $(G/C_{x},\sigma _{x})$ is called\textsl{%
\ }an $r$\textsl{--symmetric space} \textsl{of }$G$ in \cite{[DL]}. In the
special case of $r=2$, they are the global Riemannian symmetric spaces of $G$
in the sense of E. Cartan \cite{[He]}.$\square $

\bigskip

\noindent \textbf{Example 1.2. }If the minimal geodesic joining the group
unit $e\in G$ to $x$ is unique, the centralizer $C_{x}$ is a \textsl{%
parabolic subgroup} of $G$. The corresponding homogeneous space $G/C_{x}$ is
called a \textsl{flag manifold} of $G$, which is the focus of the classical
Schubert calculus \cite{[BGG],[DZ],[FP],[K],[LG]}.$\square $

\bigskip

To investigate the geometry and topology of the homogeneous space $G/C_{x}$
it is often necessary to determine explicitly the isomorphism type of the
centralizer $C_{x}$ in term of $x\in G$. However, in the existing
literatures one merely finds method to decide its local type in some special
cases \cite{[He],[BS],[R]}, see Remark 2.10. The purpose of this paper is to
give an explicit procedure for calculating the centralizer $C_{x}$ in term
of a minimal geodesic joining the unit $e$ to $x$ , see Theorem 4.3 in \S %
4.1. This result is applied to classify the isomorphism types of maximal
subgroups of maximal rank of $G$ in \S 4.2, and of parabolic subgroups of $G$
in \S 4.3.

To be precise some notation are needed. For a compact connected Lie group $K$
the identity component of the center $\mathcal{Z}(K)$ of $K$ will be denoted
by $K^{Rad}$, and will be called \textsl{the radical} \textsl{part} of $K$
(this is always a connected torus subgroup of $K$). According to Cartan's
classification on compact Lie groups, up to isomorphism, the group $K$
admits a canonical presentation of the form

\begin{enumerate}
\item[(1.1)] $K\cong (G_{1}\times \cdots \times G_{k}\times K^{Rad})/H$
\end{enumerate}

\noindent in which

i) each $G_{t}$ is one of the $1$--connected \textsl{simple Lie groups}, $%
1\leq t$ $\leq k$;

ii) the denominator $H$ is a finite subgroup of $\mathcal{Z}(G_{1})\times
\cdots \times \mathcal{Z}(G_{k})\times K^{Rad}$.

\noindent It is also known that all $1$--connected simple Lie groups $G$,
together with their centers, are classified by the types $\Phi _{G}$ of
their corresponding root systems tabulated below \cite[p.57]{[Hu]}

\begin{center}
{\footnotesize 
\begin{tabular}{l|l|l|l|l|l|l|l|l|l}
\hline\hline
${\small G}$ & ${\small SU(n)}$ & ${\small Sp(n)}$ & ${\small Spin(2n+1)}$ & 
${\small Spin(2n)}$ & ${\small G}_{2}$ & ${\small F}_{4}$ & ${\small E}_{6}$
& ${\small E}_{7}$ & ${\small E}_{8}$ \\ \hline
${\small \Phi }_{{\small G}}$ & $A_{n-1}$ & $B_{n}$ & $C_{n}$ & $D_{n}$ & $%
{\small G}_{2}$ & ${\small F}_{4}$ & ${\small E}_{6}$ & ${\small E}_{7}$ & $%
{\small E}_{8}$ \\ \hline
$\mathcal{Z}{\small (G)}$ & ${\small Z}_{n}$ & ${\small Z}_{2}$ & ${\small Z}%
_{2}$ & 
\begin{tabular}{l}
${\small Z}_{4}${\small ,\ }${\small n=2k+1}$ \\ 
${\small Z}_{2}{\small \oplus Z}_{2}${\small , }${\small n=2k}$%
\end{tabular}
& ${\small \{e\}}$ & ${\small \{e\}}$ & ${\small Z}_{3}$ & ${\small Z}_{2}$
& ${\small \{e\}}$ \\ \hline\hline
\end{tabular}
}

{\small Table 1. The types and centers of }${\small 1}${\small --connected
simple Lie groups}
\end{center}

\noindent \textbf{Definition 1.3.} In the presentation (1.1) the group $%
G_{1}\times \cdots \times G_{k}$ is called the \textsl{semisimple part} of $%
K $, and is denoted by $K^{s}$.

The obvious quotient (i.e. covering) map $\pi :$ $K^{s}\times
K^{Rad}\rightarrow K$ is called \textsl{the local type} of $K$.$\square $

\bigskip

\noindent \textbf{Corollary 1.4.} \textsl{Let }$K$\textsl{\ be a compact
connected Lie group. The following statements are equivalent.}

\textsl{i) The group }$K$\textsl{\ is semisimple;}

\textsl{ii) }$K^{Rad}=\{e\}$\textsl{;}

\textsl{iii) The local type }$\pi $\textsl{\ of }$K$\textsl{\ agrees with
the universal cover of }$K$\textsl{.}$\square $

\bigskip

The paper is arranged as the following. Section 2 contains a brief
introduction on the roots and weight systems of Lie groups, and obtains the
local type $\pi :$ $C_{x}^{s}\times C_{x}^{Rad}\rightarrow C_{x}$ of a
centralizer $C_{x}$ in term of $x\in G$ in Theorem 2.8. In section 3 we
introduce for each compact connected Lie group $K$ so called \textsl{%
extended weight lattice }$\Pi _{K}^{0}$, together with two\textsl{\
deficiency functions} on it. They are applied in Theorem 3.7 to specify the
isomorphism type of a subgroup $K$ in a semisimple Lie group $G$.
Summarizing results in Theorems 2.8 and 3.7, an explicit procedure for
calculating the isomorphism type of a centralizer $C_{x}$ in an $1$%
--connected Lie group $G$ is given in Theorem 4.3 of Section \S 4.1.
Finally, to demonstrate the use of Theorem 4.3 we determine in \S 4.2 and \S %
4.3 all the centralizers $C_{\exp (u)}$ in an $1$--connected exceptional Lie
group $G$ with $u$ a multiple of a fundamental dominant weight of $G$.

\section{The local type of a centralizer}

In this paper $K$ denotes a compact connected Lie group, and the notion $G$
is reserved for the compact semisimple ones.

Equip the Lie algebra $L(K)$ of $K$ with an inner product $(,)$ so that the
adjoint representation acts as isometries on $L(K)$. Fixing a maximal torus $%
T$ on $K$ the\textsl{\ Cartan subalgebra} of $K$ is the linear subspace $%
L(T) $ of $L(K)$. The dimension $n=\dim T$ is called the \textsl{rank} of
the group $K$.

\subsection{The root system of a compact connected Lie group}

The restriction of the exponential map $\exp :L(K)\rightarrow K$ to the
subspace $L(T)$ defines a set $\mathcal{S}(K)=\{L_{1},\cdots ,L_{m}\}$ of $m=%
\frac{1}{2}(\dim K-n)$ hyperplanes in $L(T)$, namely, the set of\textsl{\
singular hyperplanes }through the origin in $L(T)$ \cite[p.168]{[BD]}. Let $%
l_{k}\subset L(T)$ be the normal line of the plane $L_{k}$ through the
origin, $1\leq k\leq m$. Then the map $\exp $ carries $l_{k}$ onto a circle
subgroup on $T$. Let $\pm \alpha _{k}\in l_{k}$ be the non--zero vectors
with minimal length so that $\exp (\pm \alpha _{k})=e$, $1\leq k\leq m$.

\bigskip

\noindent \textbf{Definition 2.1.} The subset $\Phi _{K}=\{\pm \alpha
_{k}\in L(T)\mid 1\leq k\leq m\}$ of $L(T)$ is called the \textsl{root
system of }$K$.$\square $

\bigskip

\noindent \textbf{Remark 2.2. }We note that\textbf{\ }the root system $\Phi
_{K}$ by Definition 2.1 is \textsl{dual} to those that are commonly used in
literatures, e.g. \cite{[B],[Hu]}. In particular, the symplectic group $%
Sp(n) $ is of the type $B_{n}$, while the spinor group $Spin(2n+1)$ is of
the type $C_{n}$.$\square $

\bigskip

The planes in $\mathcal{S}(K)$ divide $L(T)$ into finitely many convex
regions, called the \textsl{Weyl chambers} of $K$. Fix a regular point $%
x_{0}\in L(T)$, and let $\mathcal{F}(x_{0})$ be the closure of the Weyl
chamber containing $x_{0}$. Assume that $L(x_{0})=\{L_{1},\cdots ,L_{h}\}$
is the subset of $\mathcal{S}(K)$ consisting of the walls of $\mathcal{F}%
(x_{0})$. Then

\begin{enumerate}
\item[(2.1)] $h\leq n$, where the equality holds if and only if $K$ is
semi--simple.
\end{enumerate}

\noindent Let $\alpha _{i}\in \Phi _{K}$ be the root normal to the wall $%
L_{i}\in $ $L(x_{0})$ and pointing toward $x_{0}$.

\bigskip

\noindent \textbf{Definition 2.3.} The subset $S(x_{0})=\{\alpha _{1},\cdots
,\alpha _{h}\}$ of the root system $\Phi _{K}$ is called \textsl{the} 
\textsl{system of} \textsl{simple roots }of $K$ relative to $x_{0}$.

The \textsl{Cartan matrix} of $K$ (relative to $x_{0}$) is the $h\times h$
matrix defined by

\begin{enumerate}
\item[(2.2)] $A=(b_{ij})_{h\times h}$, $b_{ij}=2(a_{i},\alpha _{j})/(\alpha
_{j},\alpha _{j})$.
\end{enumerate}

The lattice in $L(T)$ spanned by all simple roots is called the \textsl{root
lattice}, and is denoted by $\Lambda _{K}^{r}$. The subset of $\Lambda
_{K}^{r}$ consisting of the sums of the simple roots $\alpha _{1},\cdots
,\alpha _{h}$ is denoted by $\Lambda _{K}^{r,+}$. We shall also put

\begin{quote}
$\Phi _{K}^{+}=\Lambda _{K}^{r,+}\cap \Phi _{K}$,
\end{quote}

\noindent whose elements are called the \textsl{positive roots} of $K$.$%
\square $

\bigskip

The set $S(x_{0})$ of simple roots defines a partial order $\prec $ on $L(T)$
by the following rule: for two vectors $u,v\in L(T)$ we say $v\prec u$ if
and only if the difference $u-v$ is a sum of elements in $S(x_{0})$ (i.e.
belongs to $\Lambda _{K}^{r,+}$ \cite[p.47]{[Hu]}).

If $G$ is a simple Lie group, elements in $\Phi _{G}$ has at most two
lengths. Let $\beta \in \Phi _{G}^{+}$ (resp. $\gamma \in $ $\Phi _{G}^{+}$)
be the unique \textsl{maximal short root} (resp. unique \textsl{maximal long
root}) relative to the partial order $\prec $ on $\Phi _{G}^{+}$. From \cite[%
p.66, Table 2]{[Hu]} one gets

\bigskip

\noindent \textbf{Lemma 2.4.} \textsl{Let }$G$\textsl{\ be a simple Lie
group. We have }$\beta =\gamma $\textsl{\ unless }$G=G_{2},F_{4},B_{n}$%
\textsl{\ or }$C_{n}$\textsl{.}

\textsl{Moreover, if }$G=G_{2},F_{4},B_{n}$\textsl{\ or }$C_{n}$\textsl{\ we
have }$\beta \prec \gamma $\textsl{\ and the lengths of the three vectors }$%
\gamma ,\beta ,\delta =\gamma -\beta $\textsl{\ are given in the table below}

\begin{center}
\begin{tabular}{l|lll}
\hline\hline
& $\left\Vert \gamma \right\Vert ^{2}$ & $\left\Vert \beta \right\Vert ^{2}$
& $\left\Vert \delta \right\Vert ^{2}$ \\ \hline
$G_{2}$ & $6$ & $2$ & $2$ \\ \hline
$F_{4}$ & $2$ & $1$ & $1$ \\ \hline
$B_{n}$ & $2$ & $1$ & $1$ \\ \hline
$C_{n}$ & $4$ & $2$ & $2$ \\ \hline\hline
\end{tabular}%
.$\square $
\end{center}

\subsection{The weight system of a semisimple Lie group}

A nonzero vector $\alpha \in \Lambda _{K}^{r}$ gives rise to a linear map

\begin{enumerate}
\item[(2.3)] $\alpha ^{\ast }:L(T)\rightarrow \mathbb{R}$ by $\alpha ^{\ast
}(x)=2(x,\alpha )/(\alpha ,\alpha )$.
\end{enumerate}

\noindent If $\alpha \in \Phi _{K}$ is a root, the map $\alpha ^{\ast }$ is
called \textsl{the inverse root} of $\alpha $ (\cite[p.67]{[Hu]}).

\bigskip

\noindent \textbf{Definition 2.5.} Assume that $G$ is a semisimple Lie
group. The \textsl{weight lattice of} $G$ is the subset of $L(T)$

\begin{quote}
$\Lambda _{G}=\{x\in L(T)\mid \alpha ^{\ast }(x)\in \mathbb{Z}$ for all $%
\alpha \in \Phi _{G}\}$,
\end{quote}

\noindent whose elements are called \textsl{weights}. Its subset

\begin{enumerate}
\item[(2.4)] $\Omega _{G}=\{\omega _{i}\in L(T)\mid \alpha _{j}^{\ast
}(\omega _{i})=\delta _{i,j},$ $\alpha _{j}\in S(x_{0})\}$
\end{enumerate}

\noindent is called the set of \textsl{fundamental dominant weights} of $G$
relative to $x_{0}$, where $\delta _{i,j}$ is the Kronecker symbol.$\square $

\bigskip

\noindent \textbf{Lemma 2.6. }\textsl{Let\ }$G$ \textsl{be a semisimple Lie
group with} \textsl{Cartan matrix} $A$\textsl{,} \textsl{and let} $\Omega
_{G}=\{\omega _{1},\cdots ,\omega _{n}\}$\textsl{\ be the set of fundamental
dominant weights relative to the regular point }$x_{0}$\textsl{. Then}

\textsl{i) }$\Omega _{G}=\{\omega _{1},\cdots ,\omega _{n}\}$\textsl{\ is a
basis for }$\Lambda _{G}$ \textsl{over} $\mathbb{Z}$\textsl{;}

\textsl{ii) the fundamental dominant weights }$\omega _{1},\cdots ,\omega
_{n}$ \textsl{can be expressed in term of} \textsl{the simple roots }$%
\{\alpha _{1},\cdots ,\alpha _{n}\}$\textsl{\ as}

\begin{enumerate}
\item[(2.5)] $\left( 
\begin{tabular}{l}
$\omega _{1}$ \\ 
$\omega _{2}$ \\ 
$\vdots $ \\ 
$\omega _{n}$%
\end{tabular}%
\right) =A^{-1}\left( 
\begin{tabular}{l}
$\alpha _{1}$ \\ 
$\alpha _{2}$ \\ 
$\vdots $ \\ 
$\alpha _{n}$%
\end{tabular}%
\right) ;$
\end{enumerate}

\textsl{iii) for each }$1\leq i\leq n$ \textsl{the half line }$\{t\omega
_{i}\in L(T)\mid t\in \mathbb{R}^{+}\}$\textsl{\ is the edge of the Weyl
chamber }$\mathcal{F}(x_{0})$\textsl{\ opposite to the wall }$L_{i}$\textsl{.%
}$\square $

\bigskip

\noindent \textbf{Proof.} With the assumption that $G$ is semisimple, we
have that $h=n$ by (2.1), and that the Weyl chamber $\mathcal{F}(x_{0})$ is
a convex cone with vertex $0\in L(T)$.

Properties i) and ii) are known. By (2.4) each weight $\omega _{i}\in \Omega
_{G}$ is perpendicular to all the roots $\alpha _{j}$ (i.e. $\omega _{i}\in
L_{j}$) with $j\neq i$. This shows iii).$\square $

\bigskip

For a simple Lie group $G$ we shall adopt the convention that its
fundamental dominant weights $\omega _{1},\cdots ,\omega _{n}$ are ordered
by the order of their corresponding simple roots pictured as the vertices in
the Dynkin diagram of $G$ \cite[p.58]{[Hu]}. Let $\Pi _{G}\subseteq $ $%
\Omega _{G}=\{\omega _{1},\cdots ,\omega _{n}\}$ be the subset of \textsl{%
minimal weights} with respect to the partial order $\prec $ (see \S 2.1) on $%
\Omega _{G}$. By considering the center $\mathcal{Z}(G)$ as a finite
subgroup of $T$ one has the relation

\begin{enumerate}
\item[(2.6)] $\exp (\Pi _{G}\sqcup \{0\})=\mathcal{Z}(G)$, see \cite[p.72,
Excercise 13]{[Hu]}.
\end{enumerate}

\noindent Explicitly, for each simple Lie group $G$ with type $\Phi _{G}$
the set $\Pi _{G}$ of minimal weights is presented in Table 2 below:

\begin{center}
\begin{tabular}{l|l|l|l|l|l|l}
\hline\hline
$\Phi _{G}$ & $A_{n}$ & $B_{n}$ & $C_{n}$ & $D_{n}$ & $E_{6}$ & $E_{7}$ \\ 
\hline
$\Pi _{G}$ & $\left\{ \omega _{i}\right\} _{1\leq i\leq n}$ & $\left\{
\omega _{n}\right\} $ & $\left\{ \omega _{1}\right\} $ & $\left\{ \omega
_{1},\omega _{n-1},\omega _{n}\right\} $ & $\left\{ \omega _{1},\omega
_{6}\right\} $ & $\left\{ \omega _{7}\right\} $ \\ \hline\hline
\end{tabular}

{\small Table 2.} {\small The set }${\small \Pi }_{G}${\small \ of minimal
weights of a simple Lie group }${\small G}$
\end{center}

\subsection{Computing in the fundamental Weyl cell}

Let $G$ be an $1$--connected simple Lie group with maximal short root $\beta
\in \Phi _{G}^{+}$. The\textsl{\ fundamental Weyl cell} of $G$ is the
simplex in $\mathcal{F}(x_{0})$ defined by

\begin{quote}
$\Delta =\{u\in \mathcal{F}(x_{0})\mid \beta ^{\ast }(u)\leq 1\}$.
\end{quote}

\noindent In view of property iii) of Lemma 2.6, a vector $u\in $ $\Delta $
if and only if there is a subset $I_{u}=\{k_{1},\cdots ,k_{r}\}\subseteq
\{1,\cdots ,n\}$ so that

\begin{enumerate}
\item[(2.7)] $u=\lambda _{k_{1}}\omega _{k_{1}}+\cdots +\lambda
_{k_{r}}\omega _{k_{r}}$ with $\lambda _{k_{s}}>0$ and $\beta (u)\leq 1$.
\end{enumerate}

\noindent Let $\overline{I}_{u}$ be the complement of $I_{u}$ in $\{1,\cdots
,n\}$.

\bigskip

\noindent \textbf{Lemma 2.7.} \textsl{If} $u\in \Delta $ \textsl{is nonzero} 
\textsl{with} $u\notin \Omega _{G}$ \textsl{one has}

\begin{enumerate}
\item[(2.8)] $0\leq \alpha ^{\ast }(u)\leq 1$ \textsl{for any positive root }%
$\alpha \in \Phi _{G}^{+}$
\end{enumerate}

\noindent \textsl{Moreover}

\textsl{i) }$\alpha ^{\ast }(u)=0$\textsl{\ implies that }$\alpha $\textsl{\
is a sum of the simple roots }$\alpha _{i}$\textsl{\ with }$i\in \overline{I}%
_{u}$\textsl{;}

\textsl{ii) }$\alpha ^{\ast }(u)=1$ \textsl{implies that }$\beta ^{\ast
}(u)=1$\textsl{,} \textsl{and that} \textsl{there is an} $k\in \{1,2\}$ 
\textsl{so that} $k\beta -\alpha $ \textsl{is a sum of the simple roots }$%
\alpha _{i}$\textsl{\ with }$i\in \overline{I}_{u}$\textsl{.}

\bigskip

\noindent \textbf{Proof. }Let $d:T\times T\rightarrow \mathbb{R}$\ be the
distance function on $T$ induced from the metric on $L(G)$. Since $G$\ is $1$%
--connected, $u\in \Delta $ implies that $d(e,\exp (u))=\left\Vert
u\right\Vert $, see \cite{[C],[Cr]}. It follows that $\left\Vert
u\right\Vert \leq \left\Vert u-\alpha \right\Vert $ for any $\alpha \in
\Lambda _{G}^{r}$. In particular,

\begin{enumerate}
\item[(2.9)] $\alpha ^{\ast }(u)=\frac{2(\alpha ,u)}{(\alpha ,\alpha )}\leq
1 $ for all nonzero $\alpha \in \Lambda _{G}^{r,+}$.
\end{enumerate}

\noindent On the other side, express $\alpha \in \Lambda _{G}^{r,+}$ in term
of the simple roots as

\begin{quote}
$\alpha =k_{1}\alpha _{1}+\cdots +k_{n}\alpha _{n}$ with $k_{i}\in \mathbb{Z}%
^{+}$.
\end{quote}

\noindent With respect to the expression (2.7) we have

\begin{enumerate}
\item[(2.10)] $\alpha ^{\ast }(u)=\frac{2}{(\alpha ,\alpha )}%
(\tsum\limits_{i\in I_{u}}\lambda _{i}k_{i}(\alpha _{i},\omega _{i}))=\frac{1%
}{(\alpha ,\alpha )}(\tsum\limits_{i\in I_{u}}\lambda _{i}k_{i}(\alpha
_{i},a_{i}))\geq 0$
\end{enumerate}

\noindent Since $\Phi _{G}^{+}\subset \Lambda _{G}^{r,+}$ the relation (2.8)
has been shown by (2.9) and (2.10).

By (2.10) $\alpha ^{\ast }(u)=0$ implies that $k_{i}=0$, $i\in I_{u}$. This
shows i).

Let $\alpha \in \Phi _{G}^{+}$ be with $\alpha ^{\ast }(u)=1$. The proof of
ii) will be divided into two cases, depending on whether $\alpha $ is a
short or a long root.

\noindent \textbf{Case 1.} If $\alpha $ is short we get from $\alpha \prec
\beta $ that

\begin{quote}
$\beta -\alpha =k_{1}\alpha _{1}+\cdots +k_{n}\alpha _{n}$ with $k_{i}\in 
\mathbb{Z}^{+}$.
\end{quote}

\noindent Consequently,

\begin{quote}
$1\geq \beta ^{\ast }(u)=\alpha ^{\ast }(u)+\frac{1}{(\beta ,\beta )}%
(\tsum\limits_{i\in I_{u}}\lambda _{i}k_{i}(\alpha _{i},a_{i}))\geq 1$,
\end{quote}

\noindent where the first inequality $\geq $ comes from (2.8), and the
second follows from

\begin{quote}
$\alpha ^{\ast }(u)=1$ and $\frac{1}{(\beta ,\beta )}(\tsum\limits_{i\in
I_{u}}\lambda _{i}k_{i}(\alpha _{i},a_{i}))\geq 0$.
\end{quote}

\noindent This is possible unless $k_{i}=0$ for all $i\in I_{u}$. This shows
ii) when $\alpha $ is short.

\noindent \textbf{Case 2. }If $\alpha $ is long we get from $\alpha \prec
\gamma $ that

\begin{quote}
$\gamma -\alpha =k_{1}\alpha _{1}+\cdots +k_{n}\alpha _{n}$ with $k_{i}\in 
\mathbb{Z}^{+}$.
\end{quote}

\noindent The relation

\begin{quote}
$1\geq \gamma ^{\ast }(u)=\alpha ^{\ast }(u)+\frac{1}{(\gamma ,\gamma )}%
(\tsum\limits_{i\in I_{u}}\lambda _{i}k_{i}(\alpha _{i},a_{i}))\geq 1$
\end{quote}

\noindent and the assumption $\alpha ^{\ast }(u)=1$ force that

\begin{enumerate}
\item[(2.11)] $\gamma ^{\ast }(u)=1$ and $\gamma -\alpha =\tsum\limits_{i\in 
\overline{I}_{u}}k_{i}\alpha _{i}$, $k_{i}\in \mathbb{Z}^{+}$.
\end{enumerate}

\noindent Write $\gamma =\beta +\delta $ with $\delta =\gamma -\beta $ (note
that $\delta \in \Lambda _{G}^{r,+}$). From $\gamma ^{\ast }(u)=1$ and $%
\left\Vert \beta \right\Vert ^{2}=\left\Vert \delta \right\Vert ^{2}$ by
Lemma 2.4 we get that

\begin{quote}
$\frac{\left\Vert \gamma \right\Vert ^{2}}{\left\Vert \beta \right\Vert ^{2}}%
=\beta ^{\ast }(u)+\delta ^{\ast }(u)$,
\end{quote}

\noindent where $\beta ^{\ast }(u),\delta ^{\ast }(u)\leq 1$ by (2.9). Again
by Lemma 2.4 this is possible if and only if when $G=F_{4}$, $B_{n}$, $C_{n}$
and

\begin{enumerate}
\item[(2.12)] $\gamma ^{\ast }(u)=\beta ^{\ast }(u)=\delta ^{\ast }(u)=1$.
\end{enumerate}

\noindent Assume, apart from (2.7), that

\begin{quote}
$u=\lambda _{1}\omega _{1}+\cdots +\lambda _{n}\omega _{n}$.
\end{quote}

If $G=F_{4}$ the system (2.12) becomes

\begin{quote}
$\qquad 2\lambda _{1}+3\lambda _{2}+2\lambda _{3}+\lambda _{4}=2\lambda
_{1}+4\lambda _{2}+3\lambda _{3}+2\lambda _{4}$

$=2\lambda _{1}+2\lambda _{2}+\lambda _{3}=1$.
\end{quote}

\noindent It implies that $\lambda _{1}=\frac{1}{2}$; $\lambda _{2}=\lambda
_{3}=\lambda _{4}=0$ and consequently, $\overline{I}_{u}=\{2,3,4\}$. The
proof of ii) for $G=F_{4}$ is completed by (2.11) and

\begin{quote}
$\gamma =2\beta -\alpha _{2}-2\alpha _{3}-2\alpha _{4}$ (see \cite[p.66,
Table 2]{[Hu]}).
\end{quote}

If $G=B_{n}$ the system (2.12) gives

\begin{quote}
$\qquad \lambda _{1}+2\lambda _{2}+\cdots +2\lambda _{n-1}+\lambda
_{n}=2\lambda _{1}+\cdots +2\lambda _{n-1}+\lambda _{n}$

$=2\lambda _{2}+\cdots +2\lambda _{n-1}+\lambda _{n}=1$
\end{quote}

\noindent It implies that $\lambda _{1}=0$ and consequently $\overline{I}%
_{u}\supseteq \{1\}$. The proof of ii) for $G=B_{n}$ is completed by (2.11)
and

\begin{quote}
$\gamma =2\beta -\alpha _{1}$ (see \cite[p.66, Table 2]{[Hu]}).
\end{quote}

Finally, if $G=C_{n}$ the system (2.12) turns to be

\begin{quote}
$\lambda _{1}+\cdots +\lambda _{n}=\lambda _{1}+2\lambda _{2}+\cdots
+2\lambda _{n}=\lambda _{1}=1$.
\end{quote}

\noindent It implies that $\lambda _{1}=1$, $\lambda _{2}=\lambda
_{3}=\lambda _{4}=0$ and consequently $u=\omega _{1}$. This contradiction to 
$u\notin \Omega _{G}$ completes the proof of ii) for $G=C_{n}$.$\square $

\subsection{The local type of a centralizer}

Let $G$ be an $1$--connected simple Lie group with maximal short root $\beta
\in \Phi _{G}^{+}$. and Dynkin diagram $\Gamma _{G}$. The extended Dynkin
diagram of $G$ with respect to $-\beta $ is denoted by $\widetilde{\Gamma }%
_{G}$. For a $u\in \Delta $ given as that in (2.7) denote by $T_{u}$ (resp. $%
T_{u}^{\beta }$) for the identity component of the subgroup of $T$:

\begin{quote}
$\tbigcap\limits_{i\in \overline{I}_{u}}\ker [\widetilde{\alpha }%
_{i}:T\rightarrow S^{1}]$ (resp. $\ker [\widetilde{\beta }:T\rightarrow
S^{1}]\tbigcap\limits_{i\in \overline{I}_{u}}\ker [\widetilde{\alpha }%
_{i}:T\rightarrow S^{1}]$),
\end{quote}

\noindent where $S^{1}$ is the circle group $\{\exp (2\pi it)\in \mathbb{C}%
\mid t\in \lbrack 0,1]\}$, and where $\widetilde{\alpha }:T\rightarrow S^{1}$
is the homomorphism whose tangent map at the group unit is the inverse root $%
\alpha ^{\ast }:L(T)\rightarrow \mathbb{R}$ of $\alpha \in \Phi _{G}$ (i.e.
(2.3)). Let $\Gamma _{u}\subseteq \Gamma _{G}$ (resp. $\Gamma _{u}^{\beta
}\subseteq \widetilde{\Gamma }_{G}$) be the subdiagram obtained by deleting
all the vertices $\alpha _{k}$\textsl{\ }with\textsl{\ }$k\in I_{u}$, as
well as the edges adjoining to it, from $\Gamma _{G}$ (resp. $\widetilde{%
\Gamma }_{G}$).

Since each $x\in G$ is conjugate in $G$ to an element of the form $\exp (u)$
with $u\in \Delta $, and since the isomorphism type of a subgroup of $G$
remains invariant under conjugation, the study the isomorphism type of a
centralizer $C_{x}$, $x\in G$, can be reduced to the cases $x=\exp (u)$, $%
u\in \Delta $ . Geometrically, the path $\gamma _{u}(t)=\exp (tu)$, $t\in
\lbrack 0,1]$, is a minimal geodesic on $G$ joining the unit $e$ to $x$.

Let $C_{x}^{s}$ and $C_{x}^{Rad}$\ be respectively the semisimple part and
the radical part of the centralizer $C_{x}$. In view of the fact that the
semisimple part of a group is classified by its Dynkin diagram \cite[p.56]%
{[Hu]}, the next result specifies the local type of the centralizer $C_{\exp
(u)}$ in term of $u\in $ $\Delta $.

\bigskip

\noindent \textbf{Theorem 2.8.} \textsl{Let }$G$\textsl{\ be a compact and }$%
1$\textsl{--connected Lie group and let }$x=\exp (u)\in G$\textsl{\ be with }%
$u\in \Delta $ \textsl{but} $u\notin \Omega _{G}\sqcup \{0\}$\textsl{. Then
the centralizer }$C_{x}$ \textsl{is a compact, connected and proper subgroup
of }$G$\textsl{. }

\textsl{Moreover, }

\begin{quote}
\textsl{i) if }$\beta (u)<1$\textsl{,\ then }$\Gamma _{C_{x}^{s}}=\Gamma
_{u} $\textsl{, }$C_{x}^{Rad}=T_{u}$\textsl{;}

\textsl{ii) if }$\beta (u)=1$\textsl{,\ then }$\Gamma _{C_{x}^{s}}=\Gamma
_{u}^{\beta }$\textsl{, }$C_{x}^{Rad}=T_{u}^{\beta }.$
\end{quote}

\bigskip

\noindent \textbf{Proof. }According to Borel\textbf{\ }\cite[Corollary 3.4,
p.101]{[Bo]} the centralizer $C_{x}$ in an $1$--connected Lie group $G$ is
always connected. Furthermore, with the assumption that $u\notin \Omega
_{G}\sqcup \{0\}$, the group $C_{x}$ must be a proper subgroup of $G$.

To show i) and ii) assume that the Cartan decomposition of $L(G)$ is

\begin{quote}
$L(G)=L(T)\oplus \tbigoplus\limits_{\alpha \in \Phi _{G}^{+}}L_{\alpha }$,
\end{quote}

\noindent where $L_{\alpha }$ is the root space belonging to the root $%
\alpha \in \Phi _{G}^{+}$ (\cite[p.35]{[Hu]}). According to \cite[p.189]%
{[BD]}, for $x=\exp (u)$ with $u\in L(T)$ the Cartan decomposition of the
Lie algebra $L(C_{x})$ is

\begin{enumerate}
\item[(2.13)] $L(C_{x})=L(T)\oplus \tbigoplus\limits_{\alpha \in \Psi
_{u}}L_{\alpha }$, where $\Psi _{u}=\{\alpha \in \Phi _{G}^{+}\mid \alpha
^{\ast }(u)\in \mathbb{Z\}}$.
\end{enumerate}

\noindent In view of (2.13) the set $\Phi _{u}=\{\pm \alpha \mid \alpha \in
\Psi _{u}\}$ can be identified with the root system of the semisimple part $%
C_{x}^{s}$ of $C_{x}$.

If $u\in \Delta $ with $\beta (u)<1$ we have by Lemma 2.7 that $\Psi
_{u}=\{\alpha \in \Phi _{G}^{+}\mid \alpha ^{\ast }(u)=0\mathbb{\}}$ and that

\begin{quote}
a) the set of simple roots $\alpha _{i}$ with $i\in \overline{I}_{u}$ is a
base of $\Phi _{u}$ \cite[p.47]{[Hu]}.
\end{quote}

\noindent Consequently, from the definition of the subgroup $T_{u}$ one gets

\begin{quote}
b) $T_{u}\subseteq $ $C_{x}^{Rad}$.
\end{quote}

\noindent The relation $\Gamma _{C_{x}^{s}}=\Gamma _{u}$ is shown by a). For
the dimension reason $\dim T_{u}+$ rank$C_{x}^{s}=\dim T$ we get from b)
that $T_{u}=$ $C_{x}^{Rad}$. This finishes the proof of i).

Similarly, if $u\in \Delta $ with $\beta (u)=1$ we have by Lemma 2.7 that $%
\Psi _{u}=\{\alpha \in \Phi _{G}^{+}\mid \alpha ^{\ast }(u)=0,1\mathbb{\}}$
and that

\begin{quote}
c) any element in $\Phi _{u}$ is a linear combination of the simple roots $%
\alpha _{i}$ with $i\in \overline{I}_{u}$ and $-\beta $ with coefficients
all nonnegative or nonpositive.
\end{quote}

\noindent As a result, we get from c) and the definition of the subgroup $%
T_{u}^{\beta }\subset T$ that

\begin{quote}
d) $T_{u}^{\beta }\subseteq $ $C_{x}^{Rad}$.
\end{quote}

\noindent With the assumptions $\overline{I}_{u}\neq \emptyset $ and $%
u\notin \Omega _{G}$, we conclude by c) that the set $\{\alpha _{i},-\beta
\mid i\in \overline{I}_{u}\}$ is a base of $\Phi _{u}$ \cite[p.47]{[Hu]} and
therefore, $\Gamma _{C_{x}^{s}}=\Gamma _{u}^{\beta }$. Again, for dimension
reason we get from d) that $T_{u}^{\beta }=C_{x}^{Rad}$. This completes the
proof.$\square $

\bigskip

For an $1$--connected simple Lie group $G$ with rank $n$ assume that the
expression of the maximal short root $\beta \in \Phi _{G}^{+}$ in terms of
the simple roots is

\begin{quote}
$\beta =m_{1}\alpha _{1}+\cdots +m_{n}\alpha _{n}$.
\end{quote}

\noindent The set of vertices of the Weyl cell $\Delta $ is clearly given by

\begin{quote}
$\mathcal{V}_{G}=\{0,X_{i}=\frac{(\beta ,\beta )}{m_{i}(\alpha _{i},\alpha
_{i})}\omega _{i}\in \Delta \mid 1\leq i\leq n\}$.
\end{quote}

\noindent Let us put $\mathcal{F}_{G}=\{u_{i}\in \Delta \mid 1\leq i\leq n\}$
with

\begin{quote}
$u_{i}=\left\{ 
\begin{tabular}{l}
$\frac{1}{2}X_{i}$ if $\alpha _{i}$ is short and $m_{i}=1$; \\ 
$X_{i}$ otherwise.%
\end{tabular}%
\right. $
\end{quote}

According to Theorem 2.8, for two vectors $u,u^{\prime }\in \Delta $ one has 
$L(C_{\exp (u)})\supseteq L(C_{\exp (u^{\prime })})$ if either $%
I_{u}\subseteq I_{u^{\prime }}$ and $\beta (u),\beta (u^{\prime })<1$, or $%
I_{u}=I_{u^{\prime }}$ and $\beta (u^{\prime })\leq \beta (u)=1$. Since a
maximal connected subgroup of maximal rank of $G$ must be the centralizer of
some element in $G$ (\cite[Theorem 5]{[BS]}), we get from Theorem 2.8 the
classical result due to Borel and Siebenthal \cite[\S 7]{[BS]}:

\bigskip

\noindent \textbf{Corollary 2.9. }\textsl{For an }$1$\textsl{--connected
simple Lie group }$G$\textsl{\ with rank }$n$\textsl{,\ the set of
centralizers }$\{C_{\exp (u_{i})}\mid 1\leq i\leq n\}$\textsl{\ contains all
the isomorphism types of maximal subgroups of maximal rank of }$G$\textsl{.}$%
\square $

\bigskip

\noindent \textbf{Remark 2.10. }In the classical paper \cite{[BS]} Borel and
Siebenthal intended to find all maximal subgroups of maximal rank of compact
connected Lie groups. For the $1$--connected simple Lie groups they give the
answers only up to local types. As application of our Theorem 4.3, the
isomorphism types of these groups will be determined in Theorem 4.4, compare
Table 4 in \S 4.2 with the table in \cite[\S 7]{[BS]}.

In \cite{[R]} M. Reeder gives a description of the Lie algebra $L(C_{\exp
(u)})$ under the assumption that $m\cdot u\in \Lambda _{G}^{e}$ for some
multiple $m$. We emphasis that Theorem 2.8 is not obvious, in view of the
crucial use of Lemma 2.7 in its proof.

In a recent Web discussion J. Newman suggested the problem of finding an
algorithm for computing the isomorphism type of the centralizer $C_{H}$ of a
finite subgroup $H$ of a simple Lie group $G$ \cite{[NK]}. If $%
\{u_{1},\cdots ,u_{k}\}\subset \Delta $ is a set of vectors in the cell $%
\Delta $ and if $H$ is the subgroup of $G$ generated by $\exp (u_{i})$, $%
1\leq i\leq k$. then the proof of Theorem 2.8, together with a comment of A.
Knutson in the discussion \cite{[NK]}, implies the next Cartan decomposition
of the Lie algebra of $C_{H}$

\begin{quote}
$L(C_{H})=L(T)\oplus \tbigoplus\limits_{\alpha (u_{1})=\cdots =\alpha
(u_{k})=0}L_{\alpha }\oplus \tbigoplus\limits_{\alpha (u_{1})=\cdots =\alpha
(u_{k})=1}L_{\alpha }$, $\alpha \in \Phi _{G}^{+}$.
\end{quote}

\noindent More precisely, a base for the root system of the group $C_{H}$ is
either

i) $\{\alpha _{i},-\beta \mid i\in \overline{I}_{u_{1}}\cap \cdots \cap 
\overline{I}_{u_{k}}\}$ if $\beta (u_{i})=1$ for all $i$, or

ii) $\{\alpha _{i}\mid i\in \overline{I}_{u_{1}}\cap \cdots \cap \overline{I}%
_{u_{k}}\}$ if $\beta (u_{i})<1$ for some $i$.

\noindent In addition, Theorem 4.3 in \S 4 is applicable to determine the
isomorphism type of the identity component of the group $C_{H}$.$\square $

\section{Deficiency functions and their properties}

According to (1.1) a centralizer $C_{\exp (u)}$ with $u\in \Delta $ admits
the presentation

\begin{enumerate}
\item[(3.1)] $C_{\exp (u)}\cong (G_{1}\times \cdots \times G_{k}\times
C_{\exp (u)}^{Rad})/H$.
\end{enumerate}

\noindent Moreover, its local type $\pi :G_{1}\times \cdots \times
G_{k}\times C_{\exp (u)}^{Rad}\rightarrow $ $C_{\exp (u)}$ can be read from
the expression of $u$ in (2.7) by Theorem 2.8. To complete this work it
remains for us to decide in term of $u$ the finite subgroup $H\subseteq 
\mathcal{Z}(G_{1})\times \cdots \times \mathcal{Z}(G_{k})\times C_{\exp
(u)}^{Rad}$ appearing as the denominator in (3.1). The main idea to do so is
to introduce \textsl{the reduced weight system} of $C_{\exp (u)}$, as well
as two \textsl{deficiency functions} on it, which play the role to clarify
the difference between the group $C_{\exp (u)}$ and its local type $%
G_{1}\times \cdots \times G_{k}\times C_{\exp (u)}^{Rad}$.

\subsection{Deficiency functions for semisimple Lie groups}

We begin by introducing the reduced weight system\textsl{\ }and\textsl{\ }%
the deficiency functions for semisimple Lie groups, and demonstrate their
use in specifying the isomorphism type of such a group.

Assume that $G$ is a semisimple Lie group with local type $\pi :$ $%
G_{1}\times \cdots \times G_{k}\rightarrow G$. Fix a maximal tours $T_{i}$
on each $G_{i}$ and take $T=\pi (T_{1}\times \cdots \times T_{k})$ as the
fixed maximal torus on $G$. Then $L(T)=L(T_{1})\oplus \cdots \oplus L(T_{k})$
and the tangent map of $\pi $ at the group unit induces a partition

\begin{enumerate}
\item[(3.3)] $\Omega _{G}=\Omega _{G_{1}}\sqcup \cdots \sqcup \Omega
_{G_{k}} $
\end{enumerate}

\noindent where $\Omega _{G}$ (resp. $\Omega _{G_{i}}$) is the set of
fundamental dominant weights of $G$ (resp. of $G_{i}$) with respect to a
fixed regular point $(x_{1},\cdots ,x_{k})\in L(T)$, $x_{i}\in L(T_{i})$.

\bigskip

\noindent \textbf{Definition 3.1. }With respect to the partition (3.3) 
\textsl{the reduced weight system} of the semisimple group $G$ is the subset
of its weight lattice $\Lambda _{G}$:

\begin{quote}
$\Pi _{G}^{0}=\{\theta _{1}\oplus \cdots \oplus \theta _{k}\in \Lambda
_{G}\mid \theta _{i}\in \Pi _{G_{i}}\sqcup \{0\}\}$.
\end{quote}

\noindent where $\Pi _{G_{i}}$ is the set of minimal weights of the simple
group $G_{i}$ (see Table 2).

The\textbf{\ }integer valued function $\delta _{G}:\Pi _{G}^{0}\rightarrow 
\mathbb{Z}$ defined by

\begin{quote}
$\delta _{G}(\theta )=$ the order of the element $\exp (\theta )$ in the
group $\mathcal{Z}(G)$, $\theta \in \Pi _{G}^{0}$,
\end{quote}

\noindent is called the \textsl{deficiency function }of\textsl{\ }$G$.$%
\square $

\bigskip

Let $\Lambda _{G}^{e}=\exp ^{-1}(e)$ be the \textsl{unit lattice }of $G$. In
the Euclidean space $L(T)$ one has three lattices $\Lambda _{G}^{r}$, $%
\Lambda _{G}^{e}$ and $\Lambda _{G}$ that are subject to the relations

\begin{enumerate}
\item[(3.4)] $\Lambda _{G}^{r}\subseteq \Lambda _{G}^{e}\subseteq \Lambda
_{G}$.
\end{enumerate}

\noindent Immediate, but useful properties of the function $\delta _{G}$ are

\bigskip

\noindent \textbf{Corollary 3.2. }\textsl{Let }$G$\textsl{\ be a semisimple
Lie group with center }$\mathcal{Z}(G)$\textsl{. Then the exponential map }$%
\exp :L(T)\rightarrow T$\textsl{\ satisfies }$\exp (\Pi _{G}^{0})=\mathcal{Z}%
(G)$\textsl{. Moreover,}

\textsl{i) the value }$\delta _{G}(\theta )$\textsl{\ is the least positive
multiple so that }$\delta _{G}(\theta )\cdot \theta \in \Lambda _{G}^{e}$%
\textsl{;}

\textsl{ii)} \textsl{if }$G=G_{1}\times G_{2}$\textsl{,} $\delta
_{G_{1}\times G_{2}}(\theta _{1}\oplus \theta _{2})=$\textsl{\ l.c.m. }$%
\{\delta _{G_{1}}(\theta _{1}),\delta _{G_{2}}(\theta _{2})\}$\textsl{,}

\noindent \textsl{where} $\theta _{i}\in \Pi _{G_{i}}^{0}$\textsl{, }$i=1,2$%
\textsl{, and where l.c.m. means the least common multiple of the indicated
set of integers.}

\bigskip

\noindent \textbf{Proof.} Property i) comes from the fact that the
exponential map $\exp $ induces an one to one correspondences $\Lambda
_{G}/\Lambda _{G}^{e}\cong \mathcal{Z}(G)$.

The item ii), together with the relation $\exp (\Pi _{G}^{0})=\mathcal{Z}(G)$%
, follows from $\exp (\Pi _{G_{i}}\sqcup \{0\})=$ $\mathcal{Z}(G_{i})$ by
(2.6), and $\mathcal{Z}(G_{1}\times G_{2})=\mathcal{Z}(G_{1})\times \mathcal{%
Z}(G_{2})$.$\square $

\bigskip

\noindent \textbf{Example 3.3.} Let $G$ be an $1$--connected semisimple Lie
group with reduced weight system $\Pi _{G}^{0}$. Properties i) and ii) of
Corollary 3.2 is sufficient to evaluate the function $\delta _{G}:\Pi
_{G}^{0}\rightarrow \mathbb{Z}$.

a) If $G$ is simple with Cartan matrix $A$, then the fundamental dominant
weights $\omega _{1},\cdots ,\omega _{n}$ can be expressed by the simple
roots as (by (2.5))

\begin{enumerate}
\item[(3.5)] $\omega _{i}=r_{i,1}\alpha _{1}+\cdots +r_{i,n}\alpha _{n}$
with $r_{i,k}\in \mathbb{Q}$, $A^{-1}=(r_{i,j})_{n\times n}$.
\end{enumerate}

\noindent Since $\Lambda _{G}^{r}=\Lambda _{G}^{e}$ for the $1$--connected
Lie group $G$, the value $\delta _{G}(\omega _{i})$ with $\omega _{i}\in $ $%
\Pi _{G}$ is the least positive integer so that $\delta _{G}(\omega
_{i})\cdot r_{i,k}\in \mathbb{Z}$ for all $1\leq k\leq n$.

For all $1$--connected simple Lie groups the expressions (3.5) can be found
in \cite[p.69]{[Hu]}, from which we can read off the function $\delta
_{G}:\Pi _{G}^{0}\rightarrow \mathbb{Z}$ and tabulate its values in the
third column of the table below:

\begin{center}
\begin{tabular}{l|l|l}
\hline\hline
$G$ & ${\small \Pi }_{G}^{{\small 0}}=\Pi _{G}\amalg \{0\}$ & $\delta
_{G}(\theta )$, $\theta \in \Lambda _{G}^{0}$ \\ \hline
$SU(n+1)$ & $\{\omega _{1},\cdots ,\omega _{n}\}\amalg \{0\}$ & ${\small \{}%
\frac{n+1}{(n+1,k)}{\small \}}_{1\leq k\leq n}\amalg \{1\}$ \\ 
$Sp(n)$ & $\{\omega _{n}\}\amalg \{0\}$ & ${\small \{2\}}\amalg \{1\}$ \\ 
$Spin(2n+1)$ & $\{\omega _{1}\}\amalg \{0\}$ & ${\small \{2\}}\amalg \{1\}$
\\ 
$Spin(4n)$ & $\{\omega _{1},\omega _{2n-1},\omega _{2n}\}\amalg \{0\}$ & $%
{\small \{2,2,2\}}\amalg \{1\}$ \\ 
$Spin(4n+2)$ & $\{\omega _{1},\omega _{2n-1},\omega _{2n}\}\amalg \{0\}$ & $%
{\small \{2,4,4\}}\amalg \{1\}$ \\ 
$E_{6}$ & $\{\omega _{1},\omega _{6}\}\amalg \{0\}$ & ${\small \{3,3\}}%
\amalg \{1\}$ \\ 
$E_{7}$ & $\{\omega _{7}\}\amalg \{0\}$ & ${\small \{2\}}\amalg \{1\}$ \\ 
$G_{2},F_{4},E_{8}$ & $\{0\}$ & $\{1\}$ \\ \hline
\end{tabular}

{\small Table 3.} {\small The deficiency function }${\small \delta }_{%
{\small G}}:{\small \Pi }_{{\small G}}^{{\small 0}}\rightarrow {\small Z}$%
{\small \ of }${\small 1}${\small --connected simple Lie groups}.
\end{center}

b) If $G$ is $1$--connected with local type $G=G_{1}\times \cdots \times
G_{k}$, by ii) of Corollary 3.2 the function $\delta _{G}:\Pi
_{G}^{0}\rightarrow \mathbb{Z}$ is evaluated by

\begin{quote}
$\delta _{G}(\theta _{1}\oplus \cdots \oplus \theta _{k})=$ l.c.m. $\{\delta
_{G_{1}}(\theta _{1}),\cdots ,\delta _{G_{k}}(\theta _{k})\}$,
\end{quote}

\noindent where the values $\delta _{G_{j}}(\theta _{i})$ with $\theta
_{i}\in \Pi _{G_{i}}^{0}$ are given in Table 3.$\square $

\bigskip

In general, assume that $G$ is semisimple with local type

\begin{quote}
$\pi :$ $G^{s}=G_{1}\times \cdots \times G_{k}\rightarrow G$.
\end{quote}

\noindent The tangent map of $\pi $ at the group unit of $\widetilde{T}%
=T_{1}\times \cdots \times T_{k}$ induces the canonical identifications $L(%
\widetilde{T})=L(T),$ $\quad \Pi _{G}^{0}=\Pi _{G^{s}}^{0}$. It is in this
sense that the reduced weight system $\Pi _{G}^{0}$ of $G$ possesses two
deficiency functions:

\begin{quote}
$\delta _{G}$, $\widetilde{\delta }_{G}:\Pi _{G}^{0}\rightarrow \mathbb{Z}$,
\end{quote}

\noindent where $\widetilde{\delta }_{G}=:\delta _{G^{s}}$. The next result
tells that, comparison between these two functions enables one to specify $%
\ker \pi $, which determines the isomorphism type of $G$. For a finite group 
$H$ write $H^{+}$ for the set of all nontrivial elements in $H$.

\bigskip

\noindent \textbf{Lemma 3.4.} \textsl{Let }$G$\textsl{\ be semisimple with
local type }$\pi :$\textsl{\ }$G^{s}\rightarrow G$\textsl{, and let} $%
\widetilde{\exp }:L(\widetilde{T})\rightarrow \widetilde{T}$ \textsl{be the
exponential map of the maximal torus }$\widetilde{T}$\textsl{\ on} $G^{s}$.%
\textsl{\ Then}

\begin{quote}
$\ker \pi ^{+}=\{\widetilde{\exp }(\theta )\in G^{s}\mid \widetilde{\delta }%
_{G}(\theta )>\delta _{G}(\theta )=1,\theta \in \Pi _{G}^{0}\}$\textsl{,}
\end{quote}

\noindent \textsl{where the order of an element }$\widetilde{\exp }(\theta
)\in \ker \pi ^{+}$ \textsl{is }$\widetilde{\delta }_{G}(\theta )$\textsl{, }%
$\theta \in \Pi _{G}^{0}$\textsl{.}

\bigskip

\noindent \textbf{Proof. }Since $\ker \pi \subseteq \mathcal{Z}(G^{s})=%
\widetilde{\exp }(\Pi _{G^{s}}^{0})$ by Corollary 3.2, any element of $\ker
\pi ^{+}$ has the form $\widetilde{\exp }(\theta )$ for some $\theta \in \Pi
_{G}^{0}$.

On the other hand, for an element $\theta \in \Pi _{G}^{0}$ the statements $%
\widetilde{\delta }_{G}(\theta )>1$ and $\delta _{G}(\theta )=1$ are clearly
equivalent to $\widetilde{\exp }(\theta )\in \ker \pi ^{+}$.$\square $

\subsection{Deficiency functions for connected Lie groups}

We need to extend the deficiency functions to all compact connected Lie
groups, so that an analogue of Lemma 3.4 holds for such a group. Assume
therefore that $K$ is a compact connected Lie group with local type

\begin{quote}
$\pi :K^{s}\times K^{Rad}\rightarrow K$, $K^{s}=G_{1}\times \cdots \times
G_{k}$
\end{quote}

\noindent Taking a maximal torus $T_{i}$ on each factor group $G_{i}$ the
quotient homomorphism $\pi $ then carries the maximal torus $\widetilde{T}%
=T_{1}\times \cdots \times T_{k}\times K^{Rad}$ on $K^{s}\times K^{Rad}$
onto the maximal torus $T=\pi (\widetilde{T})$ on $K$. Moreover, the tangent
map of $\pi $ at the group unit induces an identification

\begin{quote}
$L(\widetilde{T})=L(T)=L(T_{1})\oplus \cdots \oplus L(T_{k})\tbigoplus
L(K^{Rad})$
\end{quote}

\noindent where $L(K^{Rad})$ is the Lie algebra of the radical part $K^{Rad}$%
. Let $\Pi _{G_{1}\times \cdots \times G_{k}}^{0}\subset L(T)$ be the
reduced weight system of the semisimple part $G_{1}\times \cdots \times
G_{k} $.

The unit lattice $\Lambda _{K^{Rad}}^{e}\subset L(K^{Rad})$ of the radical
part $K^{Rad}$ is obviously a subset of $\Lambda _{K}^{e}$. Let $\Lambda
_{K^{Rad}}^{e}(\mathbb{Q)}$ be the vector space spanned by the elements in $%
\Lambda _{K^{Rad}}^{e}$ over the rationals $\mathbb{Q}$, regarded as a
subset of $L(T)$.

\bigskip

\noindent \textbf{Definition 3.5. }The\textsl{\ reduced weight lattice} of $%
K $ is the subset of $L(T)$

\begin{quote}
$\Pi _{K}^{0}=\{\omega \oplus \gamma \in L(T)\mid \omega \in \Pi
_{G_{1}\times \cdots \times G_{k}}^{0},\gamma \in \Lambda _{K^{Rad}}^{e}(%
\mathbb{Q)}\}$.
\end{quote}

\noindent The \textsl{deficiency function} of $K$ is the integer valued map $%
\delta _{K}:\Pi _{K}^{0}\rightarrow \mathbb{Z}$ defined by

\begin{quote}
$\delta _{K}(\omega \oplus \gamma )=$ the least multiple so that $\delta
_{K}(\omega \oplus \gamma )\cdot (\omega \oplus \gamma )\in \Lambda _{K}^{e}$%
.$\square $
\end{quote}

Again, the tangent map of $\pi $ at the group unit induces a canonical
identification $\Pi _{K^{s}\times K^{Rad}}^{0}=\Pi _{K}^{0}$ and therefore,
in analogue to the semisimple cases, the set $\Pi _{K}^{0}$ possesses two
deficiency functions:

\begin{quote}
$\qquad \widetilde{\delta }_{K}$, $\delta _{K}:\Pi _{K}^{0}\rightarrow 
\mathbb{Z}$, $\widetilde{\delta }_{K}=:\delta _{K^{s}\times K^{Rad}}$.
\end{quote}

\noindent Clearly, $\widetilde{\delta }_{K}:\Pi _{K}^{0}\rightarrow \mathbb{Z%
}$ depends only on the local type of $K$ in the sense that

\begin{quote}
$\qquad \widetilde{\delta }_{K}(\theta _{1}\oplus \cdots \oplus \theta
_{k}\oplus \gamma )=$ l.c.m. $\{\delta _{G_{1}}(\theta _{1}),\cdots ,\delta
_{G_{k}}(\theta _{k}),\delta _{K^{Rad}}(\gamma )\}$,
\end{quote}

\noindent where $\theta _{i}\in \Pi _{G_{i}}^{0}$, $\gamma \in \Lambda
_{K^{Rad}}^{e}(\mathbb{Q)}$. The next result generalizes Lemma 3.4 from the
semisimple Lie groups to all compact connected ones.

\bigskip

\noindent \textbf{Lemma 3.6.} \textsl{If }$K$\textsl{\ is compact connected
Lie group with local type }$\pi :$\textsl{\ }$K^{s}\times K^{Rad}\rightarrow
K$\textsl{. Then}

\begin{quote}
$\ker \pi ^{+}=\{\widetilde{\exp }(\theta )\in K^{s}\times K^{Rad}\mid 
\widetilde{\delta }_{K}(\theta )>\delta _{K}(\theta )=1,\theta \in \Pi
_{K}^{0}\}$\textsl{,}
\end{quote}

\noindent \textsl{where the order of an element }$\widetilde{\exp }(\theta
)\in \ker \pi ^{+}$\textsl{\ is }$\widetilde{\delta }_{K}(\theta )$\textsl{, 
}$\theta \in \Pi _{G}^{0}$\textsl{.}$\square $

\subsection{The isomorphism type of a subgroup}

Let $K$ be a compact, connected subgroup of a semisimple Lie group $G$ with
inclusion $h:K\rightarrow G$ and local type $\pi :K^{s}\times
K^{Rad}\rightarrow K$. Assume that $h$ carries a maximal torus $T^{\prime }$
of $K$ into that $T$ of $G$, and let $h_{\ast }:L(T^{\prime })\rightarrow
L(T)$ be the\textsl{\ }tangent map of $h$\ at the group unit.

Since $h$ is monomorphic, we have $h_{\ast }^{-1}(\Lambda _{G}^{e})=\Lambda
_{K}^{e}$. It follows that the condition $\delta _{K}(\theta )=1$, $\theta
\in \Pi _{K}^{0}$, is equivalent to $h_{\ast }\theta \in \Lambda _{G}^{e}$.
Therefore, one gets from Theorem 3.6 that

\bigskip\ 

\noindent \textbf{Theorem 3.7.} \textsl{Let }$K$\textsl{\ be a compact,
connected subgroup of a semisimple Lie group }$G$\textsl{\ with inclusion }$%
h:K\rightarrow G$\textsl{\ and local type }$\pi :K^{s}\times
K^{Rad}\rightarrow K$\textsl{. Then}

\begin{enumerate}
\item[(3.6)] $\ker \pi ^{+}=\{\widetilde{\exp }(\theta )\in K^{s}\times
K^{Rad}\mid \widetilde{\delta }_{K}(\theta )>1,h_{\ast }(\theta )\in \Lambda
_{G}^{e},\theta \in \Pi _{K}^{0}\}$.$\square $
\end{enumerate}

\section{The isomorphism types of centralizers in an 1--connected Lie group}

To avoid case by case discussion we assume in this section that $G$ is an $1$%
--connected simple Lie group. Summarizing results in sections \S 2 and \S 3,
our main result is presented in Theorem 4.3, which gives an explicit
procedure for calculating the isomorphism type of a centralizer $C_{\exp
(u)} $.

As applications we determine in \S 4.2 and \S 4.3 the isomorphism type of
those centralizers $C_{\exp (u)}$ with $u\in \Delta $ a multiple of a
fundamental dominant weight.

\subsection{A procedure for calculating a centralizer $C_{x}$}

For a group element $x=\exp (u)\in G$ with $u\in \Delta $ given as that in
(2.7), Theorem 2.8 specifies the local type $\pi $ of the centralizer $C_{x}$%
, hence the deficiency function $\widetilde{\delta }_{C_{x}}:\Pi
_{C_{x}}^{0}\rightarrow \mathbb{Z}$, see discussion in Section \S 3.2. In
order to apply the formula (3.6) to compute the group $\ker \pi $ we need to
know the expressions of $h_{\ast }(\theta )$, $\theta \in \Pi _{C_{x}}^{0}$,
in term of simple roots (or equivalently, the fundamental dominant weights)
of the group $G$, where $h:$ $C_{\exp (u)}\rightarrow G$ is the inclusion.

Let $A$ (resp. $\widetilde{A}$) be the Cartan matrix (resp. the extended
Cartan matrix) of the group $G$ with respect to the system $\{\alpha
_{1},\cdots ,\alpha _{n}\}$ of simple roots (resp. the extended system $%
\{\alpha _{1},\cdots ,\alpha _{n},-\beta \}$ of simple roots). As in (2.7)
for a vector $u\in \Delta $ assume that $I_{u}=\{k_{1},\cdots ,k_{r}\}$ and
let $\overline{I}_{u}=\{j_{1},\cdots ,j_{n-r}\}$ be the complement of the
ordered sequence $I_{u}$ in $\{1,\cdots ,n\}$. Let $A_{u}$ (resp. $%
\widetilde{A}_{u}$) be the matrix obtained from the matrix $A$ (resp. $%
\widetilde{A}$) by deleting all the $i^{th}$ columns and rows with $i\in
I_{u}$.

Assume that the set $\Omega _{C_{\exp (u)}^{s}}$ of fundamental dominant
weights of the semisimple part $C_{\exp (u)}^{s}$ is $\{\omega _{1}^{\prime
},\cdots ,\omega _{s}^{\prime }\}$, where

\begin{quote}
$s=n-r$ if $\beta (u)<1$, and $n-r+1$ if $\beta (u)=1$
\end{quote}

\noindent As direct consequences of Theorem 2.8 we have

\bigskip

\noindent \textbf{Lemma 4.1. }\textsl{The tangent map }$h_{\ast }$\textsl{\
of }$h$\textsl{\ at the group unit satisfies that}

\textsl{i) if }$\beta (u)<1$\textsl{\ then}

\begin{quote}
$\left( 
\begin{array}{c}
h_{\ast }(\omega _{1}^{\prime }) \\ 
\vdots \\ 
h_{\ast }(\omega _{n-r}^{\prime })%
\end{array}%
\right) =A_{u}^{-1}\left( 
\begin{array}{c}
\alpha _{j_{1}} \\ 
\vdots \\ 
\alpha _{j_{n-r}}%
\end{array}%
\right) $
\end{quote}

\noindent \textsl{and}

\begin{quote}
$i_{\ast }(\Lambda _{C_{\exp (u)}^{Rad}}^{e}(\mathbb{Q}))=\{a_{1}\omega
_{k_{1}}+\cdots +a_{r}\omega _{k_{r}}\mid a_{i}\in \mathbb{Q}\}$\textsl{;}
\end{quote}

\textsl{ii) if }$\beta (u)=1$\textsl{\ then}

\begin{quote}
$\left( 
\begin{array}{c}
h_{\ast }(\omega _{1}^{\prime }) \\ 
\vdots \\ 
h_{\ast }(\omega _{n-r}^{\prime }) \\ 
h_{\ast }(\omega _{n-r+1}^{\prime })%
\end{array}%
\right) =\widetilde{A}_{u}^{-1}\left( 
\begin{array}{c}
\alpha _{j_{1}} \\ 
\vdots \\ 
\alpha _{j_{n-r}} \\ 
-\beta%
\end{array}%
\right) $
\end{quote}

\noindent \textsl{and}

\begin{quote}
$h_{\ast }(\Lambda _{C_{\exp (u)}^{Rad}}^{e}(\mathbb{Q}))=\{a_{1}\omega
_{k_{1}}+\cdots +a_{r}\omega _{k_{r}}\mid \Sigma a_{i}\beta ^{\ast }(\omega
_{k_{i}})=0$, $a_{i}\in \mathbb{Q}\}$.
\end{quote}

\noindent \textbf{Proof.} The formula of $h_{\ast }(\omega _{i}^{\prime })$
comes from (2.5). The expressions of $h_{\ast }(\Lambda _{C_{\exp
(u)}^{Rad}}^{e}(\mathbb{Q}))$ follows from the relation $C_{\exp
(u)}^{Rad}=T_{u}$ or $T_{u}^{\beta }$ by Theorem 2.8, as well as the
definition of the groups $T_{u}$ and $T_{u}^{\beta }$ in \S 2.4.$\square $

\bigskip

In general a centralizer $C_{\exp (u)}$ may not be semisimple. As a result
its extended weight system $\Pi _{C_{\exp (u)}}^{0}$ might contain the
infinite factor $\Lambda _{C_{x}^{Rad}}^{e}(\mathbb{Q)\subset }\Pi _{C_{\exp
(u)}}^{0}$, see Definition 3.5. This raises the question whether the
deficiency function $\widetilde{\delta }_{C_{x}}:\Pi _{C_{x}}^{0}\rightarrow 
\mathbb{Z}$ can be effectively calculated. The next result allows us to
reduce the determination of $\ker \pi ^{+}$ by dealing with the finite set

\begin{enumerate}
\item[(4.1)] $H_{u}=\{\theta \in \Pi _{C_{\exp (u)}^{s}}^{0}\mid \delta
_{C_{\exp (u)}^{s}}(\theta )>1$, $h_{\ast }(\theta )\in h_{\ast }(\Lambda
_{C_{\exp (u)}^{Rad}}^{e}(\mathbb{Q}))\func{mod}\Lambda _{G}^{r}\}$,
\end{enumerate}

\noindent which, in practice, can be easily decided from the concrete
expressions of $h_{\ast }(\theta )$ \ (with $\theta \in \Pi _{C_{\exp
(u)}^{s}}^{0}$) and $h_{\ast }(\Lambda _{C_{\exp (u)}^{Rad}}^{e}(\mathbb{Q}%
)) $ by Lemma 4.1. Note that with the assumption that $G$ is $1$--connected,
one has $\Lambda _{G}^{r}=\Lambda _{G}^{e}$.

\bigskip

\noindent \textbf{Theorem 4.2.} \textsl{Let }$\pi $\textsl{\ be the local
type of the centralizer }$C_{\exp (u)}$\textsl{. Then}

\begin{enumerate}
\item[(4.2)] $\ker \pi ^{+}=\{\widetilde{\exp }(\theta -\gamma _{\theta
})\mid \theta \in H_{u}\}$\textsl{,}
\end{enumerate}

\noindent \textsl{where} $\gamma _{\theta }\in \Lambda _{C_{\exp
(u)}^{Rad}}^{e}(\mathbb{Q})\QTR{sl}{\ is}$ \textsl{an arbitrary element
satisfying }

\begin{quote}
$h_{\ast }(\theta )\equiv h_{\ast }(\gamma _{\theta })\func{mod}\Lambda
_{G}^{r}$.
\end{quote}

\noindent \textbf{Proof.} The exponential map $\widetilde{\exp }$ of the
local type $C_{\exp (u)}^{s}\times C_{\exp (u)}^{Rad}$ of $C_{\exp (u)}$
will be written as $\exp _{1}\times \exp _{2}$, where $\exp _{1}$ and $\exp
_{2}$ are the exponential maps of the factors $C_{\exp (u)}^{s}$ and $%
C_{\exp (u)}^{Rad}$, respectively. Let

\begin{quote}
$\pi _{2}:C_{\exp (u)}^{Rad}\rightarrow C_{\exp (u)}$
\end{quote}

\noindent be the restriction of $\pi $ on the second factor $C_{\exp
(u)}^{Rad}$. By the definition of the radical part $K^{Rad}$ of a Lie group $%
K$ (see \S 1) the map $\pi _{2}$ (hence the composition $h\circ \pi _{2}$)
is injective.

For a $\theta \in H_{u}$ let $\gamma _{\theta }\in \Lambda _{C_{\exp
(u)}^{Rad}}^{e}(\mathbb{Q})$ be a vector with $h_{\ast }(\theta )\equiv
h_{\ast }(\gamma _{\theta })\func{mod}\Lambda _{G}^{r}$. We get from

\begin{center}
$\widetilde{\delta }_{C_{\exp (u)}}(\theta -\gamma _{\theta })=$ l.c.m. $%
\{\delta _{C_{\exp (u)}^{s}}(\theta ),\delta _{C_{\exp (u)}^{Rad}}(-\gamma
_{\theta })\}\geq \delta _{C_{\exp (u)}^{s}}(\theta )>1$
\end{center}

\noindent and

\begin{quote}
$h_{\ast }(\theta -\gamma _{\theta })\in \Lambda _{G}^{r}$ (since $h_{\ast
}(\theta )\equiv h_{\ast }(\gamma _{\theta })\func{mod}\Lambda _{G}^{r}$)
\end{quote}

\noindent that $\widetilde{\exp }(\theta -\gamma _{\theta })\in \ker \pi
^{+} $ by Theorem 3.7. Furthermore, the element $\widetilde{\exp }(\theta
-\gamma _{\theta })$ is independent of the choice of $\gamma _{\theta }$
since if $\gamma _{\theta }^{\prime }\in \Lambda _{C_{\exp (u)}^{Rad}}^{e}(%
\mathbb{Q})$ is a second one with $h_{\ast }(\theta )\equiv h_{\ast }(\gamma
_{\theta }^{\prime })\func{mod}\Lambda _{G}^{r}$, then the relation

\begin{quote}
$h_{\ast }(\gamma _{\theta })=h_{\ast }(\gamma _{\theta }^{\prime })\func{mod%
}\Lambda _{G}^{r}$
\end{quote}

\noindent and the injectivity of $h\circ \pi _{2}$ imply that

\begin{quote}
$\exp _{2}(\gamma _{\theta })=\exp _{2}(\gamma _{\theta }^{\prime })$ (in $%
C_{\exp (u)}^{Rad}$).
\end{quote}

\noindent Consequently,

\begin{quote}
$\qquad \widetilde{\exp }(\theta -\gamma _{\theta })=\exp _{1}(\theta
)\times \exp _{2}(-\gamma _{\theta })$

$=\exp _{1}(\theta )\times \exp _{2}(-\gamma _{\theta }^{\prime })=%
\widetilde{\exp }(\theta -\gamma _{\theta }^{\prime })$ (in $C_{\exp
(u)}^{s}\times C_{\exp (u)}^{Rad}$).
\end{quote}

Conversely, for any pair $(\theta ,\gamma )\in \Pi _{C_{\exp (u)}^{s}}^{0}$\ 
$\times \Lambda _{C_{\exp (u)}^{Rad}}^{e}(\mathbb{Q})$ ($=\Pi _{C_{\exp
(u)}}^{0}$) with $\widetilde{\exp }(\theta -\gamma )\in \ker \pi ^{+}$,\ we
get from $h_{\ast }(\theta -\gamma )\in \Lambda _{G}^{r}$ that $h_{\ast
}(\theta )\equiv h_{\ast }(\gamma )\func{mod}\Lambda _{G}^{r}$, and from

\begin{quote}
$\widetilde{\delta }_{C_{\exp (u)}}(\theta -\gamma )=$ l.c.m. $\{\delta
_{C_{\exp (u)}^{s}}(\theta ),\delta _{C_{\exp (u)}^{Rad}}(-\gamma )\}>1$
\end{quote}

\noindent and the injectivity of $h\circ \pi _{2}$ that $\delta _{C_{\exp
(u)}^{s}}(\theta )>1$. This completes the proof.$\square $

\bigskip

Summarizing the results in Theorems 2.8, 4.2 and Lemma 4.1, we obtain the
next explicit procedure for calculating $C_{\exp (u)}$ in term of $u\in
\Delta $.

\bigskip

\noindent \textbf{Theorem 4.3. }\textsl{Let\ }$G$ \textsl{be an }$1$\textsl{%
--connected simple Lie group with fundamental Weyl cell }$\Delta $\textsl{.} 
\textsl{For an vector }$u\in \Delta $ \textsl{with} $u\notin \Omega _{G}$ 
\textsl{the isomorphism type of }$C_{\exp (u)}$\textsl{\ can be obtained by
the procedure below:}

\bigskip

\noindent \textbf{Step 1. }\textsl{Apply Theorem 2.8\ to get the local type
of }$C_{\exp (u)}$\textsl{\ in the form }$G_{1}\times \cdots \times
G_{k}\times C_{x}^{Rad}$ \textsl{with each }$G_{i}$\textsl{\ an }$1$\textsl{%
--connected and semisimple Lie group. Accordingly, write the reduced weight
system of the semisimple part }$C_{\exp (u)}^{s}$ \textsl{as}

\begin{quote}
$\qquad \Pi _{C_{\exp (u)}^{s}}^{0}=\{\theta _{1}\oplus \cdots \oplus \theta
_{k}\mid \theta _{i}\in \Pi _{G_{i}}\sqcup \{0\}\}$\textsl{.}
\end{quote}

\noindent \textbf{Step 2. }\textsl{Apply Lemma 4.1 to get the expressions of
the vectors }$h_{\ast }(\theta )$ \textsl{(with} $\theta \in \Pi
_{G_{i}}^{0})$\textsl{\ and the subspace }$h_{\ast }(\Lambda _{C_{\exp
(u)}^{Rad}}^{e}(\mathbb{Q}))$\textsl{\ in }$\Lambda _{G}(\mathbb{Q})$\textsl{%
. Accordingly, specify the finite set }

\begin{center}
$H_{u}=\{\theta \in \Pi _{C_{\exp (u)}^{s}}^{0}\mid \delta _{C_{\exp
(u)}^{s}}(\theta )>1$, $h_{\ast }(\theta )\in h_{\ast }(\Lambda _{C_{\exp
(u)}^{Rad}}^{e}(\mathbb{Q}))\func{mod}\Lambda _{G}^{r}\}$\textsl{.}
\end{center}

\noindent \textbf{Step 3. }\textsl{The group }$C_{\exp (u)}$\textsl{\ is
isomorphic to }$G_{1}\times \cdots \times G_{k}\times C_{x}^{Rad}/\ker \pi $%
\textsl{\ with}

\begin{quote}
$\ker \pi ^{+}=\{\widetilde{\exp }(\theta -\gamma _{\theta })\in
C_{x}^{s}\times C_{x}^{Rad}\mid \theta \in H_{u}\}$\textsl{,}
\end{quote}

\noindent \textsl{where} $\gamma _{\theta }\in \Lambda _{C_{\exp
(u)}^{Rad}}^{e}(\mathbb{Q})\QTR{sl}{\ }$\textsl{an element satisfying} $%
h_{\ast }(\theta )\equiv h_{\ast }(\gamma _{\theta })\func{mod}\Lambda
_{G}^{r}$\textsl{.}$\square $

\bigskip

Finally, concerning the structure of $\ker \pi $ as a group, the next
observation from the relation (4.2)\textsl{\ }will be repeatedly used in the
forthcoming calculation:

\begin{enumerate}
\item[(4.3)] \textsl{if the set }$H_{u}$\textsl{\ contains }$p-1$\textsl{\
elements and if }$\theta _{0}\in H_{u}$\textsl{\ is an element with }$\delta
_{C_{\exp (u)}^{s}}(\theta _{0})=p$\textsl{, then }$\ker \pi $\textsl{\ is
the cyclic group }$\mathbb{Z}_{p}$ \textsl{of order }$p$\textsl{\ generated
by} \textsl{the element} $\widetilde{\exp }(\theta _{0}-\gamma _{\theta
_{0}})\in G_{1}\times \cdots \times G_{k}\times C_{x}^{Rad}$.$\square $
\end{enumerate}

\subsection{The maximal subgroups of maximal rank in a Lie group}

Let $G$ be an $1$--connected exceptional Lie group with rank $n$. From the
expression of the maximal short root $\beta $ given in \cite[p.66]{[Hu]},
the set $\mathcal{F}_{G}$ (see \S 2.4)) is determined and presented in the
second column of Table 4 below. Applying Theorem 2.8 one obtains the local
types of the centralizers $C_{\exp (u)}$, $u\in \mathcal{F}_{G}$, that are
presented in the third column of the table.

In view of the explicit presentation of $\mathcal{F}_{G}$ in the second
column of the table we note that elements in $\mathcal{F}_{G}$ are of the
form $u_{i}=\frac{\omega _{i}}{p_{i}}$ with $p_{i}>0$ an integer, $1\leq
i\leq n$. According to Borel and Siebenthal \cite[Theorem 6]{[BS]} we have

\begin{enumerate}
\item[(4.4)] the centralizer $C_{\exp (u_{i})}$ with $u_{i}=\frac{\omega _{i}%
}{p_{i}}\in \mathcal{F}_{G}$ is a maximal subgroup of maximal rank of $G$ if
and only if $p_{i}$ is a prime.
\end{enumerate}

Carrying on discussion in Corollary 2.9 and Remark 2.10 we determine the
isomorphism types of the centralizer $C_{\exp (u_{i})}$ for all $u_{i}\in 
\mathcal{F}_{G}$. In order to make the generators of $\ker \pi $ explicit,
for a product $K=K_{1}\times K_{2}$ of two groups we write $\exp _{1}\times
\exp _{2}$ instead of $\exp $, where $\exp $ (resp. $\exp _{i}$, $i=1,2$) is
the exponential map of the group $K$ (resp. of $K_{i}$, $i=1,2$).

\bigskip

\noindent \textbf{Theorem 4.4.} \textsl{Let }$G$\textsl{\ be an }$1$\textsl{%
--connected exceptional Lie group. The isomorphism types of the centralizers 
}$C_{\exp (u)}$\textsl{\ with }$u\in \mathcal{F}_{G}$ \textsl{are given by
the third and fourth columns of the table below, in which those are of
maximal in }$G$\textsl{\ are specified by (4.4) (compare with the table in 
\cite[\S 7]{[BS]}).}

\begin{center}
{\footnotesize 
\begin{tabular}{l|lll}
\hline\hline
$G$ & $\ \ u\in F_{G}$ & Local type of $C_{\exp (u)}$ & $\ker \pi $\
(generator) \\ \hline
$G_{2}$ & $%
\begin{tabular}{l}
$\frac{\omega _{1}}{2}$ \\ 
$\frac{\omega _{2}}{3}$%
\end{tabular}%
$ & $%
\begin{tabular}{l}
$SU(2)\times SU(2)$ \\ 
$SU(3)$%
\end{tabular}%
$ & $%
\begin{tabular}{l}
$Z_{2}$($\exp _{1}(\omega _{1}^{1})\times \exp _{2}(\omega _{1}^{2})$) \\ 
$Z_{3}$($\exp _{1}(\omega _{1}^{1})$)%
\end{tabular}%
$ \\ \hline
$F_{4}$ & 
\begin{tabular}{l}
$\frac{\omega _{1}}{2}$ \\ 
$\frac{\omega _{2}}{4}$ \\ 
$\frac{\omega _{3}}{3}$ \\ 
$\frac{\omega _{4}}{2}$%
\end{tabular}
& 
\begin{tabular}{l}
$Spin(9)$ \\ 
$SU(2)\times SU(4)$ \\ 
$SU(3)\times SU(3)$ \\ 
$Sp(3)\times SU(2)$%
\end{tabular}
& 
\begin{tabular}{l}
$\{0\}$ \\ 
$Z_{2}$($\exp _{1}(\omega _{1}^{1})\times \exp _{2}(\omega _{2}^{2})$) \\ 
$Z_{3}$($\exp _{1}(\omega _{1}^{1})\times \exp _{2}(\omega _{2}^{2})$) \\ 
$Z_{2}$($\exp _{1}(\omega _{3}^{1})\times \exp _{2}(\omega _{1}^{2})$)%
\end{tabular}
\\ \hline
$E_{6}$ & $%
\begin{tabular}{l}
$\frac{\omega _{1}}{2},\frac{\omega _{6}}{2}$ \\ 
$\frac{\omega _{2}}{2},\frac{\omega _{3}}{2},\frac{\omega _{5}}{2}$ \\ 
$\frac{\omega _{4}}{3}$%
\end{tabular}%
$ & $%
\begin{tabular}{l}
$Spin(10)\times S^{1}$ \\ 
$SU(2)\times SU(6)$ \\ 
$SU(3)\times SU(3)\times SU(3)$%
\end{tabular}%
$ & $%
\begin{tabular}{l}
$Z_{4}$($\exp _{1}(\omega _{5}^{1})\times \exp _{2}(-\frac{9}{4}\omega _{1%
\text{(}6\text{)}})$) \\ 
$Z_{2}$($\exp _{1}(\omega _{1}^{1})\times \exp _{2}(\omega _{3}^{2})$) \\ 
$Z_{3}$($\exp _{1}(\omega _{2}^{1})\times \exp _{2}(\omega _{1}^{2})\times
\exp _{3}(\omega _{1}^{3})$)%
\end{tabular}%
$ \\ \hline
$E_{7}$ & $%
\begin{tabular}{l}
$\frac{\omega _{1}}{2},\frac{\omega _{6}}{2}$ \\ 
$\frac{\omega _{2}}{2}$ \\ 
$\frac{\omega _{3}}{3},\frac{\omega _{5}}{3}$ \\ 
$\frac{\omega _{4}}{4}$ \\ 
$\frac{\omega _{7}}{2}$%
\end{tabular}%
$ & $%
\begin{tabular}{l}
$Spin(12)\times SU(2)$ \\ 
$SU(8)$ \\ 
$SU(3)\times SU(6)$ \\ 
$SU(2)\times SU(4)\times SU(4)$ \\ 
$E_{6}\times S^{1}$%
\end{tabular}%
$ & $%
\begin{tabular}{l}
$Z_{2}$($\exp _{1}(\omega _{5}^{1})\times \exp _{2}(\omega _{1}^{2})$) \\ 
$Z_{2}$($\exp (\omega _{4}^{1})$) \\ 
$Z_{3}$($\exp _{1}(\omega _{1}^{1})\times \exp _{2}(\omega _{4}^{2})$) \\ 
$Z_{4}$($\exp _{1}(\omega _{1}^{1})\times \exp _{2}(\omega _{1}^{2})\times
\exp _{2}(\omega _{3}^{3})$) \\ 
$Z_{3}$($\exp _{1}(\omega _{1}^{1})\times \exp _{2}(-\frac{4}{3}\omega _{7})$%
)%
\end{tabular}%
$ \\ \hline
$E_{8}$ & $%
\begin{tabular}{l}
$\frac{\omega _{1}}{2}$ \\ 
$\frac{\omega _{2}}{3}$ \\ 
$\frac{\omega _{3}}{4}$ \\ 
$\frac{\omega _{4}}{6}$ \\ 
$\frac{\omega _{5}}{5}$ \\ 
$\frac{\omega _{6}}{4}$ \\ 
$\frac{\omega _{7}}{3}$ \\ 
$\frac{\omega _{8}}{2}$%
\end{tabular}%
$ & $%
\begin{tabular}{l}
$Spin(16)$ \\ 
$SU(9)$ \\ 
$SU(8)\times SU(2)$ \\ 
$SU(2)\times SU(3)\times SU(6)$ \\ 
$SU(5)\times SU(5)$ \\ 
$Spin(10)\times SU(4)$ \\ 
$E_{6}\times SU(3)$ \\ 
$E_{7}\times SU(2)$%
\end{tabular}%
$ & $%
\begin{tabular}{l}
$Z_{2}$($\exp (\omega _{7}^{1})$) \\ 
$Z_{3}$($\exp (\omega _{3}^{1})$) \\ 
$Z_{4}$($\exp _{1}(\omega _{2}^{1})\times \exp _{2}(\omega _{1}^{2})$) \\ 
$Z_{6}$($\exp _{1}(\omega _{1}^{1})\times \exp _{2}(\omega _{1}^{2})\times
\exp _{3}(\omega _{5}^{3})$) \\ 
$Z_{5}$($\exp _{1}(\omega _{1}^{1})\times \exp _{2}(\omega _{2}^{2})$) \\ 
$Z_{4}$($\exp _{1}(\omega _{4}^{1})\times \exp _{2}(\omega _{3}^{2})$) \\ 
$Z_{3}$($\exp _{1}(\omega _{1}^{1})\times \exp _{2}(\omega _{2}^{2})$) \\ 
$Z_{2}$($\exp _{1}(\omega _{7}^{1})\times \exp _{2}(\omega _{1}^{2})$)%
\end{tabular}%
$ \\ \hline
\end{tabular}
}

{\small Table 4.} {\small The maximal subgroups of the maximal rank of
exceptional Lie groups}
\end{center}

\noindent \textbf{Proof. }It suffices to establish the results in fifth
columns of the table. These will be done by applying the steps 2 and 3
entailed in Theorem 4.3. The calculations will be divided into five cases in
accordance to $G=G_{2},F_{4},E_{6},E_{7}$ and $E_{8}$.

If $G=E_{6}$, $u=\frac{\omega _{1}}{2},\frac{\omega _{6}}{2}$ (resp. $%
G=E_{7} $, $u=\frac{\omega _{7}}{2}$), the centralizer $C_{\exp (u)}$ are
parabolic. The proofs of these cases will be postponed to the next section.

In the remaining cases the centralizers $C_{\exp (u)}$ are always
semisimple. Therefore, the constraint $h_{\ast }(\theta )\in h_{\ast
}(\Lambda _{C_{\exp (u)}^{Rad}}^{e}(\mathbb{Q}))\func{mod}\Lambda _{G}^{r}$
on the set $H_{u}$ (see (4.1)) is equivalent to $h_{\ast }(\theta )\equiv 0%
\func{mod}\Lambda _{G}^{r}$.

It is more convenient for us to express $h_{\ast }(\theta )$ with $\theta
\in \Pi _{C_{\exp (u)}}$ in term of the weights of $G$ (instead the simple
roots). For the transitions from weights to roots we refer to the table in 
\cite[p.69]{[Hu]}.

\bigskip

\noindent \textbf{Case 1.} $G=G_{2}$.

i) If $u=\frac{\omega _{1}}{2}$, the local type of the centralizer $C_{\exp
(u)}$ is $SU(2)\times SU(2)$. Accordingly, assume that the set of
fundamental dominant weights of $C_{\exp (u)}$ is $\Omega =\{\omega
_{1}^{1}\}\amalg \{\omega _{1}^{2}\}$. Applying Lemma 4.1 we get the
expressions of $h_{\ast }(\omega )$ with $\omega \in \Pi _{SU(2)}\sqcup \Pi
_{SU(2)}$ by fundamental dominant weights $\omega _{1},\omega _{2}$ of $%
G_{2} $:

\begin{quote}
$h_{\ast }(\omega _{1}^{1})=-\frac{\omega _{1}}{2}$;$\qquad h_{\ast }(\omega
_{1}^{2})=\omega _{2}-\frac{3}{2}\omega _{1}$.
\end{quote}

\noindent It follows that the set $H_{u}$ consists of the single element $%
\omega _{1}^{1}\oplus \omega _{1}^{2}$ whose deficiency in the group $%
SU(2)\times SU(2)$ is $2$.

Consequently, $\ker \pi =\mathbb{Z}_{2}$ with generator $\exp _{1}(\omega
_{1}^{1})\times \exp _{2}(\omega _{1}^{2})$.

ii) If $u=\frac{\omega _{2}}{3}$, the local type of the centralizer $C_{\exp
(u)}$ is $SU(3)$. Accordingly, assume that the set of fundamental dominant
weights of $C_{\exp (u)}$ is $\Omega =\{\omega _{1}^{1},\omega _{2}^{1}\}$.
Applying Lemma 4.1 we get the expressions of $h_{\ast }(\omega )$ with $%
\omega \in \Pi _{SU(3)}$ by the fundamental dominant weights $\omega
_{1},\omega _{2}$ of $G_{2}$:

\begin{quote}
$h_{\ast }(\omega _{1}^{1})=-\omega _{2}$; $\qquad h_{\ast }(\omega
_{2}^{1})=\omega _{1}-2\omega _{2}$.
\end{quote}

\noindent It follows that the set $H_{u}$ consists of two elements $\omega
_{1}^{1}$ and $\omega _{2}^{1}$ whose deficiencies in the group $SU(3)$ are
both $3$.

Consequently, $\ker \pi =\mathbb{Z}_{3}$ with generator $\exp _{1}(\omega
_{1}^{1})$ by (4.3).

\bigskip

\noindent \textbf{Case 2.} $G=F_{4}$

i) If $u=\frac{\omega _{1}}{2}$, the local type of the centralizer $C_{\exp
(u)}$ is $Spin(9)$. Accordingly, assume that the set of fundamental dominant
weights of $C_{\exp (u)}$ is $\{\omega _{1}^{1},\omega _{2}^{1},\omega
_{3}^{1},\omega _{4}^{1}\}$. Applying Lemma 4.1 we get the expressions of $%
h_{\ast }(\omega )$ with $\omega \in \Pi _{Spin(9)}=\{\omega _{1}^{1}\}$ by
fundamental dominant weights of $F_{4}$

\begin{quote}
$h_{\ast }(\omega _{1}^{1})=-\frac{\omega _{1}}{2}$.
\end{quote}

\noindent It follows that $H_{u}=\emptyset .$ Consequently, $\ker \pi =\{0\}$%
.

\noindent\ \ ii) If $u=\frac{\omega _{2}}{4}$, the local type of the
centralizer $C_{\exp (u)}$ is $SU(2)\times SU(4)$. Accordingly, assume that
the set of fundamental dominant weights of $C_{\exp (u)}$ is $\{\omega
_{1}^{1}\}\amalg \{\omega _{1}^{2},\omega _{2}^{2},\omega _{3}^{2}\}$.
Applying Lemma 4.1 we get the expressions of $h_{\ast }(\omega )$ with $%
\omega \in \Pi _{SU(2)}\sqcup \Pi _{SU(4)}$ by the fundamental dominant
weights of $F_{4}$

\begin{quote}
$h_{\ast }(\omega _{1}^{1})=\omega _{1}-\frac{1}{2}\omega _{2}$

$h_{\ast }(\omega _{1}^{2})=\omega _{3}-\frac{3}{4}\omega _{2}$

$h_{\ast }(\omega _{2}^{2})=\omega _{4}-\frac{1}{2}\omega _{2}$

$h_{\ast }(\omega _{3}^{2})=-\frac{1}{4}\omega _{2}$.
\end{quote}

\noindent It follows that the set $H_{u}$ consists of the single element $%
\omega _{1}^{1}\oplus \omega _{2}^{2}$ whose deficiency in the group $%
SU(2)\times SU(4)$ is $2$.

Consequently, $\ker \pi =\mathbb{Z}_{2}$ with generator $\exp _{1}(\omega
_{1}^{1})\times \exp _{2}(\omega _{2}^{2})$.

iii) If $u=\frac{\omega _{3}}{3}$, the local type of the centralizer $%
C_{\exp (u)}$ is $SU(3)\times SU(3)$. Accordingly, assume that the set of
fundamental dominant weights of $C_{\exp (u)}$ is $\Omega =\{\omega
_{1}^{1},\omega _{2}^{1}\}\amalg \{\omega _{1}^{2},\omega _{2}^{2}\}$.
Applying Lemma 4.1 we get the expressions of $h_{\ast }(\omega )$ with $%
\omega \in \Pi _{SU(3)}\sqcup \Pi _{SU(3)}$ by the fundamental dominant
weights of $F_{4}$:

\begin{quote}
$h_{\ast }(\omega _{1}^{1})=\omega _{1}-\frac{2}{3}\omega _{3}$

$h_{\ast }(\omega _{2}^{1})=\omega _{2}-\frac{4}{3}\omega _{3}$

$h_{\ast }(\omega _{1}^{2})=\omega _{4}-\frac{2}{3}\omega _{3}$

$h_{\ast }(\omega _{2}^{2})=-\frac{1}{3}\omega _{3}$.
\end{quote}

\noindent It follows that the set $H_{u}$ consists of the two elements $%
\omega _{1}^{1}\oplus \omega _{2}^{2}$ and $\omega _{2}^{1}\oplus \omega
_{1}^{2}$ whose deficiencies in the group $SU(3)\times SU(3)$ are both $3$.

Consequently, $\ker \pi =\mathbb{Z}_{3}$ with generator $\exp _{1}(\omega
_{1}^{1})\times \exp _{2}(\omega _{2}^{2})$ by (4.3).

iv) If $u=\frac{\omega _{4}}{2}$, the local type of the centralizer $C_{\exp
(u)}$ is $Sp(3)\times SU(2)$. Accordingly, assume that the set of
fundamental dominant weights of $C_{\exp (u)}$ is $\Omega =\{\omega
_{1}^{1},\omega _{2}^{1},\omega _{3}^{1}\}\amalg \{\omega _{1}^{2}\}$.
Applying Lemma 4.1 we get the expressions of $h_{\ast }(\omega )$ with $%
\omega \in \Pi _{Sp(3)}\sqcup \Pi _{SU(2)}$ by the fundamental dominant
weights of $F_{4}$:

\begin{quote}
$h_{\ast }(\omega _{3}^{1})=\omega _{3}-\frac{3}{2}\omega _{4}$; $\qquad
h_{\ast }(\omega _{1}^{2})=-\frac{1}{2}\omega _{4}$.
\end{quote}

\noindent It follows that the set $H_{u}$ consists of the single element $%
\omega _{3}^{1}\oplus \omega _{1}^{2}$ whose deficiency in the group $%
Sp(3)\times SU(2)$ is $2$.

Consequently, $\ker \pi =\mathbb{Z}_{2}$ with generator $\exp _{1}(\omega
_{3}^{1})\times \exp _{2}(\omega _{1}^{2})$.

\bigskip

\noindent \textbf{Case 3.} $G=E_{6}$.

i) $u=\frac{\omega _{1}}{2}$ (resp. $\frac{\omega _{6}}{2}$). See in the
proof of Theorem 4.6.

\noindent\ \ \ ii) If $u=\frac{\omega _{2}}{2}$ (resp. $\frac{\omega _{3}}{2}%
,\frac{\omega _{5}}{2}$), the local type of the centralizer $C_{\exp (u)}$
is $SU(2)\times SU(6)$. Accordingly, assume that the set of fundamental
dominant weights of $C_{\exp (u)}$ is $\{\omega _{1}^{1}\}\amalg \{\omega
_{1}^{2},\omega _{2}^{2},\omega _{3}^{2},\omega _{4}^{2},\omega _{5}^{2}\}$.
Applying Lemma 4.1 we get the expressions of $h_{\ast }(\omega )$ with $%
\omega \in \Pi _{SU(2)}\sqcup \Pi _{SU(6)}$ by the weights of $E_{6}$:

\begin{quote}
${\small h}_{\ast }{\small (\omega }_{1}^{1}{\small )=-}\frac{1}{2}{\small %
\omega }_{2}${\small \ (resp. }${\small \omega }_{1}{\small -}\frac{1}{2}%
{\small \omega }_{3}${\small ; }${\small \omega }_{6}{\small -}\frac{1}{2}%
{\small \omega }_{5}${\small )}

${\small h}_{\ast }{\small (\omega }_{1}^{2}{\small )=\omega }_{1}{\small -}%
\frac{1}{2}{\small \omega }_{2}${\small \ (resp. }${\small \omega }_{6}%
{\small -}\frac{1}{2}{\small \omega }_{3}${\small ; }${\small \omega }_{1}%
{\small -}\frac{1}{2}{\small \omega }_{5}${\small )}

${\small h}_{\ast }{\small (\omega }_{2}^{2}{\small )=\omega }_{3}{\small %
-\omega }_{2}${\small \ (resp. }${\small \omega }_{5}{\small -\omega }_{3}$%
{\small ; }${\small \omega }_{3}{\small -\omega }_{5}${\small )}

${\small h}_{\ast }{\small (\omega }_{3}^{2}{\small )=\omega }_{4}{\small -}%
\frac{3}{2}{\small \omega }_{2}${\small \ (resp. }${\small \omega }_{4}-%
\frac{3}{2}{\small \omega }_{3}${\small ; }${\small \omega }_{4}{\small -}%
\frac{3}{2}{\small \omega }_{5}${\small )}

${\small h}_{\ast }{\small (\omega }_{4}^{2}{\small )=\omega }_{5}{\small %
-\omega }_{2}${\small \ (resp. }${\small \omega }_{2}{\small -\omega }_{3}$%
{\small ; }${\small \omega }_{2}{\small -\omega }_{5}${\small )}

${\small h}_{\ast }{\small (\omega }_{5}^{2}{\small )=\omega }_{6}{\small -}%
\frac{1}{2}{\small \omega }_{2}${\small \ (resp. }${\small -}\frac{1}{2}%
{\small \omega }_{3}${\small ; }${\small -}\frac{1}{2}{\small \omega }_{5}$%
{\small )}
\end{quote}

\noindent It follows that the set $H_{u}$ consists of the single element $%
\omega _{1}^{1}\oplus \omega _{3}^{2}$ whose deficiency in the group $%
SU(2)\times SU(6)$ is $2$.

Consequently, $\ker \pi =\mathbb{Z}_{2}$ with generator $\exp _{1}(\omega
_{1}^{1})\times \exp _{2}(\omega _{3}^{2})$.

iii) If $u=\frac{\omega _{4}}{3}$, the local type of the centralizer $%
C_{\exp (u)}$ is $SU(3)\times SU(3)\times SU(3)$. Accordingly, assume that
the set of fundamental dominant weights of $C_{\exp (u)}$ is $\Omega
=\{\omega _{1}^{1},\omega _{2}^{1}\}\amalg \{\omega _{1}^{2},\omega
_{2}^{2}\}\amalg \{\omega _{1}^{3},\omega _{2}^{3}\}$. Applying Lemma 4.1,
we get the expressions of $h_{\ast }(\omega )$ with $\omega \in \Pi
_{SU(3)}\sqcup \Pi _{SU(3)}\sqcup \Pi _{SU(3)}$ by the weights of $E_{6}$:

\begin{quote}
$h_{\ast }(\omega _{1}^{1})=\omega _{1}-\frac{1}{3}\omega _{4}$

$h_{\ast }(\omega _{2}^{1})=\omega _{3}-\frac{2}{3}\omega _{4}$

$h_{\ast }(\omega _{1}^{2})=\omega _{2}-\frac{2}{3}\omega _{4}$

$h_{\ast }(\omega _{2}^{2})=-\frac{1}{3}\omega _{4}$

$h_{\ast }(\omega _{1}^{3})=\omega _{5}-\frac{2}{3}\omega _{4}$

$h_{\ast }(\omega _{2}^{3})=\omega _{6}-\frac{1}{3}\omega _{4}$.
\end{quote}

\noindent It follows that the set $H_{u}$ consists of the two elements $%
\omega _{2}^{1}\oplus \omega _{1}^{2}\oplus \omega _{1}^{3}$ and $\omega
_{1}^{1}\oplus \omega _{2}^{2}\oplus \omega _{2}^{3}$ whose deficiency in
the group $SU(3)\times SU(3)\times SU(3)$ are both $3$.

Consequently, $\ker \pi =\mathbb{Z}_{3}$ with generator $\exp _{1}(\omega
_{2}^{1})\times \exp _{2}(\omega _{1}^{2})\times \exp _{3}(\omega _{1}^{3})$
by (4.3).

\bigskip

\noindent \textbf{Case 4.} $G=E_{7}$

i) If $u=\frac{\omega _{1}}{2}$ (resp.$\frac{\omega _{6}}{2}$), the local
type of the centralizer $C_{\exp (u)}$ is $Spin(12)\times SU(2)$.
Accordingly, assume the set of fundamental dominant weights of $C_{\exp (u)}$
is $\Omega =\{\omega _{1}^{1},\omega _{2}^{1},\omega _{3}^{1},\omega
_{4}^{1},\omega _{5}^{1},\omega _{6}^{1}\}\amalg \{\omega _{1}^{2}\}$.
Applying Lemma 4.1 we get the expressions of $h_{\ast }(\omega )$ with $%
\omega \in \Pi _{Spin(12)}\sqcup \Pi _{SU(2)}$ by the weights of $E_{7}$:

\begin{quote}
${\small h}_{\ast }{\small (\omega }_{1}^{1}{\small )=\omega }_{7}{\small -}%
\frac{1}{2}{\small \omega }_{1}${\small \ (resp. }${\small -}\frac{1}{2}%
{\small \omega }_{6}${\small )}

${\small h}_{\ast }{\small (\omega }_{5}^{1}{\small )=\omega }_{3}{\small -}%
\frac{3}{2}{\small \omega }_{1}${\small \ (resp. }${\small \omega }_{5}%
{\small -}\frac{3}{2}{\small \omega }_{6}${\small )}

${\small h}_{\ast }{\small (\omega }_{6}^{1}{\small )=\omega }_{2}{\small %
-\omega }_{1}${\small \ (resp. }${\small \omega }_{2}{\small -\omega }_{6}$%
{\small )}

${\small h}_{\ast }{\small (\omega }_{1}^{2}{\small )=-}\frac{1}{2}{\small %
\omega }_{1}${\small \ (resp. }${\small \omega }_{7}{\small -}\frac{1}{2}%
{\small \omega }_{6}${\small ).}
\end{quote}

\noindent It follows that the set $H_{u}$ consists of the single element $%
\omega _{5}^{1}\oplus \omega _{1}^{2}$ whose deficiency in the group $%
Spin(12)\times SU(2)$ is $2$.

Consequently, $\ker \pi =\mathbb{Z}_{2}$ with generator $\exp _{1}(\omega
_{5}^{1})\times \exp _{2}(\omega _{1}^{2})$.

\noindent\ \ \ ii) If $u=\frac{\omega _{2}}{2}$, the local type of the
centralizer $C_{\exp (u)}$ is $SU(8)$. Accordingly, assume that the set of
fundamental dominant weights of $C_{\exp (u)}$ is $\Omega =\{\omega
_{1}^{1},\omega _{2}^{1},\omega _{3}^{1},\omega _{4}^{1},\omega
_{5}^{1},\omega _{6}^{1},\omega _{7}^{1}\}$. Applying Lemma 4.1 we get the
expressions of $h_{\ast }(\omega )$ with $\omega \in \Pi _{SU(8)}$ by the
weights of $E_{7}$:

\begin{quote}
$h_{\ast }(\omega _{1}^{1})=-\frac{1}{2}\omega _{2}$

$h_{\ast }(\omega _{2}^{1})=\omega _{1}-\omega _{2}$

$h_{\ast }(\omega _{3}^{1})=\omega _{3}-\frac{3}{2}\omega _{2}$

$h_{\ast }(\omega _{4}^{1})=\omega _{4}-2\omega _{2}$

$h_{\ast }(\omega _{5}^{1})=\omega _{5}-\frac{3}{2}\omega _{2}$

$h_{\ast }(\omega _{6}^{1})=\omega _{6}-\omega _{2}$

$h_{\ast }(\omega _{7}^{1})=\omega _{7}-\frac{1}{2}\omega _{2}$
\end{quote}

\noindent It follows that the set $H_{u}$ consists of the single element $%
\omega _{4}^{1}$ whose deficiency in the group $SU(8)$ is $2$.

Consequently, $\ker \pi =\mathbb{Z}_{2}$ with generator $\exp (\omega
_{4}^{1})$.

iii) If $u=\frac{\omega _{3}}{3}$ (resp. $\frac{\omega _{5}}{3}$), the local
type of the centralizer $C_{\exp (u)}$ is $SU(3)\times SU(6)$. Accordingly,
assume that the set of fundamental dominant weights of $C_{\exp (u)}$ is $%
\Omega =\{\omega _{1}^{1},\omega _{2}^{1}\}\amalg \{\omega _{1}^{2},\omega
_{2}^{2},\omega _{3}^{2},\omega _{4}^{2},\omega _{5}^{2}\}$. Applying Lemma
4.1 we get the expressions of $h_{\ast }(\omega )$ with $\omega \in \Pi
_{SU(3)}\sqcup \Pi _{SU(6)}$ by the weights of $E_{7}$:

\begin{quote}
$h_{\ast }(\omega _{1}^{1})=-\frac{1}{3}\omega _{3}$ (resp. $\omega _{6}-%
\frac{2}{3}\omega _{5}$)

$h_{\ast }(\omega _{2}^{1})=\omega _{1}-\frac{2}{3}\omega _{3}$ (resp. $%
\omega _{7}-\frac{1}{3}\omega _{5}$)

$h_{\ast }(\omega _{1}^{2})=\omega _{2}-\frac{2}{3}\omega _{3}$ (resp. $-%
\frac{1}{3}\omega _{5}$)

$h_{\ast }(\omega _{2}^{2})=\omega _{4}-\frac{4}{3}\omega _{3}$ (resp. $%
\omega _{1}-\frac{2}{3}\omega _{5}$)

$h_{\ast }(\omega _{3}^{2})=\omega _{5}-\omega _{3}$ (resp. $\omega
_{3}-\omega _{5}$)

$h_{\ast }(\omega _{4}^{2})=\omega _{6}-\frac{2}{3}\omega _{3}$ (resp. $%
\omega _{4}-\frac{4}{3}\omega _{5}$)

$h_{\ast }(\omega _{5}^{2})=\omega _{7}-\frac{1}{3}\omega _{3}$ (resp. $%
\omega _{2}-\frac{2}{3}\omega _{5}$).
\end{quote}

\noindent It follows that the set $H_{u}$ consists of the two elements $%
\omega _{1}^{1}\oplus \omega _{4}^{2}$ and $\omega _{2}^{1}\oplus \omega
_{2}^{2}$ whose deficiencies in the group $SU(3)\times SU(6)$ is $3$.

Consequently, $\ker \pi =\mathbb{Z}_{3}$ with generator $\exp _{1}(\omega
_{1}^{1})\times \exp _{2}(\omega _{4}^{2})$ by (4.3).

iv) If $u=\frac{\omega _{4}}{4}$, the local type of the centralizer $C_{\exp
(u)}$ is $SU(2)\times SU(4)\times SU(4)$. Accordingly, assume that the set
of fundamental dominant weights of $C_{\exp (u)}$ is $\Omega =\{\omega
_{1}^{1}\}\amalg \{\omega _{1}^{2},\omega _{2}^{2},\omega _{3}^{2}\}\amalg
\{\omega _{1}^{3},\omega _{2}^{3},\omega _{3}^{3}\}$. Applying Lemma 4.1 we
get the expressions of $h_{\ast }(\omega )$ with $\omega \in \Pi
_{SU(2)}\sqcup \Pi _{SU(4)}\sqcup \Pi _{SU(4)}$ by the weights of $E_{7}$:

\begin{quote}
$h_{\ast }(\omega _{1}^{1})=\omega _{2}-{\small \frac{1}{2}\omega }_{4}$

$h_{\ast }(\omega _{1}^{2})=-\frac{1}{4}\omega _{4}$

$h_{\ast }(\omega _{2}^{2})=\omega _{1}-\frac{1}{2}\omega _{4}$

$h_{\ast }(\omega _{3}^{2})=\omega _{3}-\frac{3}{4}\omega _{4}$

$h_{\ast }(\omega _{1}^{3})=\omega _{5}-\frac{3}{4}\omega _{4}$

$h_{\ast }(\omega _{2}^{3})=\omega _{6}-\frac{1}{2}\omega _{4}$

$h_{\ast }(\omega _{3}^{3})=\omega _{7}-\frac{1}{4}\omega _{4}$.
\end{quote}

\noindent It follows that the set $H_{u}$ consists of the $3$ elements

\begin{quote}
$\omega _{2}^{2}\oplus \omega _{2}^{3},$ $\omega _{1}^{1}\oplus \omega
_{1}^{2}\oplus \omega _{3}^{3},\omega _{1}^{1}\oplus \omega _{3}^{2}\oplus
\omega _{1}^{3}$
\end{quote}

\noindent whose deficiencies in the group $SU(2)\times SU(4)\times SU(4)$
are $2,4,4$.

Consequently, $\ker \pi =\mathbb{Z}_{4}$ with generator $\exp _{1}(\omega
_{1}^{1})\times \exp _{2}(\omega _{1}^{2})\times \exp _{2}(\omega _{3}^{3})$
by (4.3).

vi) $u=\frac{\omega _{7}}{2}$. See in the proof of Theorem 4.6.

\bigskip

\noindent \textbf{Case 5.} $G=E_{8}$

i) If $u=\frac{\omega _{1}}{2}$, the local type of the centralizer $C_{\exp
(u)}$ is $Spin(16)$. Accordingly, assume that the set of fundamental
dominant weights of $C_{\exp (u)}$ is $\Omega =\{\omega _{1}^{1},\omega
_{2}^{1},\omega _{3}^{1},\omega _{4}^{1},\omega _{5}^{1},\omega
_{6}^{1},\omega _{7}^{1},\omega _{8}^{1}\}$. Applying Lemma 4.1 we get the
expressions of $h_{\ast }(\omega )$ with $\omega \in \Pi
_{Spin(16)}=\{\omega _{1}^{1},\omega _{7}^{1},\omega _{8}^{1}\}$ by the
weights of $E_{8}$

\begin{quote}
$h_{\ast }(\omega _{1}^{1})=-\frac{1}{2}\omega _{1}$

$h_{\ast }(\omega _{7}^{1})=\omega _{3}-2\omega _{1}$

$h_{\ast }(\omega _{8}^{1})=\omega _{2}-\frac{3}{2}\omega _{1}$
\end{quote}

\noindent It follows that the set $H_{u}$ consists of the single element $%
\omega _{7}^{1}$ whose deficiency in the group $Spin(16)$ is $2$.

Consequently, $\ker \pi =\mathbb{Z}_{2}$ with generator $\exp (\omega
_{7}^{1})$.

\noindent\ \ \ ii) If $u=\frac{\omega _{2}}{3}$, the local type of the
centralizer $C_{\exp (u)}$ is $SU(9)$. Accordingly, assume that the set of
fundamental dominant weights of $C_{\exp (u)}$ is $\Omega =\{\omega
_{1}^{1},\omega _{2}^{1},\omega _{3}^{1},\omega _{4}^{1},\omega
_{5}^{1},\omega _{6}^{1},\omega _{7}^{1},\omega _{8}^{1}\}$. Applying Lemma
4.1 we get the expressions of $h_{\ast }(\omega )$ with $\omega \in \Pi
_{SU(9)}$ by the weights of $E_{8}$:

\begin{quote}
$h_{\ast }(\omega _{1}^{1})=\omega _{1}-\frac{2}{3}\omega _{2}$

$h_{\ast }(\omega _{2}^{1})=\omega _{3}-\frac{4}{3}\omega _{2}$

$h_{\ast }(\omega _{3}^{1})=\omega _{4}-2\omega _{2}$

$h_{\ast }(\omega _{4}^{1})=\omega _{5}-\frac{5}{3}\omega _{2}$

$h_{\ast }(\omega _{5}^{1})=\omega _{6}-\frac{4}{3}\omega _{2}$

$h_{\ast }(\omega _{6}^{1})=\omega _{7}-\omega _{2}$

$h_{\ast }(\omega _{7}^{1})=\omega _{8}-\frac{2}{3}\omega _{2}$

$h_{\ast }(\omega _{8}^{1})=-\frac{1}{3}\omega _{2}$.
\end{quote}

\noindent It follows that the set $H_{u}$ consists of the two elements $%
\omega _{3}^{1}$ and $\omega _{6}^{1}$ whose deficiencies in the group $%
SU(9) $ are both $3$.

Consequently, $\ker \pi =\mathbb{Z}_{3}$ with generator $\exp (\omega
_{3}^{1})$ by (4.3).

iii) If $u=\frac{\omega _{3}}{4}$, the local type of the centralizer $%
C_{\exp (u)}$ is $SU(8)\times SU(2)$. Accordingly, assume that the set of
fundamental dominant weights of $C_{\exp (u)}$ is $\Omega =\{\omega
_{1}^{1},\omega _{2}^{1},\omega _{3}^{1},\omega _{4}^{1},\omega
_{5}^{1},\omega _{6}^{1},\omega _{7}^{1}\}\amalg \{\omega _{1}^{2}\}$.
Applying Lemma 4.1 we get the expressions of $h_{\ast }(\omega )$ with $%
\omega \in \Pi _{SU(8)}\sqcup \Pi _{SU(2)}$ by the weights of $E_{8}$:

\begin{quote}
$h_{\ast }(\omega _{1}^{1})=\omega _{2}-\frac{3}{4}\omega _{3}$

$h_{\ast }(\omega _{2}^{1})=\omega _{4}-\frac{3}{2}\omega _{3}$

$h_{\ast }(\omega _{3}^{1})=\omega _{5}-\frac{5}{4}\omega _{3}$

$h_{\ast }(\omega _{4}^{1})=\omega _{6}-\omega _{3}$

$h_{\ast }(\omega _{5}^{1})=\omega _{7}-\frac{3}{4}\omega _{3}$

$h_{\ast }(\omega _{6}^{1})=\omega _{8}-\frac{1}{2}\omega _{3}$

$h_{\ast }(\omega _{7}^{1})=-\frac{1}{4}\omega _{3}$

$h_{\ast }(\omega _{1}^{2})=\omega _{1}-\frac{1}{2}\omega _{3}$
\end{quote}

\noindent It follows that the set $H_{u}$ consists of the $3$ elements $%
\omega _{4}^{1},\omega _{2}^{1}\oplus \omega _{1}^{2},\omega _{6}^{1}\oplus
\omega _{1}^{2}$ in which $\omega _{2}^{1}\oplus \omega _{1}^{2}$ has
deficiency $4$ in the group $SU(8)\times SU(2)$.

Consequently, $\ker \pi =\mathbb{Z}_{4}$ with generator $\exp _{1}(\omega
_{2}^{1})\times \exp _{2}(\omega _{1}^{2})$ by (4.3).

iv) If $u=\frac{\omega _{4}}{6}$, the local type of the centralizer $C_{\exp
(u)}$ is $SU(2)\times SU(3)\times SU(6)$. Accordingly, assume that the set
of fundamental dominant weights of $C_{\exp (u)}$ is $\Omega =\{\omega
_{1}^{1}\}\amalg \{\omega _{1}^{2},\omega _{2}^{2}\}\amalg \{\omega
_{1}^{3},\omega _{2}^{3},\omega _{3}^{3},\omega _{4}^{3},\omega _{5}^{3}\}$.
Applying Lemma 4.1 we get the expressions of $h_{\ast }(\omega )$ with $%
\omega \in \Pi _{SU(2)}\sqcup \Pi _{SU(3)}\sqcup \Pi _{SU(6)}$ by the
weights of $E_{8}$:

\begin{quote}
$h_{\ast }(\omega _{1}^{1})=\omega _{2}-\frac{1}{2}\omega _{4}$

$h_{\ast }(\omega _{1}^{2})=\omega _{1}-\frac{1}{3}\omega _{4}$

$h_{\ast }(\omega _{2}^{2})=\omega _{3}-\frac{2}{3}\omega _{4}$

$h_{\ast }(\omega _{1}^{3})=\omega _{5}-\frac{5}{6}\omega _{4}$

$h_{\ast }(\omega _{2}^{3})=\omega _{6}-\frac{2}{3}\omega _{4}$

$h_{\ast }(\omega _{3}^{3})=\omega _{7}-\frac{1}{2}\omega _{4}$

$h_{\ast }(\omega _{4}^{3})=\omega _{8}-\frac{1}{3}\omega _{4}$

$h_{\ast }(\omega _{5}^{3})=-\frac{1}{6}\omega _{4}$
\end{quote}

\noindent It follows that the set $H_{u}$ consists of the $5$ elements

\begin{quote}
$\omega _{1}^{1}\oplus \omega _{1}^{2}\oplus \omega _{5}^{3},\omega
_{1}^{1}\oplus \omega _{3}^{3},$ $\omega _{1}^{2}\oplus \omega _{2}^{3}$, $%
\omega _{2}^{2}\oplus \omega _{4}^{3}$, $\omega _{1}^{3}+\omega _{5}^{3}$
\end{quote}

\noindent in which $\omega _{1}^{1}+\omega _{1}^{2}+\omega _{5}^{3}$ has
deficiency $6$ in the group $SU(2)\times SU(3)\times SU(6)$.

Consequently, $\ker \pi =\mathbb{Z}_{6}$ with generator $\exp _{1}(\omega
_{1}^{1})\times \exp _{2}(\omega _{1}^{2})\times \exp _{3}(\omega _{5}^{3})$.

v) If $u=\frac{\omega _{5}}{5}$, the local type of the centralizer $C_{\exp
(u)}$ is $SU(5)\times SU(5)$. Accordingly, assume that the set of
fundamental dominant weights of $C_{\exp (u)}$ is $\Omega =\{\omega
_{1}^{1},\omega _{2}^{1},\omega _{3}^{1},\omega _{4}^{1}\}\amalg \{\omega
_{1}^{2},\omega _{2}^{2},\omega _{3}^{2},\omega _{4}^{2}\}$. Applying Lemma
4.1 we get the expressions of $h_{\ast }(\omega )$ with $\omega \in \Omega
=\Pi _{SU(5)}\sqcup \Pi _{SU(5)}$ by the weights of $E_{8}$:

\begin{quote}
$h_{\ast }(\omega _{1}^{1})=\omega _{1}-\frac{2}{5}\omega _{5}$

$h_{\ast }(\omega _{2}^{1})=\omega _{3}-\frac{4}{5}\omega _{5}$

$h_{\ast }(\omega _{3}^{1})=\omega _{4}-\frac{6}{5}\omega _{5}$

$h_{\ast }(\omega _{4}^{1})=\omega _{2}-\frac{3}{5}\omega _{5}$

$h_{\ast }(\omega _{1}^{2})=\omega _{6}-\frac{4}{5}\omega _{5}$

$h_{\ast }(\omega _{2}^{2})=\omega _{7}-\frac{3}{5}\omega _{5}$

$h_{\ast }(\omega _{3}^{2})=\omega _{8}-\frac{2}{5}\omega _{5}$

$h_{\ast }(\omega _{4}^{2})=-\frac{1}{5}\omega _{5}$
\end{quote}

\noindent It follows that the set $H_{u}$ consists of four elements

\begin{quote}
$\omega _{1}^{1}\oplus \omega _{2}^{2},$ $\omega _{2}^{1}\oplus \omega
_{4}^{2}$, $\omega _{3}^{1}\oplus \omega _{1}^{2}$ and $\omega
_{4}^{1}\oplus \omega _{3}^{2}$
\end{quote}

\noindent whose deficiencies in the group $SU(9)$ are all $5$.

Consequently, $\ker \pi =\mathbb{Z}_{5}$ with generator $\exp _{1}(\omega
_{1}^{1})\times \exp _{2}(\omega _{2}^{2})$ by (4.3).

vi) If $u=\frac{\omega _{6}}{4}$, the local type of the centralizer $C_{\exp
(u)}$ is $Spin(10)\times SU(4)$. Accordingly, assume that the set of
fundamental dominant weights of $C_{\exp (u)}$ is $\Omega =\{\omega
_{1}^{1},\omega _{2}^{1},\omega _{3}^{1},\omega _{4}^{1},\omega
_{5}^{1}\}\amalg \{\omega _{1}^{2},\omega _{2}^{2},\omega _{3}^{2}\}$.
Applying Lemma 4.1 we get the expressions of $h_{\ast }(\omega )$ with $%
\omega \in \Pi _{Spin(10)}\sqcup \Pi _{SU(4)}$ by the weights of $E_{8}$:

\begin{quote}
$h_{\ast }(\omega _{1}^{1})=\omega _{1}-\frac{1}{2}\omega _{6}$

$h_{\ast }(\omega _{4}^{1})=\omega _{2}-\frac{3}{4}\omega _{6}$

$h_{\ast }(\omega _{5}^{1})=\omega _{5}-\frac{5}{4}\omega _{6}$

$h_{\ast }(\omega _{1}^{2})=\omega _{7}-\frac{3}{4}\omega _{6}$

$h_{\ast }(\omega _{2}^{2})=\omega _{8}-\frac{1}{2}\omega _{6}$

$h_{\ast }(\omega _{3}^{2})=-\frac{1}{4}\omega _{6}$
\end{quote}

\noindent It follows that the set $H_{u}$ consists of $3$ elements $\omega
_{1}^{1}\oplus \omega _{2}^{2}$, $\omega _{4}^{1}\oplus \omega _{3}^{2}$, $%
\omega _{5}^{1}\oplus \omega _{1}^{2}$ in which $\omega _{4}^{1}+\omega
_{3}^{2}$ has deficiency $4$ in the group $Spin(10)\times SU(4)$.

Consequently, $\ker \pi =\mathbb{Z}_{4}$ with generator $\exp _{1}(\omega
_{4}^{1})\times \exp _{2}(\omega _{3}^{2})$ by (4.3).

vii) If $u=\frac{\omega _{7}}{3}$, the local type of the centralizer $%
C_{\exp (u)}$ is $E_{6}\times SU(3)$. Accordingly, assume that the set of
fundamental dominant weights of $C_{\exp (u)}^{s}$ is $\Omega =\{\omega
_{1}^{1},\omega _{2}^{1},\omega _{3}^{1},\omega _{4}^{1},\omega
_{5}^{1},\omega _{6}^{1}\}\amalg \{\omega _{1}^{2},\omega _{2}^{2}\}$.
Applying Lemma 4.1 we get the expressions of $h_{\ast }(\omega )$ with $%
\omega \in \Omega =\Pi _{E_{6}}\sqcup \Pi _{SU(3)}$ by the weights of $E_{8}$%
:

\begin{quote}
$h_{\ast }(\omega _{1}^{1})=\omega _{1}-\frac{2}{3}\omega _{7}$,

$h_{\ast }(\omega _{6}^{1})=\omega _{6}-\frac{4}{3}\omega _{7}$,

$h_{\ast }(\omega _{1}^{2})=\omega _{8}-\frac{2}{3}\omega _{7}$,

$h_{\ast }(\omega _{2}^{2})=-\frac{1}{3}\omega _{7}$.
\end{quote}

\noindent It follows that the set $H_{u}$ consists of the two elements $%
\omega _{1}^{1}+\omega _{2}^{2}$ and $\omega _{6}^{1}+\omega _{2}^{2}$ whose
deficiencies in the group $E_{6}\times SU(3)$ are both $3$.

Consequently, $\ker \pi =\mathbb{Z}_{3}$ with generator $\exp _{1}(\omega
_{1}^{1})\times \exp _{2}(\omega _{2}^{2})$ by (4.3).

viii) If $u=\frac{\omega _{8}}{2}$, the local type of the centralizer $%
C_{\exp (u)}$ is $E_{7}\times SU(2)$. Accordingly, assume that the set of
fundamental dominant weights of $C_{\exp (u)}$ is $\Omega =\{\omega
_{1}^{1},\omega _{2}^{1},\omega _{3}^{1},\omega _{4}^{1},\omega
_{5}^{1},\omega _{6}^{1},\omega _{7}^{1}\}\amalg \{\omega _{1}^{2}\}$.
Applying Lemma 4.1 we get the expressions of $h_{\ast }(\omega )$ with $%
\omega \in \Omega =\Pi _{E_{7}}\sqcup \Pi _{SU(2)}$ by the weights of $E_{8}$%
:

\begin{quote}
$h_{\ast }(\omega _{7}^{1})=\omega _{7}-\frac{3}{2}\omega _{8}$, $h_{\ast
}(\omega _{1}^{2})=-\frac{1}{2}\omega _{8}$
\end{quote}

\noindent It follows that the set $H_{u}$ consists of a single element $%
\omega _{7}^{1}+\omega _{1}^{2}$ whose deficiency in the group $E_{7}\times
SU(2)$ is $2$.

Consequently, $\ker \pi =\mathbb{Z}_{2}$ with generator $\exp _{1}(\omega
_{7}^{1})\times \exp _{2}(\omega _{1}^{2})$.$\square $

\bigskip

\noindent \textbf{Remark 4.5.} Based on concrete constructions of the $1$%
--connected exceptional Lie groups, Yokota obtained also the isomorphism
types of maximal subgroups of maximal rank in the recent book \cite{[Y]}. In
comparison, our approach is free of the types of the Lie groups.

In our sequel work \cite{[DL]} certain results of Theorem 4.4 are applied to
determine the fixed set of the inverse involution $G\rightarrow G$, $%
g\rightarrow g^{-1}$ on an exceptional simple Lie group $G$.$\square $

\subsection{The isomorphism types of parabolic subgroups}

If $u\in \Delta $ is a vector with $\beta (u)<1$, the centralizer $C_{\exp
(u)}$ is a \textsl{parabolic subgroup} of $G$ whose isomorphism type depends
only on the subset $I_{u}\subseteq \{1,\cdots ,n\}$ by Theorem 2.8. The
corresponding homogenous space $G/C_{\exp (u)}$ is a smooth projective
variety, called a \textsl{flag manifold} of $G$ \cite%
{[BGG],[DZ],[FP],[K],[LG]}.

In the case where the set $I_{u}$ is a singleton and $\beta (u)<1$, the
vector $u$ is an interior point on an edge of the cell $\Delta $ from the
origin $0$. The homogeneous space $G/C_{\exp (u)}$ is also known as \textsl{%
a generalized Grassmannian} of $G$ \cite{[DZ]}.

Theorem 4.3 is ready to apply to determine the isomorphism types of all
parabolic subgroups of in a given $1$--connected Lie group $G$. This is
demonstrated in the proof of the next result.

\bigskip

\noindent \textbf{Theorem 4.6. }\textsl{Let }$G$\textsl{\ be an }$1$\textsl{%
--connected exceptional Lie group. For each }$u\in \Delta $\textsl{\ with }$%
\beta (u)<1$\textsl{\ and with }$I_{u}=\{i\}$\textsl{\ a singleton, the
isomorphism type of the centralizer }$C_{\exp (u)}$\textsl{\ is given in
Table 4 below:}

\begin{center}
{\footnotesize 
\begin{tabular}{l|lll}
\hline\hline
$G$ & $\ \ I_{u}$ & the local type of $C_{\exp (u)}$ & $\ker \pi $\
(generator of $\ker \pi $) \\ \hline
$G_{2}$ & $%
\begin{tabular}{l}
$\{1\}$ \\ 
$\{2\}$%
\end{tabular}%
$ & $%
\begin{tabular}{l}
$SU(2)\times S^{1}$ \\ 
$SU(2)\times S^{1}$%
\end{tabular}%
$ & $%
\begin{tabular}{l}
$Z_{2}$($\exp _{1}(\omega _{1}^{1})\times \exp _{2}(-\frac{1}{2}\omega _{1})$%
) \\ 
$Z_{2}$($\exp _{1}(\omega _{1}^{1})\times \exp _{2}(-\frac{1}{2}\omega _{2})$%
)%
\end{tabular}%
$ \\ \hline
$F_{4}$ & $%
\begin{tabular}{l}
$\{1\}$ \\ 
$\{2\}$ \\ 
$\{3\}$ \\ 
$\{4\}$%
\end{tabular}%
$ & $%
\begin{tabular}{l}
$Spin(7)\times S^{1}$ \\ 
$SU(2)\times SU(3)\times S^{1}$ \\ 
$SU(2)\times SU(3)\times S^{1}$ \\ 
$Sp(3)\times S^{1}$%
\end{tabular}%
$ & $%
\begin{tabular}{l}
$Z_{2}$($\exp _{1}(\omega _{3}^{1})\times \exp _{2}(-\frac{1}{2}\omega _{1})$%
) \\ 
$Z_{6}$($\exp _{1}(\omega _{1}^{1})\times \exp _{2}(\omega _{1}^{2})\times
\exp _{3}(-\frac{5}{6}\omega _{2})$) \\ 
$Z_{6}$($\exp _{1}(\omega _{1}^{1})\times \exp _{2}(\omega _{1}^{2})\times
\exp _{3}(-\frac{5}{6}\omega _{3})$) \\ 
$Z_{2}$($\exp _{1}(\omega _{3}^{1})\times \exp _{2}(-\frac{1}{2}\omega _{4})$%
)%
\end{tabular}%
$ \\ \hline
$E_{6}$ & $%
\begin{tabular}{l}
$\{1\}$ \\ 
$\{2\}$ \\ 
$\{3\}$ \\ 
$\{4\}$ \\ 
$\{5\}$ \\ 
$\{6\}$%
\end{tabular}%
$ & $%
\begin{tabular}{l}
$Spin(10)\times S^{1}$ \\ 
$SU(6)\times S^{1}$ \\ 
$SU(2)\times SU(5)\times S^{1}$ \\ 
$SU(2)\times SU(3)\times SU(3)\times S^{1}$ \\ 
$SU(2)\times SU(5)\times S^{1}$ \\ 
$Spin(10)\times S^{1}$%
\end{tabular}%
$ & $%
\begin{tabular}{l}
$Z_{4}$($\exp _{1}(\omega _{4}^{1})\times \exp _{2}(-\frac{3}{4}\omega _{1})$%
) \\ 
$Z_{2}$($\exp _{1}(\omega _{3}^{1})\times \exp _{2}(-\frac{1}{2}\omega _{2})$%
) \\ 
$Z_{10}$($\exp _{1}(\omega _{1}^{1})\times \exp _{2}(\omega _{1}^{2})\times
\exp _{3}(-\frac{9}{10}\omega _{3})$) \\ 
$Z_{6}$($\exp (\omega _{1}^{1})\times \exp _{2}(\omega _{1}^{2})\times \exp
_{3}(\omega _{2}^{3})\times \exp _{4}(\frac{1}{6}\omega _{4})$) \\ 
$Z_{10}$($\exp _{1}(\omega _{1}^{1})\times \exp _{2}(\omega _{1}^{2})\times
\exp _{3}(-\frac{21}{10}\omega _{5})$) \\ 
$Z_{4}$($\exp _{1}(\omega _{1}^{1})\times \exp _{2}(-\frac{9}{4}\omega _{6})$%
)%
\end{tabular}%
$ \\ \hline
$E_{7}$ & $%
\begin{tabular}{l}
$\{1\}$ \\ 
$\{2\}$ \\ 
$\{3\}$ \\ 
$\{4\}$ \\ 
$\{5\}$ \\ 
$\{6\}$ \\ 
$\{7\}$%
\end{tabular}%
$ & $%
\begin{tabular}{l}
$Spin(12)\times S^{1}$ \\ 
$SU(7)\times S^{1}$ \\ 
$SU(2)\times SU(6)\times S^{1}$ \\ 
$SU(2)\times SU(3)\times SU(4)\times S^{1}$ \\ 
$SU(3)\times SU(5)\times S^{1}$ \\ 
$SU(2)\times Spin(10)\times S^{1}$ \\ 
$E_{6}\times S^{1}$%
\end{tabular}%
$ & $%
\begin{tabular}{l}
$Z_{2}$($\exp _{1}(\omega _{5}^{1})\times \exp _{2}(-\frac{1}{2}\omega _{1})$%
) \\ 
$Z_{7}$($\exp _{1}(\omega _{3}^{1})\times \exp _{2}(-\frac{2}{7}\omega _{2})$%
) \\ 
$Z_{6}$($\exp _{1}(\omega _{1}^{1})\times \exp _{2}(\omega _{4}^{2})\times
\exp _{3}(-\frac{5}{6}\omega _{3})$) \\ 
$Z_{12}$($\exp _{1}(\omega _{1}^{1})\times \exp _{2}(\omega _{1}^{2})\times
\exp _{3}(\omega _{1}^{3})\times \exp _{3}(-\frac{7}{12}\omega _{4})$) \\ 
$Z_{15}$($\exp _{1}(\omega _{2}^{1})\times \exp _{2}(\omega _{2}^{2})\times
\exp _{3}(-\frac{16}{15}\omega _{5})$) \\ 
$Z_{4}$($\exp _{1}(\omega _{1}^{1})\times \exp _{2}(\omega _{5}^{2})\times
\exp _{3}(-\frac{3}{4}\omega _{6})$) \\ 
$Z_{3}$($\exp _{1}(\omega _{1}^{1})\times \exp _{2}(-\frac{4}{3}\omega _{7})$%
)%
\end{tabular}%
$ \\ \hline
$E_{8}$ & $%
\begin{tabular}{l}
$\{1\}$ \\ 
$\{2\}$ \\ 
$\{3\}$ \\ 
$\{4\}$ \\ 
$\{5\}$ \\ 
$\{6\}$ \\ 
$\{7\}$ \\ 
$\{8\}$%
\end{tabular}%
$ & $%
\begin{tabular}{l}
$Spin(14)\times S^{1}$ \\ 
$SU(8)\times S^{1}$ \\ 
$SU(2)\times SU(7)\times S^{1}$ \\ 
$SU(2)\times SU(3)\times SU(5)\times S^{1}$ \\ 
$SU(4)\times SU(5)\times S^{1}$ \\ 
$Spin(10)\times SU(3)\times S^{1}$ \\ 
$E_{6}\times SU(2)\times S^{1}$ \\ 
$E_{7}\times S^{1}$%
\end{tabular}%
$ & $%
\begin{tabular}{l}
$Z_{4}$($\exp _{1}(\omega _{6}^{1})\times \exp _{2}(-\frac{1}{4}\omega _{1})$%
) \\ 
$Z_{8}$($\exp _{1}(\omega _{1}^{1})\times \exp _{2}(-\frac{3}{8}\omega _{2})$%
) \\ 
$Z_{14}$($\exp _{1}(\omega _{1}^{1})\times \exp _{2}(\omega _{1}^{2})\times
\exp _{3}(-\frac{11}{14}\omega _{3})$) \\ 
$Z_{30}$($\exp _{1}(\omega _{1}^{1})\times \exp _{2}(\omega _{1}^{2})\times
\exp _{3}(\omega _{1}^{3})\times \exp _{4}(-\frac{11}{30}\omega _{4})$) \\ 
$Z_{20}$($\exp _{1}(\omega _{1}^{1})\times \exp _{2}(\omega _{1}^{2})\times
\exp _{3}(-\frac{17}{20}\omega _{5})$) \\ 
$Z_{12}$($\exp _{1}(\omega _{4}^{1})\times \exp _{2}(\omega _{1}^{2})\times
\exp _{3}(-\frac{7}{12}\omega _{6})$) \\ 
$Z_{6}$($\exp _{1}(\omega _{1}^{1})\times \exp _{2}(\omega _{1}^{2})\times
\exp _{3}(-\frac{5}{6}\omega _{7})$) \\ 
$Z_{2}$($\exp _{1}(\omega _{7}^{1})\times \exp _{2}(-\frac{1}{2}\omega _{8})$%
)%
\end{tabular}%
$ \\ \hline
\end{tabular}
}

{\small Table 5.} {\small The parabolic subgroup corresponding to a weight
in an }${\small 1}${\small --connected exceptional Lie groups}
\end{center}

\noindent \textbf{Proof. }Again,\textbf{\ }the proof will be divided into
five cases in accordance to $G=G_{2}$, $F_{4}$, $E_{6}$, $E_{7}$ and $E_{8}$.

Since the isomorphism type of the parabolic subgroup $C_{\exp (\lambda
u_{i})}$ with $u_{i}\in \mathcal{F}_{G}$ and $\lambda \in (0,1)$ is
irrelevant with the parameter $\lambda $, we can take $u=\frac{1}{2}u_{i}$
as a representative for the case $I_{u}=\{i\}$. With this convention the
radical part of $C_{\exp (u)}$ is simply the circle subgroup $S^{1}=\{\exp
(t\omega _{i})\in G\mid t\in \mathbb{R}\}$ on $G$ by Theorem 2.8. As a
result we have by Lemma 4.1 that

\begin{quote}
$h_{\ast }(\Lambda _{C_{\exp (u)}^{Rad}}^{e}(\mathbb{Q}))=\{t\omega _{i}\mid
t\in \mathbb{Q}\}$.
\end{quote}

\noindent \textbf{Case 1.} $G=G_{2}$

For $u=\frac{\omega _{1}}{4}$ (resp. $\frac{\omega _{2}}{4}$) the local type
of the centralizer $C_{\exp (u)}$ is $SU(2)\times S^{1}$. Accordingly,
assume that the set of fundamental dominant weights of the semisimple part $%
SU(2)$ is $\Omega =\{\omega _{1}^{1}\}$. Applying Lemma 4.1 we get the
expression of $h_{\ast }(\omega )$ with $\omega \in \Pi _{SU(2)}$ by the
simple roots of $G_{2}$:

\begin{quote}
$h_{\ast }(\omega _{1}^{1})=\frac{\alpha _{2}}{2}(\equiv \frac{1}{2}\omega
_{1}\func{mod}\Lambda _{G}^{r})$ (resp. $=\frac{\alpha _{1}}{2}(\equiv \frac{%
1}{2}\omega _{2}\func{mod}\Lambda _{G}^{r})$),
\end{quote}

\noindent and get

\begin{quote}
$h_{\ast }(\Lambda _{C_{\exp (u)}^{Rad}}^{e}(\mathbb{Q}))=\{\lambda \omega
_{1}\mid \lambda \in \mathbb{Q}\}$ (resp. $=\{\lambda \omega _{2}\mid
\lambda \in \mathbb{Q}\}$).
\end{quote}

\noindent It follows that the set $H_{u}$ consists of the single element $%
\omega _{1}^{1}$ whose deficiency in the group $SU(2)\times S^{1}$ is $2$.

Consequently, $\ker \pi =\mathbb{Z}_{2}$ with generator $\exp _{1}(\omega
_{1}^{1})\times \exp _{2}(-\frac{1}{2}\omega _{1})$ (resp. $\exp _{1}(\omega
_{1}^{1})\times \exp _{2}(-\frac{1}{2}\omega _{2})$).

\bigskip

\noindent \textbf{Case 2.} $G=F_{4}$

i) If $u=\frac{\omega _{1}}{4}$, the local type of the centralizer $C_{\exp
(u)}$ is $Spin(7)\times S^{1}$. Accordingly, assume that the set of
fundamental dominant weights of the semisimple part $Spin(7)$ is $\Omega
=\{\omega _{1}^{1},\omega _{2}^{1},\omega _{3}^{1}\}$. Applying Lemma 4.1 we
get the expressions of $h_{\ast }(\omega )$ with $\omega \in \Pi
_{Spin(7)}=\{\omega _{1}^{1}\}$ by the simple roots of $F_{4}$:

\begin{quote}
$h_{\ast }(\omega _{1}^{1})=\frac{3}{2}\alpha _{2}+2\alpha _{3}+\alpha
_{4}(\equiv \frac{1}{2}\omega _{1}\func{mod}\Lambda _{G}^{r})$
\end{quote}

\noindent and get

\begin{quote}
$h_{\ast }(\Lambda _{C_{\exp (u)}^{Rad}}^{e}(\mathbb{Q}))=\{\lambda \omega
_{1}\mid \lambda \in \mathbb{Q}\}$.
\end{quote}

\noindent It follows that the set $H_{u}$ consists of the single element $%
\omega _{1}^{1}$ whose deficiency in the group $Spin(7)$ is $2$.

Consequently, $\ker \pi =\mathbb{Z}_{2}$ with generator $\exp _{1}(\omega
_{1}^{1})\times \exp _{2}(-\frac{1}{2}\omega _{1})$.

ii) If $u=\frac{\omega _{2}}{8}$ (resp. $\frac{\omega _{3}}{4}$), the local
type of the centralizer $C_{\exp (u)}$ is $SU(2)\times SU(3)\times S^{1}$.
Accordingly, assume that the set of fundamental dominant weights of the
semisimple part $SU(2)\times SU(3)$ is $\Omega =\{\omega _{1}^{1}\}\amalg
\{\omega _{1}^{2},\omega _{2}^{2}\}$. Applying Lemma 4.1 we get the
expressions of $h_{\ast }(\omega )$ with $\omega \in \Pi _{SU(2)}\sqcup \Pi
_{SU(3)}$ by the simple roots of $F_{4}$:

\begin{quote}
$h_{\ast }(\omega _{1}^{1})=\frac{1}{2}\alpha _{1}(\equiv \frac{1}{2}\omega
_{2}\func{mod}\Lambda _{G}^{r})$

$h_{\ast }(\omega _{1}^{2})=\frac{2}{3}\alpha _{3}+\frac{1}{3}\alpha
_{4}(\equiv \frac{1}{3}\omega _{2}\func{mod}\Lambda _{G}^{r})$

$h_{\ast }(\omega _{2}^{2})=\frac{1}{3}\alpha _{3}+\frac{2}{3}\alpha
_{4}(\equiv \frac{2}{3}\omega _{2}\func{mod}\Lambda _{G}^{r})$
\end{quote}

(resp.

\begin{quote}
$h_{\ast }(\omega _{1}^{1})=\frac{1}{2}\alpha _{4}(\equiv \frac{1}{2}\omega
_{3}\func{mod}\Lambda _{G}^{r})$

$h_{\ast }(\omega _{1}^{2})=\frac{2}{3}\alpha _{1}+\frac{1}{3}\alpha
_{2}(\equiv \frac{1}{3}\omega _{3}\func{mod}\Lambda _{G}^{r})$

$h_{\ast }(\omega _{2}^{2})=\frac{1}{3}\alpha _{1}+\frac{2}{3}\alpha
_{2}(\equiv \frac{2}{3}\omega _{3}\func{mod}\Lambda _{G}^{r})$)
\end{quote}

\noindent and that

\begin{quote}
${\small h}_{\ast }{\small (\Lambda }_{C_{\exp (u)}^{Rad}}^{e}{\small %
(Q))=\{\lambda \omega }_{2}{\small \mid \lambda \in Q\}}$ {\small (resp. }$%
{\small =\{\lambda \omega }_{3}{\small \mid \lambda \in Q\}}${\small )}.
\end{quote}

\noindent It follows that the set $H_{u}$ consists of $5$ elements in which
the one $\omega _{1}^{1}\oplus \omega _{1}^{2}$ has deficiency $6$ in the
group $SU(2)\times SU(3)$.

Consequently, $\ker \pi =\mathbb{Z}_{6}$ with generator $\exp _{1}(\omega
_{1}^{1})\times \exp _{2}(\omega _{1}^{2})\times \exp _{3}(-\frac{5}{6}%
\omega _{2})$ by (4.3).

iii) If $u=\frac{\omega _{4}}{4}$, the local type of the centralizer $%
C_{\exp (u)}$ is $Sp(3)\times S^{1}$. Accordingly, assume that the set of
fundamental dominant weights of the semisimple part $Sp(3)$ is $\Omega
=\{\omega _{1}^{1},\omega _{2}^{1},\omega _{3}^{1}\}$. Applying Lemma 4.1 we
get the expressions of $h_{\ast }(\omega )$ with $\omega \in \Pi
_{Sp(3)}=\{\omega _{3}^{1}\}$ by the simple roots of $F_{4}$:

\begin{quote}
$h_{\ast }(\omega _{3}^{1})=\frac{1}{2}\alpha _{1}+\alpha _{2}+\frac{3}{2}%
\alpha _{3}(\equiv \frac{1}{2}\omega _{4}\func{mod}\Lambda _{G}^{r})$
\end{quote}

\noindent and that

\begin{quote}
$h_{\ast }(\Lambda _{C_{\exp (u)}^{Rad}}^{e}(\mathbb{Q}))=\{\lambda \omega
_{4}\mid \lambda \in \mathbb{Q}\}$.
\end{quote}

\noindent It follows that the set $H_{u}$ consists of a single element $%
\omega _{3}^{1}$ whose deficiency in the group $Sp(3)$ is $2$.

Consequently, $\ker \pi =\mathbb{Z}_{2}$ with generator $\exp _{1}(\omega
_{3}^{1})\times \exp _{2}(-\frac{1}{2}\omega _{4})$.

\bigskip

\noindent \textbf{Case 3.} $G=E_{6}$

i) If $u=\frac{\omega _{1}}{2}$ (resp. $\frac{\omega _{6}}{2}$), the local
type of the centralizer $C_{\exp (u)}$ is $Spin(10)\times S^{1}$.
Accordingly, assume that the set of fundamental dominant weights of the
semisimple part $Spin(10)$ is $\{\omega _{1}^{1},\omega _{2}^{1},\omega
_{3}^{1},\omega _{4}^{1},\omega _{5}^{1}\}$. Applying Lemma 4.1 we get the
expressions of $h_{\ast }(\omega )$ with $\omega \in \Pi
_{Spin(10)}=\{\omega _{1}^{1},\omega _{4}^{1},\omega _{5}^{1}\}$ by the
simple roots of $E_{6}$:

\begin{quote}
${\small h}_{\ast }{\small (\omega }_{1}^{1}{\small )=}\frac{1}{2}{\small %
\alpha }_{2}{\small +}\frac{1}{2}{\small \alpha }_{3}{\small +\alpha }_{4}%
{\small +\alpha }_{5}{\small +\alpha }_{6}{\small (\equiv }\frac{3}{2}%
{\small \omega }_{1}\func{mod}{\small \Lambda }_{G}^{r}{\small )}$

${\small h}_{\ast }{\small (\omega }_{4}^{1}{\small )=}\frac{3}{4}{\small %
\alpha }_{2}{\small +}\frac{5}{4}{\small \alpha }_{3}{\small +}\frac{3}{2}%
{\small \alpha }_{4}{\small +\alpha }_{5}{\small +}\frac{1}{2}{\small \alpha 
}_{6}{\small (\equiv }\frac{3}{4}{\small \omega }_{1}\func{mod}{\small %
\Lambda }_{G}^{r}{\small )}$

${\small h}_{\ast }{\small (\omega }_{5}^{1}{\small )=}\frac{5}{4}{\small %
\alpha }_{2}{\small +}\frac{3}{4}{\small \alpha }_{3}{\small +}\frac{3}{2}%
{\small \alpha }_{4}{\small +\alpha }_{5}{\small +}\frac{1}{2}{\small \alpha 
}_{6}{\small (\equiv }\frac{9}{4}{\small \omega }_{1}\func{mod}{\small %
\Lambda }_{G}^{r}{\small )}$
\end{quote}

(resp.

\begin{quote}
${\small h}_{\ast }{\small (\omega }_{1}^{1}{\small )=\alpha }_{1}{\small +}%
\frac{1}{2}{\small \alpha }_{2}{\small +\alpha }_{3}{\small +\alpha }_{4}%
{\small +}\frac{1}{2}{\small \alpha }_{5}{\small (\equiv }\frac{3}{2}{\small %
\omega }_{6}\func{mod}{\small \Lambda }_{G}^{r}{\small )}$

${\small h}_{\ast }{\small (\omega }_{4}^{1}{\small )=}\frac{1}{2}{\small %
\alpha }_{1}{\small +}\frac{5}{4}{\small \alpha }_{2}{\small +\alpha }_{3}%
{\small +}\frac{3}{2}{\small \alpha }_{4}{\small +}\frac{3}{4}{\small \alpha 
}_{5}{\small (\equiv }\frac{9}{4}{\small \omega }_{6}\func{mod}{\small %
\Lambda }_{G}^{r}{\small )}$

${\small h}_{\ast }{\small (\omega }_{5}^{1}{\small )=}\frac{1}{2}{\small %
\alpha }_{1}{\small +}\frac{3}{4}{\small \alpha }_{2}{\small +\alpha }_{3}%
{\small +}\frac{3}{2}{\small \alpha }_{4}{\small +}\frac{5}{4}{\small \alpha 
}_{5}{\small (\equiv }\frac{3}{4}{\small \omega }_{6}\func{mod}{\small %
\Lambda }_{G}^{r}{\small )}${\small )}
\end{quote}

\noindent and get

\begin{quote}
${\small h}_{\ast }{\small (\Lambda }_{C_{\exp (u)}^{Rad}}^{e}{\small %
(Q))=\{\lambda \omega }_{1}{\small \mid \lambda \in Q\}}${\small \ (resp. }$%
{\small \{\lambda \omega }_{6}{\small \mid \lambda \in Q\}}${\small )}
\end{quote}

\noindent It follows that the set $H_{u}$ consists of three elements $\omega
_{1}^{1},\omega _{4}^{1},\omega _{5}^{1}$ whose deficiencies in the group $%
Spin(10)$ are $2,4,4$, respectively.

Therefore, $\ker \pi =\mathbb{Z}_{4}$ with generator $\exp _{1}(\omega
_{5}^{1})\times \exp _{2}(-\frac{9}{4}\omega _{1})$ (resp. $\exp _{1}(\omega
_{5}^{1})\times \exp _{2}(-\frac{9}{4}\omega _{6})$) by (4.3).

ii) If $u=\frac{\omega _{2}}{4}$, the local type of the centralizer $C_{\exp
(u)}$ is $SU(6)\times S^{1}$. Accordingly, assume that the set of
fundamental dominant weights of the semisimple part $SU(6)$ is $\Omega
=\{\omega _{1}^{1},\omega _{2}^{1},\omega _{3}^{1},\omega _{4}^{1},\omega
_{5}^{1}\}$. Applying Lemma 4.1 we get the expressions of $h_{\ast }(\omega
) $ with $\omega \in \Pi _{SU(6)}$ by the simple roots of $E_{6}$:

\begin{quote}
${\small h}_{\ast }{\small (\omega }_{1}^{1}{\small )=}\frac{5}{6}{\small %
\alpha }_{1}{\small +}\frac{2}{3}{\small \alpha }_{3}{\small +}\frac{1}{2}%
{\small \alpha }_{4}{\small +}\frac{1}{3}{\small \alpha }_{5}{\small +}\frac{%
1}{6}{\small \alpha }_{6}{\small (\equiv \omega }_{1}{\small -}\frac{1}{2}%
{\small \omega }_{2}\func{mod}{\small \Lambda }_{G}^{r}{\small )}$

${\small h}_{\ast }{\small (\omega }_{2}^{1}{\small )=}\frac{2}{3}{\small %
\alpha }_{1}{\small +}\frac{4}{3}{\small \alpha }_{3}{\small +\alpha }_{4}%
{\small +}\frac{2}{3}{\small \alpha }_{5}{\small +}\frac{1}{3}{\small \alpha 
}_{6}{\small (\equiv \omega }_{3}\func{mod}{\small \Lambda }_{G}^{r}{\small )%
}$

${\small h}_{\ast }{\small (\omega }_{3}^{1}{\small )=}\frac{1}{2}{\small %
\alpha }_{1}{\small +\alpha }_{3}{\small +}\frac{3}{2}{\small \alpha }_{4}%
{\small +\alpha }_{5}{\small +}\frac{1}{2}{\small \alpha }_{6}{\small %
(\equiv }\frac{1}{2}{\small \omega }_{2}\func{mod}{\small \Lambda }_{G}^{r}%
{\small )}$

${\small h}_{\ast }{\small (\omega }_{4}^{1}{\small )=}\frac{1}{3}{\small %
\alpha }_{1}{\small +}\frac{2}{3}{\small \alpha }_{3}{\small +\alpha }_{4}%
{\small +}\frac{4}{3}{\small \alpha }_{5}{\small +}\frac{2}{3}{\small \alpha 
}_{6}{\small (\equiv \omega }_{5}\func{mod}{\small \Lambda }_{G}^{r}{\small )%
}$

${\small h}_{\ast }{\small (\omega }_{5}^{1}{\small )=}\frac{1}{6}{\small %
\alpha }_{1}{\small +}\frac{1}{3}{\small \alpha }_{3}{\small +}\frac{1}{2}%
{\small \alpha }_{4}{\small +}\frac{2}{3}{\small \alpha }_{5}{\small +}\frac{%
5}{6}{\small \alpha }_{6}{\small (\equiv \omega }_{6}{\small -}\frac{1}{2}%
{\small \omega }_{2}\func{mod}{\small \Lambda }_{G}^{r}{\small )}$
\end{quote}

\noindent and get

\begin{quote}
$h_{\ast }(\Lambda _{C_{\exp (u)}^{Rad}}^{e}(\mathbb{Q}))=\{\lambda \omega
_{2}\mid \lambda \in \mathbb{Q}\}$.
\end{quote}

\noindent It follows that the set $H_{u}$ consists of a single elements $%
\omega _{3}^{1}$ whose deficiency in the group $SU(6)$ is $2$.

Therefore, $\ker \pi =\mathbb{Z}_{2}$ with generator $\exp _{1}(\omega
_{3}^{1})\times \exp _{2}(-\frac{1}{2}\omega _{2})$.

iii) If $u=\frac{\omega _{3}}{4}$ (resp. $\frac{\omega _{5}}{4}$), the local
type of the centralizer $C_{\exp (u)}$ is $SU(2)\times SU(5)\times S^{1}$.
Accordingly, assume that the set of fundamental dominant weights of the
semisimple part $SU(2)\times SU(5)$ is $\Omega =\{\omega _{1}^{1}\}\amalg
\{\omega _{1}^{2},\omega _{2}^{2},\omega _{3}^{2},\omega _{4}^{2}\}$.
Applying Lemma 4.1 we get the expressions of $h_{\ast }(\omega )$ with $%
\omega \in \Pi _{SU(2)}\sqcup \Pi _{SU(5)}$ by the simple roots of $E_{6}$:

\begin{quote}
${\small h}_{\ast }{\small (\omega }_{1}^{1}{\small )=}\frac{1}{2}{\small %
\alpha }_{1}{\small (\equiv }\frac{3}{2}{\small \omega }_{3}\func{mod}%
{\small \Lambda }_{G}^{r}{\small )}$

${\small h}_{\ast }{\small (\omega }_{1}^{2}{\small )=}\frac{4}{5}{\small %
\alpha }_{2}{\small +}\frac{3}{5}{\small \alpha }_{4}{\small +}\frac{2}{5}%
{\small \alpha }_{5}{\small +}\frac{1}{5}{\small \alpha }_{6}{\small (\equiv 
}\frac{12}{5}{\small \omega }_{3}\func{mod}{\small \Lambda }_{G}^{r}{\small )%
}$

${\small h}_{\ast }{\small (\omega }_{2}^{2}{\small )=}\frac{3}{5}{\small %
\alpha }_{2}{\small +}\frac{6}{5}{\small \alpha }_{4}{\small +}\frac{4}{5}%
{\small \alpha }_{5}{\small +}\frac{2}{5}{\small \alpha }_{6}{\small (\equiv 
}\frac{9}{5}{\small \omega }_{3}\func{mod}{\small \Lambda }_{G}^{r}{\small )}
$

${\small h}_{\ast }{\small (\omega }_{3}^{2}{\small )=}\frac{2}{5}{\small %
\alpha }_{2}{\small +}\frac{4}{5}{\small \alpha }_{4}{\small +}\frac{6}{5}%
{\small \alpha }_{5}{\small +}\frac{3}{5}{\small \alpha }_{6}{\small (\equiv 
}\frac{6}{5}{\small \omega }_{3}\func{mod}{\small \Lambda }_{G}^{r}{\small )}
$

${\small h}_{\ast }{\small (\omega }_{4}^{2}{\small )=}\frac{1}{5}{\small %
\alpha }_{2}{\small +}\frac{2}{5}{\small \alpha }_{4}{\small +}\frac{3}{5}%
{\small \alpha }_{5}{\small +}\frac{4}{5}{\small \alpha }_{6}{\small (\equiv 
}\frac{3}{5}{\small \omega }_{3}\func{mod}{\small \Lambda }_{G}^{r}{\small )}
$
\end{quote}

(resp.

\begin{quote}
${\small h}_{\ast }{\small (\omega }_{1}^{1}{\small )=}\frac{1}{2}{\small %
\alpha }_{6}{\small (\equiv }\frac{3}{2}{\small \omega }_{5}\func{mod}%
{\small \Lambda }_{G}^{r}{\small )}$

${\small h}_{\ast }{\small (\omega }_{1}^{2}{\small )=}\frac{4}{5}{\small %
\alpha }_{1}{\small +}\frac{1}{5}{\small \alpha }_{2}{\small +}\frac{3}{5}%
{\small \alpha }_{3}{\small +}\frac{2}{5}{\small \alpha }_{4}{\small (\equiv 
}\frac{3}{5}{\small \omega }_{5}\func{mod}{\small \Lambda }_{G}^{r}{\small )}
$

${\small h}_{\ast }{\small (\omega }_{2}^{2}{\small )=}\frac{3}{5}{\small %
\alpha }_{1}{\small +}\frac{2}{5}{\small \alpha }_{2}{\small +}\frac{6}{5}%
{\small \alpha }_{3}{\small +}\frac{4}{5}{\small \alpha }_{4}{\small (\equiv 
}\frac{6}{5}{\small \omega }_{5}\func{mod}{\small \Lambda }_{G}^{r}{\small )}
$

${\small h}_{\ast }{\small (\omega }_{1}^{3}{\small )=}\frac{2}{5}{\small %
\alpha }_{1}{\small +}\frac{3}{5}{\small \alpha }_{2}{\small +}\frac{4}{5}%
{\small \alpha }_{3}{\small +}\frac{6}{5}{\small \alpha }_{4}{\small (\equiv 
}\frac{9}{5}{\small \omega }_{5}\func{mod}{\small \Lambda }_{G}^{r}{\small )}
$

${\small h}_{\ast }{\small (\omega }_{2}^{3}{\small )=}\frac{1}{5}{\small %
\alpha }_{1}{\small +}\frac{4}{5}{\small \alpha }_{2}{\small +}\frac{2}{5}%
{\small \alpha }_{3}{\small +}\frac{3}{5}{\small \alpha }_{4}{\small (\equiv 
}\frac{12}{5}{\small \omega }_{5}\func{mod}{\small \Lambda }_{G}^{r}{\small )%
}${\small )}
\end{quote}

\noindent and get

\begin{quote}
${\small h}_{\ast }{\small (\Lambda }_{C_{\exp (u)}^{Rad}}^{e}{\small %
(Q))=\{\lambda \omega }_{3}{\small \mid \lambda \in Q\}}${\small \ (resp. }$%
{\small =\{\lambda \omega }_{5}{\small \mid \lambda \in Q\}}${\small ).}
\end{quote}

\noindent It follows that the set $H_{u}$ contains $9$ elements in which the
one $\omega _{1}^{1}\oplus \omega _{1}^{2}$ has deficiency $10$ in the group 
$SU(2)\times SU(5)$.

Consequently, $\ker \pi =\mathbb{Z}_{10}$ with generator $\exp (\omega
_{1}^{1})\times \exp _{2}(\omega _{1}^{2})\times \exp _{3}(-\frac{9}{10}%
\omega _{3})$ (resp. $\exp _{1}(\omega _{1}^{1})\times \exp _{2}(\omega
_{1}^{2})\times \exp _{3}(-\frac{21}{10}\omega _{5})$) by (4.3).

iv) If $u=\frac{\omega _{4}}{6}$, the local type of the centralizer $C_{\exp
(u)}$ is $SU(2)\times SU(3)\times SU(3)\times S^{1}$. Accordingly, assume
that the set of fundamental dominant weights of the semisimple part $%
SU(2)\times SU(3)\times SU(3)$ is $\Omega =\{\omega _{1}^{1}\}\amalg
\{\omega _{1}^{2},\omega _{2}^{2}\}\amalg \{\omega _{1}^{3},\omega
_{2}^{3}\} $. Applying Lemma 4.1, we get the expressions of $h_{\ast
}(\omega )$ with $\omega \in \Pi _{SU(2)}\sqcup $ $\Pi _{SU(3)}\sqcup \Pi
_{SU(3)}$ by simple roots of $E_{6}$:

\begin{quote}
$h_{\ast }(\omega _{1}^{1})=\frac{1}{2}\alpha _{2}(\equiv \frac{1}{2}\omega
_{4}\func{mod}\Lambda _{G}^{r})$

$h_{\ast }(\omega _{1}^{2})=\frac{2}{3}\alpha _{1}+\frac{1}{3}\alpha
_{3}(\equiv \omega _{1}-\frac{1}{3}\omega _{4}\func{mod}\Lambda _{G}^{r})$

$h_{\ast }(\omega _{2}^{2})=\frac{1}{3}\alpha _{1}+\frac{2}{3}\alpha
_{3}(\equiv \omega _{3}-\frac{2}{3}\omega _{4}\func{mod}\Lambda _{G}^{r})$

$h_{\ast }(\omega _{1}^{3})=\frac{2}{3}\alpha _{5}+\frac{1}{3}\alpha
_{6}(\equiv \omega _{5}-\frac{2}{3}\omega _{4}\func{mod}\Lambda _{G}^{r})$

$h_{\ast }(\omega _{2}^{3})=\frac{1}{3}\alpha _{5}+\frac{2}{3}\alpha
_{6}(\equiv \omega _{6}-\frac{1}{3}\omega _{4}\func{mod}\Lambda _{G}^{r})$
\end{quote}

\noindent and get

\begin{quote}
$h_{\ast }(\Lambda _{C_{\exp (u)}^{Rad}}^{e}(\mathbb{Q}))=\{\lambda \omega
_{4}\mid \lambda \in \mathbb{Q}\}$.
\end{quote}

\noindent It follows that the set $H_{u}$ consists of the $5$ elements

\begin{quote}
$\omega _{1}^{1},$ $\omega _{1}^{2}\oplus \omega _{2}^{3},$ $\omega
_{2}^{2}\oplus \omega _{1}^{3},\omega _{1}^{1}\oplus \omega _{1}^{2}\oplus
\omega _{2}^{3}$, $\omega _{1}^{1}\oplus \omega _{2}^{2}\oplus \omega
_{1}^{3}$
\end{quote}

\noindent whose deficiencies in the group $SU(2)\times SU(3)\times SU(3)$
are $2,3,3,6,6$, respectively.

Consequently, $\ker \pi =\mathbb{Z}_{6}$ with generator $\exp (\omega
_{1}^{1})\times \exp _{2}(\omega _{1}^{2})\times \exp _{3}(\omega
_{2}^{3})\times \exp _{4}(\frac{1}{6}\omega _{4})$ by (4.3).

\bigskip

\noindent \textbf{Case 4.} $G=E_{7}$

i) If $u=\frac{\omega _{1}}{4}$, the local type of the centralizer $C_{\exp
(u)}$ is $Spin(12)\times S^{1}$. Accordingly, assume that the set of
fundamental dominant weights of the semisimple part $Spin(12)$ is $\Omega
=\{\omega _{1}^{1},\omega _{2}^{1},\omega _{3}^{1},\omega _{4}^{1},\omega
_{5}^{1},\omega _{6}^{1}\}$. Applying Lemma 4.1 we get the expressions of $%
h_{\ast }(\omega )$ with $\omega \in \Pi _{Spin(12)}=\{\omega
_{1}^{1},\omega _{5}^{1},\omega _{6}^{1}\}$ by the simple roots of $E_{7}$:

\begin{quote}
${\small h}_{\ast }{\small (\omega }_{1}^{1}{\small )=}\frac{1}{2}{\small %
\alpha }_{2}{\small +}\frac{1}{2}{\small \alpha }_{3}{\small +\alpha }_{4}%
{\small +\alpha }_{5}{\small +\alpha }_{6}{\small +\alpha }_{7}{\small %
(\equiv \omega }_{7}{\small -}\frac{1}{2}{\small \omega }_{1}\func{mod}%
{\small \Lambda }_{G}^{r}{\small )}$

${\small h}_{\ast }{\small (\omega }_{5}^{1}{\small )=\alpha }_{2}{\small +}%
\frac{3}{2}{\small \alpha }_{3}{\small +2\alpha }_{4}{\small +}\frac{3}{2}%
{\small \alpha }_{5}{\small +\alpha }_{6}{\small +}\frac{1}{2}{\small \alpha 
}_{7}{\small (\equiv }\frac{1}{2}{\small \omega }_{1}\func{mod}{\small %
\Lambda }_{G}^{r}{\small )}$

${\small h}_{\ast }{\small (\omega }_{6}^{1}{\small )=}\frac{3}{2}{\small %
\alpha }_{2}{\small +\alpha }_{3}{\small +2\alpha }_{4}{\small +}\frac{3}{2}%
{\small \alpha }_{5}{\small +\alpha }_{6}{\small +}\frac{1}{2}{\small \alpha 
}_{7}{\small (\equiv \omega }_{2}\func{mod}{\small \Lambda }_{G}^{r}{\small )%
}$
\end{quote}

\noindent and get

\begin{quote}
$h_{\ast }(\Lambda _{C_{\exp (u)}^{Rad}}^{e}(\mathbb{Q}))=\{\lambda \omega
_{1}\mid \lambda \in \mathbb{Q}\}$.
\end{quote}

\noindent It follows that the set $H_{u}$ consists of the single element $%
\omega _{5}^{1}$ whose deficiency in the group $Spin(12)$ is $2$.

Consequently, $\ker \pi =\mathbb{Z}_{2}$ generated by $\exp _{1}(\omega
_{5}^{1})\times \exp _{2}(-\frac{1}{2}\omega _{1})$.

ii) If $u=\frac{\omega _{2}}{4}$, the local type of the centralizer $C_{\exp
(u)}$ is $SU(7)\times S^{1}$. Accordingly, assume that the set of
fundamental dominant weights of the semisimple part $SU(7)$ is $\Omega
=\{\omega _{1}^{1},\omega _{2}^{1},\omega _{3}^{1},\omega _{4}^{1},\omega
_{5}^{1},\omega _{6}^{1}\}$. Applying Lemma 4.1 we get the expressions of $%
h_{\ast }(\omega )$ with $\omega \in \Pi _{SU(7)}$ by the simple roots of $%
E_{7}$:

\begin{quote}
${\small h}_{\ast }{\small (\omega }_{1}^{1}{\small )=}\frac{6}{7}{\small %
\alpha }_{1}{\small +}\frac{5}{7}{\small \alpha }_{3}{\small +}\frac{4}{7}%
{\small \alpha }_{4}{\small +}\frac{3}{7}{\small \alpha }_{5}{\small +}\frac{%
2}{7}{\small \alpha }_{6}{\small +}\frac{1}{7}{\small \alpha }_{7}{\small %
(\equiv }\frac{10}{7}{\small \omega }_{2}\func{mod}{\small \Lambda }_{G}^{r}%
{\small )}$

${\small h}_{\ast }{\small (\omega }_{2}^{1}{\small )=}\frac{5}{7}{\small %
\alpha }_{1}{\small +}\frac{10}{7}{\small \alpha }_{3}{\small +}\frac{8}{7}%
{\small \alpha }_{4}{\small +}\frac{6}{7}{\small \alpha }_{5}{\small +}\frac{%
4}{7}{\small \alpha }_{6}{\small +}\frac{2}{7}{\small \alpha }_{7}{\small %
(\equiv }\frac{6}{7}{\small \omega }_{2}\func{mod}{\small \Lambda }_{G}^{r}%
{\small )}$

${\small h}_{\ast }{\small (\omega }_{3}^{1}{\small )=}\frac{4}{7}{\small %
\alpha }_{1}{\small +}\frac{8}{7}{\small \alpha }_{3}{\small +}\frac{12}{7}%
{\small \alpha }_{4}{\small +}\frac{9}{7}{\small \alpha }_{5}{\small +}\frac{%
6}{7}{\small \alpha }_{6}{\small +}\frac{3}{7}{\small \alpha }_{7}{\small %
(\equiv }\frac{2}{7}{\small \omega }_{2}\func{mod}{\small \Lambda }_{G}^{r}%
{\small )}$

${\small h}_{\ast }{\small (\omega }_{4}^{1}{\small )=}\frac{3}{7}{\small %
\alpha }_{1}{\small +}\frac{6}{7}{\small \alpha }_{3}{\small +}\frac{9}{7}%
{\small \alpha }_{4}{\small +}\frac{12}{7}{\small \alpha }_{5}{\small +}%
\frac{8}{7}{\small \alpha }_{6}{\small +}\frac{4}{7}{\small \alpha }_{7}%
{\small (\equiv }\frac{12}{7}{\small \omega }_{2}\func{mod}{\small \Lambda }%
_{G}^{r}{\small )}$

${\small h}_{\ast }{\small (\omega }_{5}^{1}{\small )=}\frac{2}{7}{\small %
\alpha }_{1}{\small +}\frac{4}{7}{\small \alpha }_{3}{\small +}\frac{6}{7}%
{\small \alpha }_{4}{\small +}\frac{8}{7}{\small \alpha }_{5}{\small +}\frac{%
10}{7}{\small \alpha }_{6}{\small +}\frac{5}{7}{\small \alpha }_{7}{\small %
(\equiv }\frac{8}{7}{\small \omega }_{2}\func{mod}{\small \Lambda }_{G}^{r}%
{\small )}$

${\small h}_{\ast }{\small (\omega }_{6}^{1}{\small )=}\frac{1}{7}{\small %
\alpha }_{1}{\small +}\frac{2}{7}{\small \alpha }_{3}{\small +}\frac{3}{7}%
{\small \alpha }_{4}{\small +}\frac{4}{7}{\small \alpha }_{5}{\small +}\frac{%
5}{7}{\small \alpha }_{6}{\small +}\frac{6}{7}{\small \alpha }_{7}{\small %
(\equiv }\frac{4}{7}{\small \omega }_{2}\func{mod}{\small \Lambda }_{G}^{r}%
{\small )}$
\end{quote}

\noindent and get

\begin{quote}
$h_{\ast }(\Lambda _{C_{\exp (u)}^{Rad}}^{e}(\mathbb{Q}))=\{\lambda \omega
_{2}\mid \lambda \in \mathbb{Q}\}$.
\end{quote}

\noindent It follows that the set $H_{u}$ consists of $6$ elements whose
deficiencies in the group $SU(7)$ are all $7$.

Consequently, $\ker \pi =\mathbb{Z}_{7}$ with generator $\exp _{1}(\omega
_{3}^{1})\times \exp _{2}(-\frac{2}{7}\omega _{2})$ by (4.3).

iii) If $u=\frac{\omega _{3}}{6}$, the local type of the centralizer $%
C_{\exp (u)}$ is $SU(2)\times SU(6)\times S^{1}$. Accordingly, assume that
the set of fundamental dominant weights of the semisimple part $SU(2)\times
SU(6)$ is $\Omega =\{\omega _{1}^{1}\}\amalg \{\omega _{1}^{2},\omega
_{2}^{2},\omega _{3}^{2},\omega _{4}^{2},\omega _{5}^{2}\}$. Applying Lemma
4.1 we get the expressions of $h_{\ast }(\omega )$ with $\omega \in \Pi
_{SU(2)}\sqcup \Pi _{SU(6)}$ by the simple roots of $E_{7}$:

\begin{quote}
${\small h}_{\ast }{\small (\omega }_{1}^{1}{\small )=}\frac{1}{2}{\small %
\alpha }_{1}{\small (\equiv }\frac{1}{2}{\small \omega }_{3}\func{mod}%
{\small \Lambda }_{G}^{r}{\small )}$

${\small h}_{\ast }{\small (\omega }_{1}^{2}{\small )=}\frac{5}{6}{\small %
\alpha }_{2}{\small +}\frac{2}{3}{\small \alpha }_{4}{\small +}\frac{1}{2}%
{\small \alpha }_{5}{\small +}\frac{1}{3}{\small \alpha }_{6}{\small +}\frac{%
1}{6}{\small \alpha }_{7}{\small (\equiv \omega }_{2}{\small -}\frac{2}{3}%
{\small \omega }_{3}\func{mod}{\small \Lambda }_{G}^{r}{\small )}$

${\small h}_{\ast }{\small (\omega }_{2}^{2}{\small )=}\frac{2}{3}{\small %
\alpha }_{2}{\small +}\frac{4}{3}{\small \alpha }_{4}{\small +\alpha }_{5}%
{\small +}\frac{2}{3}{\small \alpha }_{6}{\small +}\frac{1}{3}{\small \alpha 
}_{7}{\small (\equiv }\frac{2}{3}{\small \omega }_{3}\func{mod}{\small %
\Lambda }_{G}^{r}{\small )}$

${\small h}_{\ast }{\small (\omega }_{3}^{2}{\small )=}\frac{1}{2}{\small %
\alpha }_{2}{\small +\alpha }_{4}{\small +}\frac{3}{2}{\small \alpha }_{5}%
{\small +\alpha }_{6}{\small +}\frac{1}{2}{\small \alpha }_{7}{\small %
(\equiv \omega }_{5}\func{mod}{\small \Lambda }_{G}^{r}{\small )}$

${\small h}_{\ast }{\small (\omega }_{4}^{2}{\small )=}\frac{1}{3}{\small %
\alpha }_{2}{\small +}\frac{2}{3}{\small \alpha }_{4}{\small +\alpha }_{5}%
{\small +}\frac{4}{3}{\small \alpha }_{6}{\small +}\frac{2}{3}{\small \alpha 
}_{7}{\small (\equiv }\frac{1}{3}{\small \omega }_{3}\func{mod}{\small %
\Lambda }_{G}^{r}{\small )}$

${\small h}_{\ast }{\small (\omega }_{5}^{2}{\small )=}\frac{1}{6}{\small %
\alpha }_{2}{\small +}\frac{1}{3}{\small \alpha }_{4}{\small +}\frac{1}{2}%
{\small \alpha }_{5}{\small +}\frac{2}{3}{\small \alpha }_{6}{\small +}\frac{%
5}{6}{\small \alpha }_{7}{\small (\equiv \omega }_{7}{\small -}\frac{1}{3}%
{\small \omega }_{3}\func{mod}{\small \Lambda }_{G}^{r}{\small )}$
\end{quote}

\noindent and

\begin{quote}
$h_{\ast }(\Lambda _{C_{\exp (u)}^{Rad}}^{e}(\mathbb{Q}))=\{\lambda \omega
_{3}\mid \lambda \in \mathbb{Q}\}$.
\end{quote}

\noindent It follows that the set $H_{u}$ consists of $5$ elements, in which
the deficiency of $\omega _{1}^{1}\oplus \omega _{4}^{2}$ in the group $%
SU(2)\times SU(6)$ is $6$.

Consequently, $\ker \pi =\mathbb{Z}_{6}$ with generator $\exp _{1}(\omega
_{1}^{1})\times \exp _{2}(\omega _{4}^{2})\times \exp _{3}(-\frac{5}{6}%
\omega _{3})$.

iv) If $u=\frac{\omega _{4}}{8}$, the local type of the centralizer $C_{\exp
(u)}$ is $SU(2)\times SU(3)\times SU(4)\times S^{1}$. Accordingly, assume
that the set of fundamental dominant weights of the semisimple part $%
SU(2)\times SU(3)\times SU(4)$ is $\Omega =\{\omega _{1}^{1}\}\amalg
\{\omega _{1}^{2},\omega _{2}^{2}\}\amalg \{\omega _{1}^{3},\omega
_{2}^{3},\omega _{3}^{3}\}$. Applying Lemma 4.1 we get the expressions of $%
h_{\ast }(\omega )$ with $\omega \in \Pi _{SU(2)}\sqcup \Pi _{SU(3)}\sqcup
\Pi _{SU(4)}$ by the simple roots of $E_{7}$:

\begin{quote}
${\small h}_{\ast }{\small (\omega }_{1}^{1}{\small )=}\frac{1}{2}{\small %
\alpha }_{2}{\small (\equiv \omega }_{2}{\small -}\frac{1}{2}{\small \omega }%
_{4}\func{mod}{\small \Lambda }_{G}^{r}{\small )}$

${\small h}_{\ast }{\small (\omega }_{1}^{2}{\small )=}\frac{2}{3}{\small %
\alpha }_{1}{\small +}\frac{1}{3}{\small \alpha }_{3}{\small (\equiv }\frac{2%
}{3}{\small \omega }_{4}\func{mod}{\small \Lambda }_{G}^{r}{\small )}$

${\small h}_{\ast }{\small (\omega }_{2}^{2}{\small )=}\frac{1}{3}{\small %
\alpha }_{1}{\small +}\frac{2}{3}{\small \alpha }_{3}{\small (\equiv }\frac{1%
}{3}{\small \omega }_{4}\func{mod}{\small \Lambda }_{G}^{r}{\small )}$

${\small h}_{\ast }{\small (\omega }_{1}^{3}{\small )=}\frac{3}{4}{\small %
\alpha }_{5}{\small +}\frac{1}{2}{\small \alpha }_{6}{\small +}\frac{1}{4}%
{\small \alpha }_{7}{\small (\equiv \omega }_{5}{\small -}\frac{3}{4}{\small %
\omega }_{4}\func{mod}{\small \Lambda }_{G}^{r}{\small )}$

${\small h}_{\ast }{\small (\omega }_{2}^{3}{\small )=}\frac{1}{2}{\small %
\alpha }_{5}{\small +\alpha }_{6}{\small +}\frac{1}{2}{\small \alpha }_{7}%
{\small (\equiv }\frac{1}{2}{\small \omega }_{4}\func{mod}{\small \Lambda }%
_{G}^{r}{\small )}$

${\small h}_{\ast }{\small (\omega }_{3}^{3}{\small )=}\frac{1}{4}{\small %
\alpha }_{5}{\small +}\frac{1}{2}{\small \alpha }_{6}{\small +}\frac{3}{4}%
{\small \alpha }_{7}{\small (\equiv \omega }_{7}{\small -}\frac{1}{4}{\small %
\omega }_{4}\func{mod}{\small \Lambda }_{G}^{r}{\small )}$
\end{quote}

\noindent and get

\begin{quote}
$h_{\ast }(\Lambda _{C_{\exp (u)}^{Rad}}^{e}(\mathbb{Q}))=\{\lambda \omega
_{4}\mid \lambda \in \mathbb{Q}\}$.
\end{quote}

\noindent It follows that the set $H_{u}$ consists of $11$ elements in which
the one $\omega _{1}^{1}\oplus \omega _{1}^{2}\oplus \omega _{1}^{3}$

\noindent has deficiency $12$ in the group $SU(2)\times SU(3)\times SU(4)$.

Consequently, $\ker \pi =\mathbb{Z}_{12}$ with generators $\exp _{1}(\omega
_{1}^{1})\times \exp _{2}(\omega _{1}^{2})\times \exp _{3}(\omega
_{1}^{3})\times \exp _{3}(-\frac{7}{12}\omega _{4})$ by (4.3).

v) If $u=\frac{\omega _{5}}{6}$, the local type of the centralizer $C_{\exp
(u)}$ is $SU(3)\times SU(5)\times S^{1}$. Accordingly, assume that the set
of fundamental dominant weights of the semisimple part $SU(3)\times SU(5)$
is $\Omega =\{\omega _{1}^{1},\omega _{2}^{1}\}\amalg \{\omega
_{1}^{2},\omega _{2}^{2},\omega _{3}^{2},\omega _{4}^{2}\}$. Applying Lemma
4.1 we get the expressions of $h_{\ast }(\omega )$ with $\omega \in \Pi
_{SU(3)}\sqcup \Pi _{SU(5)}$ by the simple roots of $E_{7}$:

\begin{quote}
${\small h}_{\ast }{\small (\omega }_{1}^{1}{\small )=}\frac{2}{3}{\small %
\alpha }_{6}{\small +}\frac{1}{3}{\small \alpha }_{7}{\small (\equiv }\frac{4%
}{3}{\small \omega }_{5}\func{mod}{\small \Lambda }_{G}^{r}{\small )}$

${\small h}_{\ast }{\small (\omega }_{2}^{1}{\small )=}\frac{1}{3}{\small %
\alpha }_{6}{\small +}\frac{2}{3}{\small \alpha }_{7}{\small (\equiv }\frac{2%
}{3}{\small \omega }_{5}\func{mod}{\small \Lambda }_{G}^{r}{\small )}$

${\small h}_{\ast }{\small (\omega }_{1}^{2}{\small )=}\frac{4}{5}{\small %
\alpha }_{1}{\small +}\frac{1}{5}{\small \alpha }_{2}{\small +}\frac{3}{5}%
{\small \alpha }_{3}{\small +}\frac{2}{5}{\small \alpha }_{4}{\small (\equiv 
}\frac{8}{5}{\small \omega }_{5}\func{mod}{\small \Lambda }_{G}^{r}{\small )}
$

${\small h}_{\ast }{\small (\omega }_{2}^{2}{\small )=}\frac{1}{5}{\small %
\alpha }_{1}{\small +}\frac{4}{5}{\small \alpha }_{2}{\small +}\frac{2}{5}%
{\small \alpha }_{3}{\small +}\frac{3}{5}{\small \alpha }_{4}{\small (\equiv 
}\frac{2}{5}{\small \omega }_{5}\func{mod}{\small \Lambda }_{G}^{r}{\small )}
$

${\small h}_{\ast }{\small (\omega }_{3}^{2}{\small )=}\frac{3}{5}{\small %
\alpha }_{1}{\small +}\frac{2}{5}{\small \alpha }_{2}{\small +}\frac{6}{5}%
{\small \alpha }_{3}{\small +}\frac{4}{5}{\small \alpha }_{4}{\small (\equiv 
}\frac{6}{5}{\small \omega }_{5}\func{mod}{\small \Lambda }_{G}^{r}{\small )}
$

${\small h}_{\ast }{\small (\omega }_{4}^{2}{\small )=}\frac{2}{5}{\small %
\alpha }_{1}{\small +}\frac{3}{5}{\small \alpha }_{2}{\small +}\frac{4}{5}%
{\small \alpha }_{3}{\small +}\frac{6}{5}{\small \alpha }_{4}{\small (\equiv 
}\frac{4}{5}{\small \omega }_{5}\func{mod}{\small \Lambda }_{G}^{r}{\small )}
$
\end{quote}

\noindent and get

\begin{quote}
$h_{\ast }(\Lambda _{C_{\exp (u)}^{Rad}}^{e}(\mathbb{Q}))=\{\lambda \omega
_{5}\mid \lambda \in \mathbb{Q}\}$.
\end{quote}

\noindent It follows that the set $H_{u}$ consists of $14$ elements, in
which the deficiency of $\omega _{2}^{1}+\omega _{2}^{2}$ in the group $%
SU(3)\times SU(5)$ is $15$.

Consequently, $\ker \pi =\mathbb{Z}_{15}$ with generator $\exp _{1}(\omega
_{2}^{1})\times \exp _{2}(\omega _{2}^{2})\times \exp _{3}(-\frac{16}{15}%
\omega _{5})$ by (4.3)

vi) If $u=\frac{\omega _{6}}{4}$, the local type of the centralizer $C_{\exp
(u)}$ is $SU(2)\times Spin(10)\times S^{1}$. Accordingly, assume that the
set of fundamental dominant weights of the semisimple part $SU(2)\times
Spin(10)$ is $\Omega =\{\omega _{1}^{1}\}\amalg \{\omega _{1}^{2},\omega
_{2}^{2},\omega _{3}^{2},\omega _{4}^{2},\omega _{5}^{2}\}$. Applying Lemma
4.1 we get the expressions of $h_{\ast }(\omega )$ with $\omega \in \Pi
_{SU(2)}\sqcup \Pi _{Spin(10)}$ by the simple roots of $E_{7}$:

\begin{quote}
${\small h}_{\ast }{\small (\omega }_{1}^{1}{\small )=}\frac{1}{2}{\small %
\alpha }_{7}{\small (\equiv \omega }_{7}{\small -}\frac{1}{2}{\small \omega }%
_{6}\func{mod}{\small \Lambda }_{G}^{r}{\small )}$

${\small h}_{\ast }{\small (\omega }_{1}^{2}{\small )=\alpha }_{1}{\small +}%
\frac{1}{2}{\small \alpha }_{2}{\small +\alpha }_{3}{\small +\alpha }_{4}%
{\small +}\frac{1}{2}{\small \alpha }_{5}{\small (\equiv }\frac{1}{2}{\small %
\omega }_{6}\func{mod}{\small \Lambda }_{G}^{r}{\small )}$

${\small h}_{\ast }{\small (\omega }_{4}^{2}{\small )=}\frac{1}{2}{\small %
\alpha }_{1}{\small +}\frac{5}{4}{\small \alpha }_{2}{\small +\alpha }_{3}%
{\small +}\frac{3}{2}{\small \alpha }_{4}{\small +}\frac{3}{4}{\small \alpha 
}_{5}{\small (\equiv \omega }_{2}{\small -}\frac{3}{4}{\small \omega }_{6}%
\func{mod}{\small \Lambda }_{G}^{r}{\small )}$

${\small h}_{\ast }{\small (\omega }_{5}^{2}{\small )=}\frac{1}{2}{\small %
\alpha }_{1}{\small +}\frac{3}{4}{\small \alpha }_{2}{\small +\alpha }_{3}%
{\small +}\frac{3}{2}{\small \alpha }_{4}{\small +}\frac{5}{4}{\small \alpha 
}_{5}{\small (\equiv \omega }_{5}{\small -}\frac{5}{4}{\small \omega }_{6}%
\func{mod}{\small \Lambda }_{G}^{r}{\small )}$
\end{quote}

\noindent and get

\begin{quote}
$h_{\ast }(\Lambda _{C_{\exp (u)}^{Rad}}^{e}(\mathbb{Q}))=\{\lambda \omega
_{6}\mid \lambda \in \mathbb{Q}\}$.
\end{quote}

\noindent It follows that the set $H_{u}$ consists of the three elements

\begin{quote}
$\omega _{1}^{2}$, $\omega _{1}^{1}\oplus \omega _{4}^{2}$, $\omega
_{1}^{1}\oplus \omega _{5}^{2}$
\end{quote}

\noindent in which the one $\omega _{1}^{1}\oplus \omega _{5}^{2}$ has
deficiency $4$ in the group $SU(2)\times Spin(10)$.

Consequently, $\ker \pi =\mathbb{Z}_{4}$ with generator $\exp _{1}(\omega
_{1}^{1})\otimes \exp _{2}(\omega _{5}^{2})\times \exp _{3}(-\frac{3}{4}%
\omega _{6})$ by (4.3).

vii) If $u=\frac{\omega _{7}}{2}$, the local type of the centralizer $%
C_{\exp (u)}$ is $E_{6}\times S^{1}$. Accordingly, assume that the set of
fundamental dominant weights of the semisimple part $E_{6}$ is $\Omega
=\{\omega _{1}^{1},\omega _{2}^{1},\omega _{3}^{1},\omega _{4}^{1},\omega
_{5}^{1},\omega _{6}^{1}\}$. Applying Lemma 4.1 we get the expressions of $%
h_{\ast }(\omega )$ with $\omega \in \Pi _{E_{6}}=\{\omega _{1}^{1},\omega
_{6}^{1}\}$ by the simple roots of $E_{7}$:

\begin{quote}
${\small h}_{\ast }{\small (\omega }_{1}^{1}{\small )=}\frac{1}{3}{\small %
(4\alpha }_{1}{\small +3\alpha }_{2}{\small +5\alpha }_{3}{\small +6\alpha }%
_{4}{\small +4\alpha }_{5}{\small +2\alpha }_{6}{\small )(\equiv }\frac{4}{3}%
{\small \omega }_{7}\func{mod}{\small \Lambda }_{G}^{r}{\small )}$

${\small h}_{\ast }{\small (\omega }_{6}^{1}{\small )=}\frac{1}{3}{\small %
(2\alpha }_{1}{\small +3\alpha }_{2}{\small +4\alpha }_{3}{\small +6\alpha }%
_{4}{\small +5\alpha }_{5}{\small +4\alpha }_{6}{\small )(\equiv }\frac{2}{3}%
{\small \omega }_{7}\func{mod}{\small \Lambda }_{G}^{r}{\small )}$
\end{quote}

\noindent and

\begin{quote}
$h_{\ast }(\Lambda _{C_{\exp (u)}^{Rad}}^{e}(\mathbb{Q}))=\{\lambda \omega
_{7}\mid \lambda \in \mathbb{Q}\}$.
\end{quote}

\noindent It follows that the set $H_{u}$ consists of the two elements $%
\omega _{1}^{1}$ and $\omega _{6}^{1}$ whose deficiencies in the group $%
E_{6} $ are both $3$.

Consequently $\ker \pi =\mathbb{Z}_{3}$ with generator $\exp _{1}(\omega
_{1}^{1})\times \exp _{2}(-\frac{4}{3}\omega _{7})$ by (4.3).

\bigskip

\noindent \textbf{Case 5.} $G=E_{8}$

i) If $u=\frac{\omega _{1}}{4}$, the local type of the centralizer $C_{\exp
(u)}$ is $Spin(14)\times S^{1}$. Accordingly, assume that the set of
fundamental dominant weights of the semisimple part $Spin(14)$ is $\Omega
=\{\omega _{1}^{1},\omega _{2}^{1},\omega _{3}^{1},\omega _{4}^{1},\omega
_{5}^{1},\omega _{6}^{1},\omega _{7}^{1}\}$. Applying Lemma 4.1 we get the
expressions of $h_{\ast }(\omega )$ with $\omega \in \Pi
_{Spin(14)}=\{\omega _{1}^{1},\omega _{6}^{1},\omega _{7}^{1}\}$ by the
simple roots of $E_{8}$:

\begin{quote}
${\small h}_{\ast }{\small (\omega }_{1}^{1}{\small )=}\frac{1}{2}{\small %
\alpha }_{2}{\small +}\frac{1}{2}{\small \alpha }_{3}{\small +\alpha }_{4}%
{\small +\alpha }_{5}{\small +\alpha }_{6}{\small +\alpha }_{7}{\small %
+\alpha }_{8}{\small (\equiv }\frac{1}{2}{\small \omega }_{1}\func{mod}%
{\small \Lambda }_{G}^{r}{\small )}$

${\small h}_{\ast }{\small (\omega }_{6}^{1}{\small )=}\frac{5}{4}{\small %
\alpha }_{2}{\small +}\frac{7}{4}{\small \alpha }_{3}{\small +}\frac{5}{2}%
{\small \alpha }_{4}{\small +2\alpha }_{5}{\small +}\frac{3}{2}{\small %
\alpha }_{6}{\small +\alpha }_{7}{\small +}\frac{1}{2}{\small \alpha }_{8}%
{\small (\equiv }\frac{1}{4}{\small \omega }_{1}\func{mod}{\small \Lambda }%
_{G}^{r}{\small )}$

${\small h}_{\ast }{\small (\omega }_{7}^{1}{\small )=}\frac{7}{4}{\small %
\alpha }_{2}{\small +}\frac{5}{4}{\small \alpha }_{3}{\small +}\frac{5}{2}%
{\small \alpha }_{4}{\small +2\alpha }_{5}{\small +}\frac{3}{2}{\small %
\alpha }_{6}{\small +\alpha }_{7}{\small +}\frac{1}{2}{\small \alpha }_{8}%
{\small (\equiv }\frac{3}{4}{\small \omega }_{1}\func{mod}{\small \Lambda }%
_{G}^{r}{\small )}$
\end{quote}

\noindent and that

\begin{quote}
$h_{\ast }(\Lambda _{C_{\exp (u)}^{Rad}}^{e}(\mathbb{Q}))=\{\lambda \omega
_{1}\mid \lambda \in \mathbb{Q}\}$.
\end{quote}

\noindent It follows that the set $H_{u}$ consists of three elements in
which the deficiency of the element $\omega _{6}^{1}$ in the group $Spin(14)$
is $4$.

Consequently, $\ker \pi =\mathbb{Z}_{4}$ with generator $\exp _{1}(\omega
_{6}^{1})\times \exp _{2}(-\frac{1}{4}\omega _{1})$ by (4.3).

ii) If $u=\frac{\omega _{2}}{6}$, the local type of the centralizer $C_{\exp
(u)}$ is $SU(8)\times S^{1}$. Accordingly, assume that the set of
fundamental dominant weights of the semisimple part $SU(8)$ is $\Omega
=\{\omega _{1}^{1},\omega _{2}^{1},\omega _{3}^{1},\omega _{4}^{1},\omega
_{5}^{1},\omega _{6}^{1},\omega _{7}^{1}\}$. Applying Lemma 4.1 we get the
expressions of $h_{\ast }(\omega )$ with $\omega \in \Pi _{SU(8)}$ by the
simple roots of $E_{8}$:

\begin{quote}
${\small h}_{\ast }{\small (\omega }_{1}^{1}{\small )=}\frac{7}{8}{\small %
\alpha }_{1}{\small +}\frac{3}{4}{\small \alpha }_{3}{\small +}\frac{5}{8}%
{\small \alpha }_{4}{\small +}\frac{1}{2}{\small \alpha }_{5}{\small +}\frac{%
3}{8}{\small \alpha }_{6}{\small +}\frac{1}{4}{\small \alpha }_{7}{\small +}%
\frac{1}{8}{\small \alpha }_{8}{\small (\equiv }\frac{3}{8}{\small \omega }%
_{2}\func{mod}{\small \Lambda }_{G}^{r}{\small )}$

${\small h}_{\ast }{\small (\omega }_{2}^{1}{\small )=}\frac{3}{4}{\small %
\alpha }_{1}{\small +}\frac{3}{2}{\small \alpha }_{3}{\small +}\frac{5}{4}%
{\small \alpha }_{4}{\small +\alpha }_{5}{\small +}\frac{3}{4}{\small \alpha 
}_{6}{\small +}\frac{1}{2}{\small \alpha }_{7}{\small +}\frac{1}{4}{\small %
\alpha }_{8}{\small (\equiv }\frac{3}{4}{\small \omega }_{2}\func{mod}%
{\small \Lambda }_{G}^{r}{\small )}$

${\small h}_{\ast }{\small (\omega }_{3}^{1}{\small )=}\frac{5}{8}{\small %
\alpha }_{1}{\small +}\frac{5}{4}{\small \alpha }_{3}{\small +}\frac{15}{8}%
{\small \alpha }_{4}{\small +}\frac{3}{2}{\small \alpha }_{5}{\small +}\frac{%
9}{8}{\small \alpha }_{6}{\small +}\frac{3}{4}{\small \alpha }_{7}{\small +}%
\frac{3}{8}{\small \alpha }_{8}{\small (\equiv }\frac{1}{8}{\small \omega }%
_{2}\func{mod}{\small \Lambda }_{G}^{r}{\small )}$

${\small h}_{\ast }{\small (\omega }_{4}^{1}{\small )=}\frac{1}{2}{\small %
\alpha }_{1}{\small +\alpha }_{3}{\small +}\frac{3}{2}{\small \alpha }_{4}%
{\small +2\alpha }_{5}{\small +}\frac{3}{2}{\small \alpha }_{6}{\small %
+\alpha }_{7}{\small +}\frac{1}{2}{\small \alpha }_{8}{\small (\equiv }\frac{%
1}{2}{\small \omega }_{2}\func{mod}{\small \Lambda }_{G}^{r}{\small )}$

${\small h}_{\ast }{\small (\omega }_{5}^{1}{\small )=}\frac{3}{8}{\small %
\alpha }_{1}{\small +}\frac{3}{4}{\small \alpha }_{3}{\small +}\frac{9}{8}%
{\small \alpha }_{4}{\small +}\frac{3}{2}{\small \alpha }_{5}{\small +}\frac{%
15}{8}{\small \alpha }_{6}{\small +}\frac{5}{4}{\small \alpha }_{7}{\small +}%
\frac{5}{8}{\small \alpha }_{8}{\small (\equiv }\frac{7}{8}{\small \omega }%
_{2}\func{mod}{\small \Lambda }_{G}^{r}{\small )}$

${\small h}_{\ast }{\small (\omega }_{6}^{1}{\small )=}\frac{1}{4}{\small %
\alpha }_{1}{\small +}\frac{1}{2}{\small \alpha }_{3}{\small +}\frac{3}{4}%
{\small \alpha }_{4}{\small +\alpha }_{5}{\small +}\frac{5}{4}{\small \alpha 
}_{6}{\small +}\frac{3}{2}{\small \alpha }_{7}{\small +}\frac{3}{4}{\small %
\alpha }_{8}{\small (\equiv }\frac{1}{4}{\small \omega }_{2}\func{mod}%
{\small \Lambda }_{G}^{r}{\small )}$

${\small h}_{\ast }{\small (\omega }_{7}^{1}{\small )=}\frac{1}{8}{\small %
\alpha }_{1}{\small +}\frac{1}{4}{\small \alpha }_{3}{\small +}\frac{3}{8}%
{\small \alpha }_{4}{\small +}\frac{1}{2}{\small \alpha }_{5}{\small +}\frac{%
5}{8}{\small \alpha }_{6}{\small +}\frac{3}{4}{\small \alpha }_{7}{\small +}%
\frac{7}{8}{\small \alpha }_{8}{\small (\equiv }\frac{5}{8}{\small \omega }%
_{2}\func{mod}{\small \Lambda }_{G}^{r}{\small )}$
\end{quote}

\noindent and get

\begin{quote}
$h_{\ast }(\Lambda _{C_{\exp (u)}^{Rad}}^{e}(\mathbb{Q}))=\{\lambda \omega
_{2}\mid \lambda \in \mathbb{Q}\}$.
\end{quote}

\noindent It follows that the set $H_{u}$ consists of $7$ elements, in which
the deficiency of $\omega _{1}^{1}$ in the group $SU(8)$ is $8$.

Consequently, $\ker \pi =\mathbb{Z}_{8}$ with generator $\exp _{1}(\omega
_{1}^{1})\times \exp _{2}(-\frac{3}{8}\omega _{2})$ by (4.3).

iii) If $u=\frac{\omega _{3}}{8}$, the local type of the centralizer $%
C_{\exp (u)}$ is $SU(2)\times SU(7)\times S^{1}$. Accordingly, assume that
the set of fundamental dominant weights of the semisimple part $SU(2)\times
SU(7)$ is $\Omega =\{\omega _{1}^{1}\}\amalg \{\omega _{1}^{2},\omega
_{2}^{2},\omega _{3}^{2},\omega _{4}^{2},\omega _{5}^{2},\omega _{6}^{2}\}$.
Applying Lemma 4.1 we get the expressions of $h_{\ast }(\omega )$ with $%
\omega \in \Pi _{SU(2)}\sqcup \Pi _{SU(7)}$ by the simple roots of $E_{8}$:

\begin{quote}
${\small h}_{\ast }{\small (\omega }_{1}^{1}{\small )=}\frac{1}{2}{\small %
\alpha }_{1}{\small (\equiv }\frac{1}{2}{\small \omega }_{3}\func{mod}%
{\small \Lambda }_{G}^{r}{\small )}$

${\small h}_{\ast }{\small (\omega }_{1}^{2}{\small )=}\frac{6}{7}{\small %
\alpha }_{2}{\small +}\frac{5}{7}{\small \alpha }_{4}{\small +}\frac{4}{7}%
{\small \alpha }_{5}{\small +}\frac{3}{7}{\small \alpha }_{6}{\small +}\frac{%
2}{7}{\small \alpha }_{7}{\small +}\frac{1}{7}{\small \alpha }_{8}{\small %
(\equiv }\frac{2}{7}{\small \omega }_{3}\func{mod}{\small \Lambda }_{G}^{r}%
{\small )}$

${\small h}_{\ast }{\small (\omega }_{2}^{2}{\small )=}\frac{5}{7}{\small %
\alpha }_{2}{\small +}\frac{10}{7}{\small \alpha }_{4}{\small +}\frac{8}{7}%
{\small \alpha }_{5}{\small +}\frac{6}{7}{\small \alpha }_{6}{\small +}\frac{%
4}{7}{\small \alpha }_{7}{\small +}\frac{2}{7}{\small \alpha }_{8}{\small %
(\equiv }\frac{4}{7}{\small \omega }_{3}\func{mod}{\small \Lambda }_{G}^{r}%
{\small )}$

${\small h}_{\ast }{\small (\omega }_{3}^{2}{\small )=}\frac{4}{7}{\small %
\alpha }_{2}{\small +}\frac{8}{7}{\small \alpha }_{4}{\small +}\frac{12}{7}%
{\small \alpha }_{5}{\small +}\frac{9}{7}{\small \alpha }_{6}{\small +}\frac{%
6}{7}{\small \alpha }_{7}{\small +}\frac{3}{7}{\small \alpha }_{8}{\small %
(\equiv }\frac{6}{7}{\small \omega }_{3}\func{mod}{\small \Lambda }_{G}^{r}%
{\small )}$

${\small h}_{\ast }{\small (\omega }_{4}^{2}{\small )=}\frac{3}{7}{\small %
\alpha }_{2}{\small +}\frac{6}{7}{\small \alpha }_{4}{\small +}\frac{9}{7}%
{\small \alpha }_{5}{\small +}\frac{12}{7}{\small \alpha }_{6}{\small +}%
\frac{8}{7}{\small \alpha }_{7}{\small +}\frac{4}{7}{\small \alpha }_{8}%
{\small (\equiv }\frac{1}{7}{\small \omega }_{3}\func{mod}{\small \Lambda }%
_{G}^{r}{\small )}$

${\small h}_{\ast }{\small (\omega }_{5}^{2}{\small )=}\frac{2}{7}{\small %
\alpha }_{2}{\small +}\frac{4}{7}{\small \alpha }_{4}{\small +}\frac{6}{7}%
{\small \alpha }_{5}{\small +}\frac{8}{7}{\small \alpha }_{6}{\small +}\frac{%
10}{7}{\small \alpha }_{7}{\small +}\frac{5}{7}{\small \alpha }_{8}{\small %
(\equiv }\frac{3}{7}{\small \omega }_{3}\func{mod}{\small \Lambda }_{G}^{r}%
{\small )}$

${\small h}_{\ast }{\small (\omega }_{6}^{2}{\small )=}\frac{1}{7}{\small %
\alpha }_{2}{\small +}\frac{2}{7}{\small \alpha }_{4}{\small +}\frac{3}{7}%
{\small \alpha }_{5}{\small +}\frac{4}{7}{\small \alpha }_{6}{\small +}\frac{%
5}{7}{\small \alpha }_{7}{\small +}\frac{6}{7}{\small \alpha }_{8}{\small %
(\equiv }\frac{5}{7}{\small \omega }_{3}\func{mod}{\small \Lambda }_{G}^{r}%
{\small )}$
\end{quote}

\noindent and

\begin{quote}
$h_{\ast }(\Lambda _{C_{\exp (u)}^{Rad}}^{e}(\mathbb{Q}))=\{\lambda \omega
_{3}\mid \lambda \in \mathbb{Q}\}$.
\end{quote}

\noindent It follows that the set $H_{u}$ consists of $13$ elements, in
which the deficiency of $\omega _{1}^{1}\oplus \omega _{1}^{2}$ in the group 
$SU(2)\times SU(7)$ is $14$.

Consequently, $\ker \pi =\mathbb{Z}_{14}$ with generator $\exp _{1}(\omega
_{1}^{1})\times \exp _{2}(\omega _{1}^{2})\times \exp _{3}(-\frac{11}{14}%
\omega _{3})$ by (4.3)$.$

iv) If $u=\frac{\omega _{4}}{12}$, the local type of the centralizer $%
C_{\exp (u)}$ is $SU(2)\times SU(3)\times SU(5)\times S^{1}$. Accordingly,
assume that the set of fundamental dominant weights of the semisimple part $%
SU(2)\times SU(3)\times SU(5)$ is $\Omega =\{\omega _{1}^{1}\}\amalg
\{\omega _{1}^{2},\omega _{2}^{2}\}\amalg \{\omega _{1}^{3},\omega
_{2}^{3},\omega _{3}^{3},\omega _{4}^{3}\}$. Applying Lemma 4.1 we get the
expressions of $h_{\ast }(\omega )$ with $\omega \in \Pi _{SU(2)}\sqcup \Pi
_{SU(3)}\sqcup \Pi _{SU(5)}$ by the simple roots of $E_{8}$:

\begin{quote}
${\small h}_{\ast }{\small (\omega }_{1}^{1}{\small )=}\frac{1}{2}{\small %
\alpha }_{2}{\small (\equiv }\frac{1}{2}{\small \omega }_{4}\func{mod}%
{\small \Lambda }_{G}^{r}{\small )}$

${\small h}_{\ast }{\small (\omega }_{1}^{2}{\small )=}\frac{2}{3}{\small %
\alpha }_{1}{\small +}\frac{1}{3}{\small \alpha }_{3}{\small (\equiv }\frac{2%
}{3}{\small \omega }_{4}\func{mod}{\small \Lambda }_{G}^{r}{\small )}$

${\small h}_{\ast }{\small (\omega }_{2}^{2}{\small )=}\frac{1}{3}{\small %
\alpha }_{1}{\small +}\frac{2}{3}{\small \alpha }_{3}{\small (\equiv }\frac{1%
}{3}{\small \omega }_{4}\func{mod}{\small \Lambda }_{G}^{r}{\small )}$

${\small h}_{\ast }{\small (\omega }_{1}^{3}{\small )=}\frac{4}{5}{\small %
\alpha }_{5}{\small +}\frac{3}{5}{\small \alpha }_{6}{\small +}\frac{2}{5}%
{\small \alpha }_{7}{\small +}\frac{1}{5}{\small \alpha }_{8}{\small (\equiv 
}\frac{1}{5}{\small \omega }_{4}\func{mod}{\small \Lambda }_{G}^{r}{\small )}
$

${\small h}_{\ast }{\small (\omega }_{2}^{3}{\small )=}\frac{3}{5}{\small %
\alpha }_{5}{\small +}\frac{6}{5}{\small \alpha }_{6}{\small +}\frac{4}{5}%
{\small \alpha }_{7}{\small +}\frac{2}{5}{\small \alpha }_{8}{\small (\equiv 
}\frac{2}{5}{\small \omega }_{4}\func{mod}{\small \Lambda }_{G}^{r}{\small )}
$

${\small h}_{\ast }{\small (\omega }_{3}^{3}{\small )=}\frac{2}{5}{\small %
\alpha }_{5}{\small +}\frac{4}{5}{\small \alpha }_{6}{\small +}\frac{6}{5}%
{\small \alpha }_{7}{\small +}\frac{3}{5}{\small \alpha }_{8}{\small (\equiv 
}\frac{3}{5}{\small \omega }_{4}\func{mod}{\small \Lambda }_{G}^{r}{\small )}
$

${\small h}_{\ast }{\small (\omega }_{4}^{3}{\small )=}\frac{1}{5}{\small %
\alpha }_{5}{\small +}\frac{2}{5}{\small \alpha }_{6}{\small +}\frac{3}{5}%
{\small \alpha }_{7}{\small +}\frac{4}{5}{\small \alpha }_{8}{\small (\equiv 
}\frac{4}{5}{\small \omega }_{4}\func{mod}{\small \Lambda }_{G}^{r}{\small %
)\allowbreak }$
\end{quote}

\noindent and get

\begin{quote}
$h_{\ast }(\Lambda _{C_{\exp (u)}^{Rad}}^{e}(\mathbb{Q}))=\{\lambda \omega
_{4}\mid \lambda \in \mathbb{Q}\}$.
\end{quote}

\noindent It follows that the set $H_{u}$ consists of $29$ elements, in
which the deficiency of $\omega _{1}^{1}\oplus \omega _{1}^{2}\oplus \omega
_{1}^{3}$ in the group $SU(2)\times SU(3)\times SU(5)$ is $30$.

Consequently, $\ker \pi =\mathbb{Z}_{30}$ with generator $\exp _{1}(\omega
_{1}^{1})\times \exp _{2}(\omega _{1}^{2})\times \exp _{3}(\omega
_{1}^{3})\times \exp _{4}(-\frac{11}{30}\omega _{4})$ by (4.3).

v) If $u=\frac{\omega _{5}}{10}$, the local type of the centralizer $C_{\exp
(u)}$ is $SU(5)\times SU(4)\times S^{1}$. Accordingly, assume that the set
of fundamental dominant weights of the semisimple part $SU(5)\times SU(4)$
is $\Omega =\{\omega _{1}^{1},\omega _{2}^{1},\omega _{3}^{1},\omega
_{4}^{1}\}\amalg \{\omega _{1}^{2},\omega _{2}^{2},\omega _{3}^{2}\}$.
Applying Lemma 4.1 we get the expressions of $h_{\ast }(\omega )$ with $%
\omega \in \Pi _{SU(5)}\sqcup \Pi _{SU(4)}$ by the simple roots of $E_{8}$:

\begin{quote}
${\small h}_{\ast }{\small (\omega }_{1}^{1}{\small )=}\frac{4}{5}{\small %
\alpha }_{1}{\small +}\frac{1}{5}{\small \alpha }_{2}{\small +}\frac{3}{5}%
{\small \alpha }_{3}{\small +}\frac{2}{5}{\small \alpha }_{4}{\small (\equiv 
}\frac{3}{5}{\small \omega }_{5}\func{mod}{\small \Lambda }_{G}^{r}{\small )}
$

${\small h}_{\ast }{\small (\omega }_{2}^{1}{\small )=}\frac{3}{5}{\small %
\alpha }_{1}{\small +}\frac{2}{5}{\small \alpha }_{2}{\small +}\frac{6}{5}%
{\small \alpha }_{3}{\small +}\frac{4}{5}{\small \alpha }_{4}{\small (\equiv 
}\frac{1}{5}{\small \omega }_{5}\func{mod}{\small \Lambda }_{G}^{r}{\small )}
$

${\small h}_{\ast }{\small (\omega }_{3}^{1}{\small )=}\frac{2}{5}{\small %
\alpha }_{1}{\small +}\frac{3}{5}{\small \alpha }_{2}{\small +}\frac{4}{5}%
{\small \alpha }_{3}{\small +}\frac{6}{5}{\small \alpha }_{4}{\small (\equiv 
}\frac{4}{5}{\small \omega }_{5}\func{mod}{\small \Lambda }_{G}^{r}{\small )}
$

${\small h}_{\ast }{\small (\omega }_{4}^{1}{\small )=}\frac{1}{5}{\small %
\alpha }_{1}{\small +}\frac{4}{5}{\small \alpha }_{2}{\small +}\frac{2}{5}%
{\small \alpha }_{3}{\small +}\frac{3}{5}{\small \alpha }_{4}{\small (\equiv 
}\frac{2}{5}{\small \omega }_{5}\func{mod}{\small \Lambda }_{G}^{r}{\small )}
$

${\small h}_{\ast }{\small (\omega }_{1}^{2}{\small )=}\frac{3}{4}{\small %
\alpha }_{6}{\small +}\frac{1}{2}{\small \alpha }_{7}{\small +}\frac{1}{4}%
{\small \alpha }_{8}{\small (\equiv }\frac{1}{4}{\small \omega }_{5}\func{mod%
}{\small \Lambda }_{G}^{r}{\small )}$

${\small h}_{\ast }{\small (\omega }_{2}^{2}{\small )=}\frac{1}{2}{\small %
\alpha }_{6}{\small +\alpha }_{7}{\small +}\frac{1}{2}{\small \alpha }_{8}%
{\small (\equiv }\frac{1}{2}{\small \omega }_{5}\func{mod}{\small \Lambda }%
_{G}^{r}{\small )}$

${\small h}_{\ast }{\small (\omega }_{3}^{2}{\small )=}\frac{1}{4}{\small %
\alpha }_{6}{\small +}\frac{1}{2}{\small \alpha }_{7}{\small +}\frac{3}{4}%
{\small \alpha }_{8}{\small (\equiv }\frac{3}{4}{\small \omega }_{5}\func{mod%
}{\small \Lambda }_{G}^{r}{\small )}$
\end{quote}

\noindent and get

\begin{quote}
$h_{\ast }(\Lambda _{C_{\exp (u)}^{Rad}}^{e}(\mathbb{Q}))=\{\lambda \omega
_{5}\mid \lambda \in \mathbb{Q}\}$.
\end{quote}

\noindent It follows that the set $H_{u}$ consists of $19$ elements in which
the deficiency of $\omega _{1}^{1}\oplus \omega _{1}^{2}$ in the group $%
SU(5)\times SU(4)$ is $20$.

Consequently, $\ker \pi =\mathbb{Z}_{20}$ with generator $\exp _{1}(\omega
_{1}^{1})\times \exp _{2}(\omega _{1}^{2})\times \exp _{3}(-\frac{17}{20}%
\omega _{5})$ by (4.3).

vi) If $u=\frac{\omega _{6}}{8}$, the local type of the centralizer $C_{\exp
(u)}$ is $Spin(10)\times SU(3)\times S^{1}$. Accordingly, assume that the
set of fundamental dominant weights of the semisimple part $Spin(10)\times
SU(3)$ is $\Omega =\{\omega _{1}^{1},\omega _{2}^{1},\omega _{3}^{1},\omega
_{4}^{1},\omega _{5}^{1}\}\amalg \{\omega _{1}^{2},\omega _{2}^{2}\}$.
Applying Lemma 4.1 we get the expressions of $h_{\ast }(\omega )$ with $%
\omega \in \Pi _{Spin(10)}\sqcup \Pi _{SU(3)}$ by the simple roots of $E_{8}$%
:

\begin{quote}
${\small h}_{\ast }{\small (\omega }_{1}^{1}{\small )=\alpha }_{1}{\small +}%
\frac{1}{2}{\small \alpha }_{2}{\small +\alpha }_{3}{\small +\alpha }_{4}%
{\small +}\frac{1}{2}{\small \alpha }_{5}{\small (\equiv }\frac{1}{2}{\small %
\omega }_{6}\func{mod}{\small \Lambda }_{G}^{r}{\small )}$

${\small h}_{\ast }{\small (\omega }_{4}^{1}{\small )=}\frac{1}{2}{\small %
\alpha }_{1}{\small +}\frac{5}{4}{\small \alpha }_{2}{\small +\alpha }_{3}%
{\small +}\frac{3}{2}{\small \alpha }_{4}{\small +}\frac{3}{4}{\small \alpha 
}_{5}{\small (\equiv }\frac{1}{4}{\small \omega }_{6}\func{mod}{\small %
\Lambda }_{G}^{r}{\small )}$

${\small h}_{\ast }{\small (\omega }_{5}^{1}{\small )=}\frac{1}{2}{\small %
\alpha }_{1}{\small +}\frac{3}{4}{\small \alpha }_{2}{\small +\alpha }_{3}%
{\small +}\frac{3}{2}{\small \alpha }_{4}{\small +}\frac{5}{4}{\small \alpha 
}_{5}{\small (\equiv }\frac{3}{4}{\small \omega }_{6}\func{mod}{\small %
\Lambda }_{G}^{r}{\small )}$

${\small h}_{\ast }{\small (\omega }_{1}^{2}{\small )=}\frac{2}{3}{\small %
\alpha }_{7}{\small +}\frac{1}{3}{\small \alpha }_{8}{\small (\equiv }\frac{1%
}{3}{\small \omega }_{6}\func{mod}{\small \Lambda }_{G}^{r}{\small )}$

${\small h}_{\ast }{\small (\omega }_{2}^{2}{\small )=}\frac{1}{3}{\small %
\alpha }_{7}{\small +}\frac{2}{3}{\small \alpha }_{8}{\small (\equiv }\frac{2%
}{3}{\small \omega }_{6}\func{mod}{\small \Lambda }_{G}^{r}{\small )}$
\end{quote}

\noindent and get

\begin{quote}
$h_{\ast }(\Lambda _{C_{\exp (u)}^{Rad}}^{e}(\mathbb{Q}))=\{\lambda \omega
_{6}\mid \lambda \in \mathbb{Q}\}$.
\end{quote}

\noindent It follows that the set $H_{u}$ consists of $11$ elements in which
the deficiency of $\omega _{4}^{1}\oplus \omega _{1}^{2}$ in the group $%
Spin(10)\times SU(3)$ is $12$.

Consequently, $\ker \pi =\mathbb{Z}_{12}$ with generator $\exp _{1}(\omega
_{4}^{1})\times \exp _{2}(\omega _{1}^{2})\times \exp _{3}(-\frac{7}{12}%
\omega _{6})$ by (4.3).

vii) If $u=\frac{\omega _{7}}{6}$, the local type of the centralizer $%
C_{\exp (u)}$ is $E_{6}\times SU(2)\times S^{1}$. Accordingly, assume that
the set of fundamental dominant weights of the semisimple part $E_{6}\times
SU(2)$ is $\Omega =\{\omega _{1}^{1},\omega _{2}^{1},\omega _{3}^{1},\omega
_{4}^{1},\omega _{5}^{1},\omega _{6}^{1}\}\amalg \{\omega _{1}^{2}\}$.
Applying Lemma 4.1 we get the expressions of $h_{\ast }(\omega )$ with $%
\omega \in \Pi _{E_{6}}\sqcup \Pi _{SU(2)}$ by the simple roots of $E_{8}$:

\begin{quote}
${\small h}_{\ast }{\small (\omega }_{1}^{1}{\small )=}\frac{4}{3}{\small %
\alpha }_{1}{\small +\alpha }_{2}{\small +}\frac{5}{3}{\small \alpha }_{3}%
{\small +2\alpha }_{4}{\small +}\frac{4}{3}{\small \alpha }_{5}{\small +}%
\frac{2}{3}{\small \alpha }_{6}{\small (\equiv }\frac{1}{3}{\small \omega }%
_{7}\func{mod}{\small \Lambda }_{G}^{r}{\small )}$

${\small h}_{\ast }{\small (\omega }_{6}^{1}{\small )=}\frac{2}{3}{\small %
\alpha }_{1}{\small +\alpha }_{2}{\small +}\frac{4}{3}{\small \alpha }_{3}%
{\small +2\alpha }_{4}{\small +}\frac{5}{3}{\small \alpha }_{5}{\small +}%
\frac{4}{3}{\small \alpha }_{6}{\small (\equiv }\frac{2}{3}{\small \omega }%
_{7}\func{mod}{\small \Lambda }_{G}^{r}{\small )}$

${\small h}_{\ast }{\small (\omega }_{1}^{2}{\small )=}\frac{1}{2}{\small %
\alpha }_{8}{\small (\equiv }\frac{1}{2}{\small \omega }_{7}\func{mod}%
{\small \Lambda }_{G}^{r}{\small )}$
\end{quote}

\noindent and get

\begin{quote}
$h_{\ast }(\Lambda _{C_{\exp (u)}^{Rad}}^{e}(\mathbb{Q}))=\{\lambda \omega
_{7}\mid \lambda \in \mathbb{Q}\}$.
\end{quote}

\noindent It follows that the set $H_{u}$ consists of $5$ elements, in which
the deficiency of $\omega _{1}^{1}\oplus \omega _{1}^{2}$ in the group $%
E_{6}\times SU(2)$ is $6$.

Consequently, $\ker \pi =\mathbb{Z}_{6}$ with generator $\exp _{1}(\omega
_{1}^{1})\times \exp _{2}(\omega _{1}^{2})\times \exp _{3}(-\frac{5}{6}%
\omega _{7})$ by (4.3).

viii) If $u=\frac{\omega _{8}}{4}$, the local type of the centralizer $%
C_{\exp (u)}$ is $E_{7}\times S^{1}$ by Theorem 2.8. Accordingly, assume
that the set of fundamental dominant weights of the semisimple part $E_{7}$
is $\Omega =\{\omega _{1}^{1},\omega _{2}^{1},\omega _{3}^{1},\omega
_{4}^{1},\omega _{5}^{1},\omega _{6}^{1},\omega _{7}^{1}\}$. Applying Lemma
4.1 we get the expressions of $h_{\ast }(\omega )$ with $\omega \in \Pi
_{E_{7}}=\{\omega _{7}^{1}\}$ by the simple roots of $E_{8}$:

${\small h}_{\ast }{\small (\omega }_{7}^{1}{\small )=\alpha }_{1}{\small +}%
\frac{3}{2}{\small \alpha }_{2}{\small +2\alpha }_{3}{\small +3\alpha }_{4}%
{\small +}\frac{5}{2}{\small \alpha }_{5}{\small +2\alpha }_{6}{\small +}%
\frac{3}{2}{\small \alpha }_{7}{\small (\equiv }\frac{1}{2}{\small \omega }%
_{8}\func{mod}{\small \Lambda }_{G}^{r}{\small )\allowbreak }$

\noindent and get

\begin{quote}
$h_{\ast }(\Lambda _{C_{\exp (u)}^{Rad}}^{e}(\mathbb{Q}))=\{\lambda \omega
_{8}\mid \lambda \in \mathbb{Q}\}$.
\end{quote}

\noindent It follows that the set $H_{u}$ consists of the single element $%
\omega _{7}^{1}$ whose deficiency in the group $E_{7}$ is $2$.

Consequently, $\ker \pi =\mathbb{Z}_{2}$ with generator $\exp _{1}(\omega
_{7}^{1})\times \exp _{2}(-\frac{1}{2}\omega _{8})$.$\square $

\end{document}